\documentclass[11pt]{article}
\usepackage{amsfonts}
\usepackage{amsfonts}
\usepackage{amsfonts}
\usepackage{amsmath}
\usepackage{mathrsfs}
\usepackage{amsmath}
\usepackage{mathrsfs}
\usepackage{mathrsfs}
\usepackage{amsmath}
\usepackage{fancyhdr}
\usepackage{mathrsfs,epsfig,morefloats}
\usepackage{amsfonts}
\usepackage{natbib}
\usepackage{amsmath}
\usepackage{amsthm}
\usepackage{amssymb}
\usepackage{color}
\usepackage{lscape}
\usepackage{rotating}
\usepackage{caption2}
\usepackage{subfigure}
\usepackage{graphicx}
\usepackage{amssymb}
\usepackage[compact]{titlesec}

\newcommand{\Date}[1]{\def\@Date{#1}}
\def\today{\number\day~\ifcase\month\or
 January\or February\or March\or April\or May\or June\or
 July\or August\or September\or October\or November\or December\fi~\number\year}

\def\be{\begin{equation}}
\def\ee{\end{equation}}
\def\bea{\begin{eqnarray}}
\def\eea{\end{eqnarray}}
\def\bd{\begin{displaymath}}
\def\ed{\end{displaymath}}
\def\bda{\begin{eqnarray*}}
\def\eda{\end{eqnarray*}}
\def\bsm{\begin{small}}
\def\esm{\end{small}}

\def\t0{\theta_0}

\def\ha1{\hat \beta_1}

\def\bnt{\begin{enumerate}}
\def\ent{\end{enumerate}}
\def\T{{ \mathrm{\scriptscriptstyle T} }}

\def\AS{A\"{\i}t-Sahalia}

\def\bsc{\begin{scriptsize}}
\def\esc{\end{scriptsize}}

\newtheorem{theorem}{Theorem}

\newtheorem{lemma}{Lemma}

\newtheorem{proposition}{Proposition}

\theoremstyle{definition}
\newtheorem{as}{Condition}

\makeatletter
\newcommand{\figcaption}{\def\@captype{figure}\caption}
\newcommand{\tabcaption}{\def\@captype{table}\caption}
\makeatother

\newcommand{\sgn}{\mbox{\rm sgn}}

\newcommand{\supp}{\mathrm{supp}}

\newcommand{\bA}{{\mathbf A}}
\newcommand{\bB}{{\mathbf B}}

\newcommand{\bJ}{{\mathbf J}}

\newcommand{\bI}{{\mathbf I}}

\newcommand{\bM}{{\mathbf M}}
\newcommand{\bQ}{{\mathbf Q}}

\newcommand{\bR}{{\mathbf R}}

\newcommand{\bU}{{\mathbf U}}
\newcommand{\bV}{{\mathbf V}}

\newcommand{\bX}{{\mathbf X}}
\newcommand{\bY}{{\mathbf Y}}
\newcommand{\bZ}{{\mathbf Z}}
\newcommand{\ba}{{\mathbf a}}

\newcommand{\bfe}{{\mathbf e}}

\newcommand{\bh}{{\mathbf h}}
\newcommand{\bg}{{\mathbf g}}

\newcommand{\bm}{{\mathbf m}}

\newcommand{\bv}{{\mathbf v}}
\newcommand{\bw}{{\mathbf w}}

\newcommand{\bz}{{\mathbf z}}
\newcommand{\balpha} {\boldsymbol{\alpha}}

\newcommand{\bfeta}  {\boldsymbol{\eta}}

\newcommand{\bdelta} {\boldsymbol{\delta}}
\newcommand{\blambda}{\boldsymbol{\lambda}}

\newcommand{\bSigma}{\boldsymbol{\Sigma}}

\newcommand{\bTheta} {\boldsymbol{\Theta}}

\newcommand{\bpsi}{\boldsymbol{\psi}}

\newcommand{\btheta} {\boldsymbol{\theta}}

\newcommand{\bzeta} {\boldsymbol{\zeta}}

\newcommand{\bkappa}{\boldsymbol{\kappa}}

\newcommand{\bzero}{{\mathbf 0}}

\def\JRSSB{{\sl Journal of the Royal Statistical Society}, {\bf B}}

\def\BKA{{\sl Biometrika}}
\def\JASA{{\sl Journal of the American Statistical Association}}

\def\AS{{\sl The Annals of Statistics}}

\def\AP{{\sl The Annals of Probability}}

\def\JOE{{\sl Journal of Econometrics}}

\def\JE{{\sl Journal of Econometrics}}

\def\ECA{{\sl Econometrica}}

\begin{document}

\title{ \bf   A New Scope of Penalized Empirical Likelihood with High-Dimensional Estimating Equations}

\author{Jinyuan Chang\\Southwestern University of \\ Finance and Economics \and Cheng Yong Tang\\Temple University \and Tong Tong Wu\\University of Rochester}

\date{}
\maketitle

\begin{abstract}

Statistical methods with empirical likelihood (EL) are appealing and effective  
especially in conjunction with estimating equations through which useful data information can be adaptively and flexibly incorporated. It is also known in the literature that EL approaches encounter difficulties when dealing with problems having high-dimensional model parameters and estimating equations. To overcome the challenges,
we  begin our study with a careful investigation on high-dimensional EL  from a new scope targeting at estimating a high-dimensional sparse model parameters.
We show that the new scope provides an opportunity for relaxing the stringent requirement on the dimensionality of the model parameter.
Motivated by the new scope, we then propose  a  new penalized EL by applying two penalty functions respectively regularizing the model parameters and the associated Lagrange multipliers in the optimizations of EL. By penalizing the Lagrange multiplier to encourage its sparsity, we show that drastic dimension reduction in the number of estimating equations can be effectively achieved without compromising the validity and consistency of the resulting estimators. Most attractively, such a reduction in dimensionality of estimating equations is actually equivalent to a selection among those high-dimensional estimating equations, resulting in a highly parsimonious and effective device for high-dimensional sparse model parameters. Allowing both the dimensionalities of model parameters and estimating equations growing exponentially with the sample size, our theory demonstrates that the estimator from our new penalized EL is sparse and consistent with asymptotically normally distributed nonzero components. Numerical simulations and a real data analysis show that the proposed penalized EL works promisingly.

\end{abstract}

\begin{quote}
\noindent
{\sl Keywords}: Empirical likelihood; Estimating equations; High-dimensional statistical methods; Moment selection; Penalized likelihood. \end{quote}

\begin{quote}
\noindent
{\sl MSC2010 subject classifications}: Primary 62G99; secondary 62F40
\end{quote}

\thispagestyle{empty}
\pagenumbering{gobble}

\newpage
\pagenumbering{arabic}

\setcounter{page}{1}

\section{Introduction}
\label{s1}

Statistical approaches using estimating equations are widely applicable to solve a broad class of practical problems. The most influential special cases of estimating equations include the fundamental maximum likelihood score equations and those from the popular generalized methods of moments \citep{Hansen1982}. The approaches of using estimating equations are particularly appealing in practice with merits from requiring less stringent distributional assumptions on the data model, yet being adaptable to flexibly incorporate suitable information and conditions extracted from practical features in various scenarios of interests.

Empirical likelihood (EL, hereinafter) \citep{Owen2001} coupled with estimating equations has been demonstrated successful since the seminal work of \cite{QinLawless1994}. It is particularly appealing that the maximum EL estimator asymptotically achieves the semiparametric efficiency bound \citep{QinLawless1994}. The properties of EL are also desirable through some higher order analyses \citep{NeweySmith2004,Chen:Cui:2006,ChenCui2007}. Moreover, the Wilks' theorems \citep{Owen2001,QinLawless1994} for EL ensure that EL ratio is asymptotically central chi-square distributed when evaluated at the truth.  Hence, EL provides an analogous device to the conventional fully parametric likelihood for statistical inferences, but without requiring a  fully parametric likelihood built upon more stringent distributional assumptions.

In recent years, high data dimensionality in practice has attracted increasing research attention and brought unprecedented challenges to approaches based on estimating equations and EL. On one hand, studies in \cite{ChenPengQin2009}, \cite{Hjortetal2009}, \cite{TangLeng2010}, \cite{LengTang2012}, and \cite{ChangChenChen_2015} reveal that  conventional asymptotic schemes and results for EL are expected to work only when both the dimensionality of the parameter $p$ and the number of the estimating equations $r$ are growing at some rate slower than the sample size $n$.  On the other hand, however, challenges due to high-dimensionality require a capacity to deal with cases where both $p$ and $r$ can be much larger than $n$. \cite{TangLeng2010}, \cite{LengTang2012}, and \cite{ChangChenChen_2015} attempt to utilize sparsity of the model parameters by applying penalty functions on those parameters. Their results show that sparse estimators with good properties are achievable. However, the restriction from the data dimensionality is not alleviated by using penalized EL in their works.

The challenges for EL  from high data  dimensionality  are well documented in the literature, and there are recent investigations on the remedies.
\cite{Tsao_2004_AOS} found that for fixed $n$  with moderately large fixed $p$,  the probability that the truth is contained in the EL based confidence region can be substantially smaller than the nominal level, resulting in the under-coverage problem.   As remedies, \cite{TsaoWu_2013_AOS, TsaoWu_2014_Bioka} propose  extended EL to address the under-coverage problems due to the constraints on the parameter space.  With a modification avoiding  equality constraints, \cite{Bart_2007_SPL} propose a penalized EL method via optimizing products of probability weights penalized by a loss function depending on the model parameter.  \cite{LahiriMuk_2012_AOS} propose a different type of loss from that in \cite{Bart_2007_SPL} and study its properties with high-dimensional model parameter and dependent data. To our best knowledge, no estimation problems have been investigated with the EL formulations of \cite{Bart_2007_SPL} and \cite{LahiriMuk_2012_AOS}.

In this paper,  from a new scope on investigating high-dimensional sparse model parameters, we  study the properties of EL by carefully examining the impacts from the data dimensionally, and exploring the opportunity from targeting at the sparse model parameter.  We find that consistently estimating high-dimensional sparse model parameter by a penalized EL is feasible with fewer number of estimating functions than the model parameter.
Such an observation motivates us to  propose a new penalized EL approach to tackle high-dimensional statistical problems where both the numbers of model parameters and estimating equations, $p$ and $r$ respectively, can grow at an exponential rate of the sample size $n$. We solve the problem by employing two penalty functions when constructing the EL with high-dimensional estimating equations. Specifically, the first penalty function is on the magnitude of the model parameters with the goal to encourage sparsity in the resulting estimator. Additionally, a second penalty function is imposed on the Lagrange multiplier to encourage its sparsity when optimizing the EL evaluated at given values of the parameters. We also observe that obtaining a sparse Lagrange multiplier in EL is equivalent to reducing the dimensionality $r$ via an effective selection among those estimating equations, which itself is an interesting problem and a new scope; see our discussions in Sections \ref{s2} and \ref{s3}.

Here we note that the effect of the sparsity encouraging penalty on the Lagrange multiplier relates to the methods for selecting moments in the GMM methods, a problem that has been extensively studied in the econometrics literature; see, among others, \cite{ChengLiao_2015_JOE} and reference therein.
Recently,  \cite{ChengLiao_2015_JOE}  and \cite{Shi_2016} study the problem with many moment conditions for estimating a fixed dimensional model parameter.     \cite{ChengLiao_2015_JOE} propose to treat the  sample averages of the moment conditions as additional parameters to be optimized, and to apply the $L_1$ penalty on them to encourage sparsity so that effective  moment selection can be achieved.    The role of the $L_1$ penalty in their approach is seen similar to ours on the Lagrange multiplier  for the purpose of moment selection.   In light of the Dantzig selector approach of \cite{CandesTao_2007},    \cite{Shi_2016} propose a new EL formulation by relaxing the equality constraints to inequality ones involving some regularization parameter, so that effective selection of the moment conditions is also achieved.
Nevertheless,  none  of \cite{ChengLiao_2015_JOE} and \cite{Shi_2016} investigates the impacts from diverging number of model parameters that potentially can be sparse.

Our investigation contributes to the area of EL with high-dimensional statistical problems from a new scope. Our approach successfully extends the EL approach with estimating functions to scenarios  allowing both $p$ and $r$ growing exponentially with the sample size $n$. As shown in Sections \ref{s2} and \ref{s3}, new results for high-dimensional penalized EL are established, and many of them are interesting in both areas of EL and estimating equations. Our analysis first reveals a result of its own interests that substantially broadens  the understanding of the relationship between the number of estimating equations $r$ and the number of model parameters $p$ with penalized EL. Surprisingly, we find that  with an appropriate penalization, a consistent and sparse estimator of the model parameter actually does not require $r\ge p$, thanks to the new scope from estimating a sparse model parameter. In particular, we show that a sparse estimator with $s$  nonzero components for the $p$-dimensional parameter technically may only require that the number of estimating equations $r$ to be no less than $s$. Such a result crucially supports the motivation in our new penalized EL approach for the second penalty function imposed on the Lagrange multiplier to reduce the effective number of estimating equations actually involved in the high-dimensional penalized EL. That is, the resulting sparse Lagrange multiplier from the penalization is equivalent to a selection among available estimating equations for the model parameters. Our theory shows that the penalized EL estimator is consistent and can estimate the zero components of the model parameters as zero with  probability tending to one. Additionally, the nonzero components of the penalized EL estimator is asymptotically normally distributed.

The rest of this paper is organized as follows.  The new scope with high-dimensional sparse model parameter on EL and penalized EL is investigated in Section \ref{s2}. The new penalized EL with an additional penalty function on the Lagrange multiplier  and its properties  for estimating high-dimensional sparse model parameters  are given in Section \ref{s3}. An algorithm using coordinate descent for solving the penalized EL is presented in Section \ref{s4}. Numerical examples with simulated and real data are shown in Section \ref{s5}.  Some discussions are given in Section \ref{s6}. All technical details are provided in Section \ref{s7}. The Supplementary Material contains more technical proofs of the theoretical results.

\section{Empirical likelihood and penalized empirical likelihood}
\label{s2}

\subsection{An overview of empirical likelihood with diverging dimensionality}

Let us define some notations first. For a $q$-dimensional vector $\ba=(a_1,\ldots,a_q)^\T$, let $|\ba|_\infty=\max_{1\leq k\leq q}|a_k|$, $|\ba|_1=\sum_{k=1}^q|a_k|$ and $|\ba|_2=(\sum_{k=1}^qa_k^2)^{1/2}$ be its $L_\infty$-norm, $L_1$-norm, and $L_2$-norm, respectively. For a $q\times q$ matrix $\bM=(m_{ij})_{q\times q}$, let $\|\bM\|_\infty=\max_{1\leq i\leq q}\sum_{j=1}^q|m_{ij}|$,
$\|\bM\|_2=\lambda_{\max}^{1/2}(\bM^\T\bM)$ and $\|\bM\|_F=(\sum_{i,j=1}^qm_{ij}^2)^{1/2}$ be the $L_\infty$-norm, $L_2$-norm and Frobenius-norm of $\bM$, respectively.

Let $\bX_1,\ldots,\bX_n$ be $d$-dimensional independent and identically distributed generic observations and $\btheta=(\theta_1,\ldots,\theta_p)^\T$ be a $p$-dimensional parameter with support $\bTheta$. For an $r$-dimensional estimating function
    $
    \bg(\bX;\btheta)=\{g_1(\bX;\btheta),\ldots,g_r(\bX;\btheta)\}^\T,
    $
the information for the model parameter $\btheta$ is collected by the unbiased moment condition
    \begin{equation}\label{eq:esteq}
    \mathbb{E}\{\bg(\bX_i;\btheta_0)\}=\bzero,
    \end{equation}
where $\btheta_0\in\bTheta$ is the unknown truth. When the sample size $n$ grows,  following \cite{Hjortetal2009} and \cite{ChangChenChen_2015},  the observations $\{\bg(\bX_i;\btheta)\}_{i=1}^n$ can be viewed as a triangular array where $r$, $p$, $d$, $\bX_i$, $\btheta$ and $\bg(\bX;\btheta)$ may all depend on the sample size $n$.
Following the idea of EL \citep{Owen1988,Owen1990}, \cite{QinLawless1994} investigate an EL with estimating equations:
    \begin{equation}\label{eq:el}
    L(\btheta)=\sup\bigg\{\prod_{i=1}^n\pi_i:\pi_i>0,~\sum_{i=1}^n\pi_i=1,~\sum_{i=1}^n\pi_i\bg(\bX_i;\btheta)=\bzero\bigg\}.
    \end{equation}
By maximizing $L(\btheta)$ with respect to $\btheta$, one obtains the so-called maximum EL estimator $ \widehat{\btheta}=\arg\max_{\btheta\in\bTheta}L(\btheta)$. Maximizing  (\ref{eq:el}) can be carried out equivalently by solving the corresponding dual problem,  implying 
    \begin{equation}\label{eq:elest1}
    \widehat{\btheta}=\arg\min_{\btheta\in\bTheta}\max_{\blambda\in\widehat{\Lambda}_n(\btheta)}\sum_{i=1}^n\log\{1+\blambda^\T\bg(\bX_i;\btheta)\},
    \end{equation}
where $\widehat{\Lambda}_n(\btheta)=\{\blambda\in\mathbb{R}^r:\blambda^\T\bg(\bX_i;\btheta)\in\mathcal {V},~i=1,\ldots,n\}$ for $\btheta\in\bTheta$ and $\mathcal{V}$ is an open interval containing zero.

In a conventional setting when $p$ and $r$ are fixed as $n\to\infty$, $r\ge p$ is required to ensure that all components of $\btheta$ are identifiable. In high-dimensional cases, however, it is documented in the literature that accommodating a diverging $r$ is a key difficulty for EL; see, among others, \cite{Hjortetal2009}, \cite{ChenPengQin2009}, \cite{LengTang2012}, and \cite{ChangChenChen_2015}. The reason is that the Lagrange multiplier $\blambda \in \mathbb{R}^r$ in (\ref{eq:elest1}) is of the same high dimensionality $r$. Since $|\blambda|_2$ is required to be $o_p(1)$ in theoretical analyses of EL, high-dimensional $r$ is clearly cumbersome. A direct consequence is that dimensionality $p$ and $r$ for EL  in (\ref{eq:el}) can only be accommodated  at some polynomial rate of the sample size $n$ .


To explore EL with high-dimensional statistical problems, let us begin with elucidating their impacts on the EL estimator synthetically  from the sample size $n$,  the number of estimating functions $r$, and the dimensionality of the model parameter $p$.
We first present a general result for the maximum EL estimator $\widehat\btheta$ with $r$ estimating equations.

\begin{proposition}\label{pn:0}
Assume that there exist uniform constants $C_1>0$, $C_2>1$ and $\gamma>2$ such that
    \begin{equation}\label{eq:c1}
    \max_{1\leq j\leq r}\mathbb{E}\bigg\{\sup_{\btheta\in\bTheta}|g_j(\bX_i;\btheta)|^{\gamma}\bigg\}\leq C_1,
    \end{equation}
and
    \begin{equation}\label{eq:c2}
    \begin{split}
    \mathbb{P}\bigg[C_2^{-1}\leq&\inf_{\btheta\in\bTheta}\lambda_{\min}\bigg\{\frac{1}{n}\sum_{i=1}^n\bg(\bX_i;\btheta)\bg(\bX_i;\btheta)^\T\bigg\}\\
    &~~~~~~~~~~~~~~~~~~\leq\sup_{\btheta\in\bTheta}\lambda_{\max}\bigg\{\frac{1}{n}\sum_{i=1}^n\bg(\bX_i;\btheta)\bg(\bX_i;\btheta)^\T\bigg\}\leq C_2\bigg]\rightarrow1.
    \end{split}
    \end{equation}
If $r=o(n^{1/2-1/\gamma})$, then $\widehat{\btheta}$ defined in {\rm(\ref{eq:elest1})} satisfies $|\bar{\bg}(\widehat{\btheta})|_2=O_p(r^{1/2}n^{-1/2})$ where $\bar{\bg}(\widehat{\btheta})=n^{-1}\sum_{i=1}^n\bg(\bX_i;\widehat{\btheta})$.
\end{proposition}

Conditions for Proposition \ref{pn:0} are conventional ones and are mild.   The requirement (\ref{eq:c1}) ensures that some moments with order larger than 2  exist  for the estimating functions, and  (\ref{eq:c2}) says that the sample covariance matrices of the estimating functions should behave reasonably well.    Consistent with the finding in \cite{Hjortetal2009} and  \cite{ChenPengQin2009},  the higher the order of the moment $\gamma$ is,  the more estimating functions can be accommodated.  When the estimating functions are bounded, $\gamma=\infty$, $r$ is allowed to be $o(n^{1/2})$.

The key implication of Proposition \ref{pn:0} is that the sample mean of the estimating functions is well behaving,  regardless the dimensionality of the model parameter $p$ is.  That is, with $r$ unbiased estimating functions,   the optimum  $|\bar{\bg}(\widehat{\btheta})|_2$  is  $O_p(r^{1/2}n^{-1/2})$. Hence  the impact on the behavior of the estimating function  is the dimensionality $r$, which cannot grow faster than $n^{1/2}$ as $n\to \infty$.

 Clearly,  the impact from $p$ on the maximum EL estimator is on the identifiability  of the model parameter.  That is, $\widehat \btheta$  in (\ref{eq:elest1})  is  not uniquely defined when $r<p$   with no further constraints,     rendering  ambiguity and inapplicability of the EL methods for estimating high-dimensional model parameters.
An example of the situation is that the identifiability issue happens in the classical linear models when the model matrix is not of full rank, so that the minimum of the least squares criterion function well exists but the ordinary least squares estimator is not uniquely defined in that case.

To solve the problem, our next objective is to illustrate that identifying a sparse $p$-dimensional model parameter is still feasible.
	


\subsection{ High-dimensional sparse model parameter}

The intuition here is that if one concerns instead a high-dimensional sparse model parameter  $\btheta$ such that most of its components are zeros, then identification and estimation of such a model parameter are feasible with fewer estimating functions by EL  with appropriate penalization.
Specifically, we write $\btheta_0=(\theta_1^0,\ldots,\theta_p^0)^\T$ and let $\mathcal{S}=\{1\leq k\leq p:\theta_k^0\neq 0\}$ with $s=|\mathcal{S}|$.
Here ${\cal S}$ is an unknown set,  and the number of nonzero components $s$ is much smaller than $p$. 
Without loss of generality, we let $\btheta_0=(\btheta_{0,\mathcal{S}}^\T,\btheta_{0,\mathcal{S}^c}^\T)^\T$ where $\btheta_{0,\mathcal{S}}\in\mathbb{R}^{s}$ being the nonzero components and $\btheta_{0,\mathcal{S}^c}=\bzero\in\mathbb{R}^{p-s}$.
For identification of the sparse model parameter, we impose the following condition.

\begin{as}\label{iden} 
Assume that
  \begin{equation}\label{eq:idenlocal}
    \inf_{\btheta\in\{\btheta=(\btheta_{\mathcal{S}}^\T,\btheta_{\mathcal{S}^c}^\T)^\T\in\bTheta:|\btheta_{\mathcal{S}}-\btheta_{0,\mathcal{S}}|_\infty>\varepsilon,\btheta_{\mathcal{S}^c}=\bzero\}}|\mathbb{E}\{\bg(\bX_i;\btheta)\}|_\infty\geq \Delta(\varepsilon) 
    \end{equation} for any $\varepsilon>0$, where $\Delta(\cdot)$ is a positive function satisfying $\lim\inf_{\varepsilon\rightarrow0^+}\varepsilon^{-\beta}\Delta(\varepsilon)\geq K_1$ for some uniform constants $K_1>0$ and $\beta>0$.
\end{as}
The identification condition  (\ref{eq:idenlocal}) can be viewed as a dedicated one for estimating sparse model parameters.  
Condition \ref{iden} is not stringent, and it ensures identifying the nonzero components of $\btheta$ locally.  Studying local optimums in high-dimensional statistical problems is common in the literature with reasonable technical conditions; see, for example,  \cite{LvFan2009} and \cite{Zhang2010}.
 Condition \ref{iden}  means that the mean values of the estimating functions at the truth adequately differ from those outside a small neighborhood of the sparse support of $\btheta_0$.
Here $\beta$ is some generic constant related to  the consistency result in Proposition \ref{pn:00}.  
For estimating a high-dimensional mean parameter with $\bg(\bX;\btheta)=\bX-\btheta$,  we can choose $\Delta(\varepsilon)=\varepsilon$ and $\beta=1$ in Condition \ref{iden}. For linear regression model, $\bg(\bX;\btheta)=\bZ(Y-\bZ^\T\btheta)$ with $\bZ$ and $Y$ being the covariates and response variable respectively, and $\bX=(Y,\bZ^\T)^\T$, we can select $\Delta(\varepsilon)=\varepsilon\|\bSigma_{\bZ,\mathcal{S}}^{-1}\|_\infty^{-1}$ in Condition \ref{iden}, where $\bSigma_{\bZ,\mathcal{S}}=\mathbb{E}(\bZ_{\mathcal{S}}\bZ_{\mathcal{S}}^\T)$. More generally, if there is a subset $\mathcal{E}\subset\{1,\ldots,r\}$ with $|\mathcal{E}|=s$ and $[\mathbb{E}\{\nabla_{\btheta_{\mathcal{S}}}\bg_{\mathcal{E}}(\bX_i;\btheta)\}]^{-1}$ exists where $\bg_{\mathcal{E}}(\cdot)$ collects the set of estimating functions indexed by $\cal E$, then we can select $\Delta(\varepsilon)=\varepsilon\inf_{\btheta\in\{\btheta=(\btheta_{\mathcal{S}}^\T,\btheta_{\mathcal{S}^c}^\T)^\T:\btheta_{\mathcal{S}^c}=\bzero\}}\|[\mathbb{E}\{\nabla_{\btheta_{\mathcal{S}}}\bg_{\mathcal{E}}(\bX_i;\btheta)\}]^{-1}\|_\infty^{-1}$ in Condition \ref{iden}.
Intuitively,  Condition \ref{iden} ensures  the identifiability of  the $s$ nonzero components of $\btheta_0$  so that a consistent sparse estimator is possible as $n\to \infty$, provided $r\geq s$, $r^{1/2}n^{-1/2}\to 0$,  and conditions in Proposition  \ref{pn:00}.


As a special case when ${\cal S}^c$ is empty,
Condition \ref{iden}  for identification becomes a global one for a dense model parameter $\btheta$:
\begin{equation}\label{eq:ginde}
\inf_{\btheta\in\{\btheta\in\bTheta:|\btheta-\btheta_{0}|_\infty>\varepsilon\}}|\mathbb{E}\{\bg(\bX_i;\btheta)\}|_\infty\geq \Delta(\varepsilon),
\end{equation}
where $\Delta(\cdot)$ is a positive function satisfying $\lim\inf_{\varepsilon\rightarrow0^+}\varepsilon^{-\beta}\Delta(\varepsilon)\geq K_1$ for some uniform constants $K_1>0$ and $\beta>0$. Similar global identification conditions can be found in \cite{Chen2007} and \cite{ChenPouzo2012} for some other models.

To estimate a sparse model parameter with unknown zero components, we consider a penalized EL estimator as
    \begin{equation}\label{eq:pest}
    \widetilde{\btheta}_n=\arg\min_{\btheta\in\bTheta}\max_{\blambda\in\widehat{\Lambda}_n(\btheta)}\bigg[\sum_{i=1}^n\log\{1+\blambda^\T\bg(\bX_i;\btheta)\}+n\sum_{k=1}^pP_{1,\pi}(|\theta_k|)\bigg],
    \end{equation}
where $\btheta=(\theta_1,\ldots,\theta_p)^\T$, and $P_{1,\pi}(\cdot)$ is a penalty function with tuning parameter $\pi$. For any penalty function $P_{\tau}(\cdot)$ with tuning parameter $\tau$, let $\rho(t;\tau)=\tau^{-1}P_{\tau}(t)$ for any $t\in[0,\infty)$ and $\tau\in(0,\infty)$.
We assume the penalty function $P_{1,\pi}(\cdot)$ belongs to the following class as considered in \cite{LvFan2009}:
 \begin{equation}\label{eq:peclass}
 \begin{split}
 \mathcal{P}=\{P_{\tau}(\cdot):&~\rho(t;\tau)~\textrm{is increasing in}~t\in[0,\infty)~\textrm{and has continuous derivative}~\rho'(t;\tau)~
 \textrm{for}\\
 &~t\in(0,\infty)~\textrm{with}~\rho'(0^+;\tau)\in(0,\infty),~\textrm{where}~\rho'(0^+;\tau)~\textrm{is independent of}~\tau\}.
 \end{split}
 \end{equation}
The class of penalty function by (\ref{eq:peclass}) is broad and general.
The commonly used $L_1$ penalty, SCAD penalty \citep{FanLi2001} and MCP penalty \citep{Zhang2010} all belong to the class $\mathcal{P}$.
For establishing the consistency of $\widetilde \btheta_n$, we also assume the following condition.
\begin{as} \label{conti}
The function $g_j(\bX;\btheta)$ is continuously differentiable with respect to $\btheta\in\bTheta$ for any $\bX$ and $j=1,\ldots,r$ satisfying the conditions
\begin{equation}\label{eq:ide1}
\max_{1\leq j\leq r}\max_{k\notin \mathcal{S}}\mathbb{E}\bigg\{\sup_{\btheta\in\bTheta}\bigg|\frac{\partial g_j(\bX_i;\btheta)}{\partial \theta_k}\bigg|\bigg\}\leq K_2
\end{equation}
for some uniform constant $K_2>0$, and
 \begin{equation}\label{eq:p2}
    \sup_{\btheta\in\bTheta}\max_{1\leq j\leq r}\max_{k\notin\mathcal{S}}\bigg\{\frac{1}{n}\sum_{i=1}^n\bigg|\frac{\partial g_j(\bX_i;\btheta)}{\partial \theta_k}\bigg|\bigg\}=O_p(\varphi_n)
    \end{equation}
holds for some $\varphi_n>0$, which may diverge with $n$.
\end{as}

Condition \ref{conti} is on the continuity of the estimating function with respect to $\btheta$.  Typically,  smooth estimating functions can be assumed to have bounded derivatives so that Condition \ref{conti} is easily satisfied.   At the sample level, considering the high-dimensionality of the problem,  we can accommodate  diverging $\varphi_n$ in (\ref{eq:p2}) so that our results hold in broad situations.
If there exist envelop functions $B_{n,jk}(\cdot)$ such that $|\partial g_j(\bX;\btheta)/\partial\theta_k|\leq B_{n,jk}(\bX)$ for any $\btheta\in\bTheta$, $j=1,\ldots,r$ and $k\notin\mathcal{S}$, and $|\mathbb{E}\{B_{n,jk}^m(\bX_i)\}|\leq Km!H^{m-2}$ for any $m\geq 2$ and $j=1,\ldots,r$ and $k\notin\mathcal{S}$, where $K$ and $H$ are two uniform positive constants independent of $j$ and $k$. Then by Theorem 2.8 of \cite{Petrov1995}, we know $\sup_{1\leq j\leq r}\sup_{k\notin\mathcal{S}}n^{-1}\sum_{i=1}^nB_{n,jk}(\bX_i)=O_p(1)$ provided that $\max\{\log r,\log p\}=o(n)$. Therefore, (\ref{eq:p2}) holds with $\varphi_n=1$, accommodating exponentially growing dimensionality $r$ and $p$.   
 Since  the identifiability condition (\ref{eq:idenlocal}) only provides a lower bound for the difference between $|\mathbb{E}\{\bg(\bX_i;\btheta)\}|_\infty$ and $0$ when $\btheta=(\btheta_{\mathcal{S}}^\T,\btheta_{\mathcal{S}^c}^\T)^\T$ satisfying $|\btheta_{\mathcal{S}}-\btheta_{0,\mathcal{S}}|_\infty>\varepsilon$ and $\btheta_{\mathcal{S}^c}=\bzero$, we make use of (\ref{eq:ide1}) to derive a lower bound for $|\mathbb{E}\{\bg(\bX_i;\btheta)\}|_\infty$ when $\btheta=(\btheta_{\mathcal{S}}^\T,\btheta_{\mathcal{S}^c}^\T)^\T$ satisfies $\btheta_{\mathcal{S}^c}\neq \bzero$ but $|\btheta_{\mathcal{S}^c}|_1$ is small, and then $\btheta_0$ is a local minimizer for $|\mathbb{E}\{\bg(\bX_i;\btheta)\}|_\infty$.
For special case with linear models,  Condition (\ref{eq:ide1}) becomes one similar to the well known  crucial irrepresentable condition \citep{Zhao07} for sparse linear regression at the population level.
We have the following proposition on the properties of the penalized EL estimator (\ref{eq:pest}).

\begin{proposition}\label{pn:00}
Let $P_{1,\pi}(\cdot)\in\mathcal{P}$ for $\mathcal{P}$ defined in {\rm(\ref{eq:peclass})}. Define $a_n=\sum_{k=1}^pP_{1,\pi}(|\theta_k^0|)$ and $b_n=\max\{rn^{-1},a_n\}$. Assume that {\rm(\ref{eq:c1})},  {\rm(\ref{eq:c2})},   Conditions {\rm\ref{iden}} and {\rm\ref{conti}} hold.
  Suppose that
    \begin{equation}\label{eq:pena1}
    \max_{k\in\mathcal{S}}\sup_{0<t<|\theta_k^0|+c_n}P_{1,\pi}'(t)=O(\chi_n)
    \end{equation}
for some $\chi_n\rightarrow0$ and $c_n\rightarrow0$ with $b_n^{1/(2\beta)}c_n^{-1}\rightarrow0$. If $r=o(n^{1/2-1/\gamma})$, $\max\{b_n,rs\chi_nb_n^{1/(2\beta)}\}=o(n^{-2/\gamma})$ and $r^{1/2}\varphi_n\max\{r^{1/2}n^{-1/2},s^{1/2}\chi_n^{1/2}b_n^{1/(4\beta)}\}=o(\pi)$, then there exists a local minimizer $\widetilde{\btheta}_n\in\bTheta$ for {\rm(\ref{eq:pest})} such that $|\widetilde{\btheta}_{n,\mathcal{S}}-\btheta_{0,\mathcal{S}}|_\infty=O_p\{b_n^{1/(2\beta)}\}$ and $\mathbb{P}(\widetilde{\btheta}_{n,\mathcal{S}^c}=\bzero)\rightarrow1$ as $n\rightarrow\infty$.
\end{proposition}


In Proposition \ref{pn:00}, $a_n$ depends on the truth of the model parameter and the tuning parameter $\pi$ in the penalty function.  For a typical penalty function belonging to (\ref{eq:peclass}) and a model parameter with $s$ nonzero components,  it is the case that $a_n=O(s\pi)\to 0$ as $n\to \infty$.
Requirements on the first derivative of the penalty function via $\chi_n$ is to control the bias introduced by the penalty function $P_{1, \pi}(\cdot)$ on $\widetilde \btheta_n$. See (\ref{eq:pdiff1}) in Section \ref{se:pn00} for details. If we propose the condition $b_n=o(\min_{k\in\mathcal{S}}|\theta_k^0|^{2\beta})$ on the magnitudes of the nonzero components of $\btheta_0$, (\ref{eq:pena1}) can be replaced by
\begin{equation}\label{eq:p1}
    \max_{k\in\mathcal{S}}\sup_{c|\theta_k^0|<t<c^{-1}|\theta_k^0|}P_{1,\pi}'(t)=O(\chi_n)
    \end{equation}
for some constant $c\in(0,1)$.
For those asymptotically unbiased penalty functions like SCAD and MCP,  $\chi_n$ is exactly 0 in (\ref{eq:p1}) for $n$ sufficiently large provided that the nonzero components of $\btheta_0$ are not too small in the sense that the signal strength does not diminish to zero too fast, i.e. $b_n=o(\min_{k\in\mathcal{S}}|\theta_k^0|^{2\beta})$; see  also \cite{FanLi2001}.   
 Hence,    if $\beta=1$ in Condition \ref{iden},   $|\widetilde{\btheta}_{n,\mathcal{S}}-\btheta_{0,\mathcal{S}}|_\infty=O_p(b_n^{1/2})\to 0$ as $n\to \infty$.
Further, if $\pi$ is chosen as $O\{(n^{-1}\log p)^{1/2}\}$, a common one  in the literature,  then $|\widetilde{\btheta}_{n,\mathcal{S}}-\btheta_{0,\mathcal{S}}|_\infty=O_p\{s^{1/2}(n^{-1}\log p)^{1/4}\}$, providing a conservative convergence rate of the estimator $\widetilde \btheta_{n,{\cal S}}$.

 Let $F_n(\btheta)=\max_{\blambda\in\widehat{\Lambda}_n(\btheta)}n^{-1}\sum_{i=1}^n\log\{1+\blambda^\T\bg(\bX_i;\btheta)\}+\sum_{k=1}^pP_{1,\pi}(|\theta_k|)$.  
The rationale of Proposition \ref{pn:00} is that for  $\btheta=(\btheta_{\mathcal{S}}^\T,\btheta_{\mathcal{S}^c}^\T)^\T$ in a small neighborhood of $\btheta_0$ such that $|\btheta_{\mathcal{S}}-\btheta_{0,\mathcal{S}}|_\infty\geq \varepsilon_n$ takes value departing from $\btheta_0$, i.e., $\Delta(\varepsilon_n)$ decays
to zero at some slow enough rate, $F_n(\btheta)$ takes a value larger than $\xi_n F_n(\btheta_0)$ for some diverging $\xi_n$ with probability tending to 1; see also \cite{ChangTangWu2013,ChangTangWu2016} for such a phenomenon
of  EL.  Then with the penalty function encouraging sparsity of $\widetilde\btheta_n$, we are able to establish the consistency of the penalized EL estimator for a sparse model parameter.

Our Proposition \ref{pn:00} shows that the penalized EL can consistently estimate a  high-dimensional model parameter with $p$ growing exponentially with $n$ provided $b_n\to 0$, though the requirement on $r$ remains in a way such that $r=o(n^{1/2})$.
The development of Proposition \ref{pn:00} is fundamentally facilitated by our motivation: to estimate a high-dimensional sparse model parameter. With the new identification condition (\ref{eq:idenlocal}), sparse and consistent estimator can be obtained  by using penalized EL.
%
The intuition of our results is clear: to identify $s$ nonzero components of a sparse $p$-dimensional model parameter,  one essentially requires $r$ $(r\geq s)$ informative estimating functions for those $s$ components.  
The practical interpretation is also clear:  given fewer estimating functions than the model parameters, a reasonable direction is to identify and estimate a sparse model parameter.
Such an observation is consistent with the ones found in \cite{GauTsy_2014} for high-dimensional instrumental variables regression with endogenity where the number of instrumental variables may be less than the model parameters in the regression problems.

\setcounter{equation}{0}
\section{A new penalized empirical likelihood}
\label{s3}

With the penalized EL estimator $\widetilde\btheta_n$ in  (\ref{eq:pest}) capable of handling high-dimensional model parameter with fewer number of estimating functions, our next goal is to accommodate a more general situation: allowing both $r$ and $p$ to grow exponentially with $n$. For such a purpose, we propose to update the penalized EL estimator  with an extra penalty encouraging sparsity in $\blambda$:
    \begin{equation}\label{eq:est}
    \widehat{\btheta}_n=\arg\min_{\btheta\in\bTheta}\max_{\blambda\in\widehat{\Lambda}_n(\btheta)}\bigg[\sum_{i=1}^n\log\{1+\blambda^\T\bg(\bX_i;\btheta)\}-n\sum_{j=1}^rP_{2,\nu}(|\lambda_j|)+n\sum_{k=1}^pP_{1,\pi}(|\theta_k|)\bigg],
    \end{equation}
where $\btheta=(\theta_1,\ldots,\theta_p)^\T$, $\blambda=(\lambda_1,\ldots,\lambda_r)^\T$, and $P_{1,\pi}(\cdot)$ and $P_{2,\nu}(\cdot)$ are two penalty functions with tuning parameters $\pi$ and $\nu$, respectively. Our motivation is that with appropriately chosen penalty function $P_{2,\nu}(\cdot)$ and tuning parameter $\nu$, the estimator $\widehat \btheta_n$ is associated with a sparse Lagrange multiplier $\blambda$.   Since sparse $\blambda$ effectively uses a subset of the estimating functions $\bg(\cdot;\cdot)$,  $r$ itself can be allowed to be large as long as the number of nonzero components in $\blambda$ is small,  essentially satisfying the requirement in Proposition \ref{pn:00}.  Hence, one expects analogous  properties of  (\ref{eq:est}) to those in Proposition \ref{pn:00},  but now being capable of accommodating high-dimensional $p$ and $r$ simultaneously.

Not surprisingly, involving the penalty $P_{2,\nu}(\cdot)$ makes the technical analysis much more challenging, especially   when we are handling exponentially diverging $p$ and $r$ with $n\to\infty$.
For $\btheta\in\bTheta$ and $\blambda\in\widehat{\Lambda}_n(\btheta)$, we define
    \[
    \begin{split}
    f(\blambda;\btheta)=&~\frac{1}{n}\sum_{i=1}^n\log\{1+\blambda^\T\bg(\bX_i;\btheta)\}-\sum_{j=1}^rP_{2,\nu}(|\lambda_j|), \\
    S_n(\btheta)=&~\max_{\blambda\in\widehat{\Lambda}_n(\btheta)}f(\blambda;\btheta)+\sum_{k=1}^pP_{1,\pi}(|\theta_k|).
    \end{split}
    \]
Here $f(\blambda;\btheta)$ is a function of $\blambda$ upon given $\btheta$. Let $ \widehat{\blambda}(\btheta)=\arg\max_{\blambda\in\widehat{\Lambda}_n(\btheta)}f(\blambda;\btheta) $ be the Lagrange multiplier defined at $\btheta\in\bTheta$. For any subset $\mathcal{A}\subset\{1,\ldots,r\}$, we denote by $\bg_{\mathcal{A}}(\bX_i;\btheta)$ the subvector of $\bg(\bX_i;\btheta)$ with components indexed by $\mathcal{A}$. We write $\bar{\bg}_{\mathcal{A}}(\btheta)=n^{-1}\sum_{i=1}^n\bg_{\mathcal{A}}(\bX_i;\btheta)$, $\widehat{\bV}_{\mathcal{A}}(\btheta)=n^{-1}\sum_{i=1}^n\bg_{\mathcal{A}}(\bX_i;\btheta)\bg_{\mathcal{A}}(\bX_i;\btheta)^\T$ and $\bV_{\mathcal{A}}(\btheta)=\mathbb{E}\{\bg_{\mathcal{A}}(\bX_i;\btheta)\bg_{\mathcal{A}}(\bX_i;\btheta)^\T\}$.
For any $\btheta\in\bTheta$ and $j=1,\ldots,r$, define $\bar{g}_j(\btheta)=n^{-1}\sum_{i=1}^ng_j(\bX_i;\btheta)$. We first characterize the properties of $\widehat{\blambda}(\btheta)$ for $\btheta$ near the truth $\btheta_0$. To do this, we assume the following condition for the existence of higher order moments,  a similar  one to the common technical conditions on the tail probability in high-dimensional statistical analysis.

\begin{as}\label{as:moment}
There exist some $K_3>0$ and $\gamma>4$ such that
\[
\max_{1\leq j\leq r}\mathbb{E}\bigg\{\sup_{\btheta\in\bTheta}|g_j(\bX_i;\btheta)|^{\gamma}\bigg\}\leq K_3.
\]
\end{as}

Let $\rho_2(t;\nu)=\nu^{-1}P_{2,\nu}(t)$. We also take $P_{2,\nu}(\cdot)\in\mathcal{P}$ for $\mathcal{P}$ defined in (\ref{eq:peclass}), so that $\rho_2'(0^+;\nu)$ is independent of $\nu$. We write it as $\rho_2'(0^+)$ for simplicity and define $\mathcal{M}_{\btheta}=\{1\leq j\leq r:|\bar{g}_j(\btheta)|\geq \nu\rho_2'(0^+)\}$ for any $\btheta\in\bTheta$. Proposition \ref{pn:2} below shows that for any $\btheta$ near the truth $\btheta_0$,  the support of the Lagrange multiplier $\widehat{\blambda}(\btheta)$ is a subset of $\mathcal{M}_{\btheta}$ with probability approaching one.

\begin{proposition}\label{pn:2}
Let $\{\btheta_n\}$ be a sequence  in $\bTheta$ and $P_{2,\nu}(\cdot)\in\mathcal{P}$ be a convex function for $\mathcal{P}$ defined in {\rm(\ref{eq:peclass})}. For some $C\in(0,1)$, define $\mathcal{M}_{\btheta_n}^*=\{1\leq j\leq r:|\bar{g}_j(\btheta_n)|\geq C\nu\rho_2'(0^+)\}$. Assume Condition {\rm\ref{as:moment}} hold. Further, for the sequence $\{\btheta_n\}$, we assume that   the eigenvalues of $\widehat{\bV}_{\mathcal{M}_{\btheta_n}}(\btheta_n)$ are uniformly bounded away from zero and infinity with probability approaching one, and $|\bar{\bg}_{\mathcal{M}_{\btheta_n}}({\btheta}_n)-\nu\rho_2'(0^+)\sgn\{\bar{\bg}_{\mathcal{M}_{\btheta_n}}(\btheta_n)\}|_2=O_p(u_n)$ for some $u_n\rightarrow0$. Let $\max_{1\leq j\leq r}n^{-1}\sum_{i=1}^n|g_j(\bX_i;\btheta_n)|^2=O_p(\varsigma_n)$ for some $\varsigma_n>0$ that may diverge with $n$. If $m_n^{1/2}u_n\varsigma_n=o(\nu)$ and $m_n^{1/2}u_n n^{1/\gamma}=o(1)$ where $m_n=|\mathcal{M}_{\btheta_n}^*|$, then with probability approaching one there exists a sparse local maximizer $\widehat{\blambda}({\btheta}_n)=(\widehat{\lambda}_{n,1},\ldots,\widehat{\lambda}_{n,r})^\T$ for $f(\blambda;\btheta_n)$ satisfying the three results: {\rm (i)} $|\widehat{\blambda}({\btheta}_n)|_2=O_p(u_n)$, {\rm (ii)}
$\supp\{\widehat{\blambda}({\btheta}_n)\}\subset\mathcal{M}_{\btheta_n}$, and {\rm (iii)} $\sgn(\widehat{\lambda}_{n,j})=\sgn\{\bar{g}_j(\btheta_n)\}$ for any $j\in\mathcal{M}_{\btheta_n}$ with $\widehat{\lambda}_{n,j}\neq 0$.
\end{proposition}

Conditions in Proposition \ref{pn:2} play roles from a few aspects. First, the sequence $\{\btheta_n\}$ can be taken as one that approaches the truth $\btheta_0$ as $n\to\infty$.
Then $\bar{\bg}_{\mathcal{M}_{\btheta_n}}({\btheta}_n)$ will be small when $n$ is large.
As shown in the proof,  $\nu\rho_2'(0^+)\sgn\{\bar{\bg}_{\mathcal{M}_{\btheta_n}}(\btheta_n)\}$ is the asymptotically leading term of $\bar{\bg}_{\mathcal{M}_{\btheta_n}}({\btheta}_n)$.
The reason is that the tuning parameter $\nu$ typically diminishes to $0$ at some slower rate than $n^{-1/2}$,     so that $\nu\rho_2'(0^+)\sgn\{\bar{\bg}_{\mathcal{M}_{\btheta_n}}(\btheta_n)\}$  leads to a non-negligible contribution in the limiting distribution of $\widehat \btheta_n$,  and our analysis shows that it leads to a correctable bias term in $\widehat \btheta_n$.
Upon removing the leading order term, we assume that $|\bar{\bg}_{\mathcal{M}_{\btheta_n}}({\btheta}_n)-\nu\rho_2'(0^+)\sgn\{\bar{\bg}_{\mathcal{M}_{\btheta_n}}(\btheta_n)\}|_2=O_p(u_n)$ with $u_n\to 0$,  which  is a condition that can be easily satisfied.   Requirement on the eigenvalues of  $\widehat{\bV}_{\mathcal{M}_{\btheta_n}}(\btheta_n)$ is natural so that we can characterize the limiting behavior of the estimator $\widehat \btheta_n$.
Furthermore, $m_n$ is taken to be an upper bound of the size of ${\cal M}_{\btheta_n}$,
the generic description such as $m_n^{1/2}u_n\varsigma_n=o(\nu)$ and $m_n^{1/2}u_n n^{1/\gamma}=o(1)$ can be viewed as  characterizing  the capacity of  the  penalized EL under which it is reliable for consistent estimators,  depending on the behavior of the estimating function $\bg(\cdot;\cdot)$ on its continuity and tail probabilistic properties.

Proposition \ref{pn:2} implies that when $\btheta$ is approaching $\btheta_0$, the sparse $\blambda$ in (\ref{eq:est}) effectively conducts a  moments selection by choosing the estimating functions in a way that $\bar{g}_j(\btheta)$ has large absolute deviation from $0$.
Let $\mu_j(\btheta)={\mathbb E}\{g_j (\bX_i;\btheta)\}$,  then we know that $\mu_j(\btheta_0)=0$ and $\bar g_j(\btheta)\to \mu_j(\btheta)$ in probability as $n\to \infty$. If we take $\btheta$ to be in the neighborhood of $\btheta_0$,  then the first order Taylor expansion gives that $\mu_j(\btheta)=\mu_j(\btheta)-\mu_j(\btheta_0)=\{\nabla_{\btheta}\mu_j(\btheta^*)\}^\T(\btheta-\btheta_0)$ for some $\btheta^*$ between $\btheta$ and $\btheta_0$.  Hence, those  components of the estimating functions with large magnitude in the derivative of their expected value with respect  to $\btheta$ will be selected.  Since  larger  derivative indicates a steeper direction towards the truth $\btheta_0$, making it easier and more informative to find the optimum. Therefore,  selecting components in ${\cal M}_{\btheta}$ is seen sensible.   However, we note that without further  strong  and likely to be unrealistic conditions on the shape of the estimating functions, ${\cal M}_{\btheta}$ cannot be controlled as  a fixed set even at the limiting case when $n\to\infty$,   so that it will depend on the value of the parameter $\btheta$.    
Instead of requiring that ${\cal M}_{\btheta}$ to be fixed, we  show  in the following that for any choice of its subset satisfying some reasonable conditions, the resulting penalized EL estimator is consistent and asymptotically normally distributed.

Let
\begin{equation}\label{eq:ell}
\ell_n=\max_{\btheta\in\{\btheta=(\btheta_{\mathcal{S}}^\T,\btheta_{\mathcal{S}^c}^\T)^\T\in\bTheta:|\btheta_{\mathcal{S}}-\btheta_{0,\mathcal{S}}|_\infty\leq c_n,\btheta_{\mathcal{S}^c}=\bzero\}}|\mathcal{M}_{\btheta}|
\end{equation}
 for some $c_n\rightarrow0$ satisfying $b_n^{1/(2\beta)}c_n^{-1}\rightarrow0$ where $b_n$ is more clearly specified in Condition \ref{as:paa1} below. 
Based on Proposition \ref{pn:2}, we know the support of Lagrange multiplier $\widehat{\blambda}(\btheta)$ is a subset of $\mathcal{M}_{\btheta}$ with probability approaching one when $\btheta$ is in a small neighborhood of $\btheta_0$.
Here $\ell_n$ is a technical device controlling the maximum number of effective estimating functions when applying the new penalized EL, and it can be viewed as a cap of the $r$ in Proposition \ref{pn:00}.
Though $\ell_n$ is a technical device, we remark that, practically, one can always achieve the control of the nonzero components of $\blambda$ by appropriately choosing the tuning parameter $\nu$.

To establish the consistency of the penalized EL estimator $\widehat{\btheta}_n$ defined in (\ref{eq:est}), we need the following extra regularity conditions on the continuity and  probabilistic behavior of the estimating functions.


\begin{as}\label{as:smalleigenvalue}
There exist uniform constants $0<K_4<K_5$ such that $K_4<\lambda_{\min}\{\bV_{\mathcal{F}}(\btheta_0)\}\leq \lambda_{\max}\{\bV_{\mathcal{F}}(\btheta_0)\}<K_5$ for any $\mathcal{F}\subset\{1,\ldots,r\}$ with $|\mathcal{F}|\leq \ell_n$, where $\ell_n$ is defined in (\ref{eq:ell}).
\end{as}

\begin{as}\label{as:pa1} 
 Assume that   \[
    \begin{split}
\sup_{\btheta\in\bTheta}\max_{1\leq j\leq r}\max_{k\notin\mathcal{S}}\bigg\{\frac{1}{n}\sum_{i=1}^n\bigg|\frac{\partial g_j(\bX_i;\btheta)}{\partial\theta_k}\bigg|^2\bigg\}=&~O_p(\xi_n),\\
\sup_{\btheta\in\bTheta}\max_{1\leq j\leq r}\max_{k\in\mathcal{S}}\bigg\{\frac{1}{n}\sum_{i=1}^n\bigg|\frac{\partial g_j(\bX_i;\btheta)}{\partial\theta_k}\bigg|^2\bigg\}=&~O_p(\omega_n),\\
\sup_{\btheta\in\bTheta}\max_{1\leq j\leq r}\bigg\{\frac{1}{n}\sum_{i=1}^n|g_j(\bX_i;\btheta)|^4\bigg\}=&~O_p(\varrho_n)
    \end{split}
    \]
for some $\xi_n>0$, $\omega_n>0$ and $\varrho_n>0$ that may diverge with $n$.
\end{as}

\begin{as}\label{as:paa1}
Let $b_n=\max\{a_n,\nu^2\}$ with $a_n=\sum_{k=1}^pP_{1,\pi}(|\theta_k^0|)$. There exist $\chi_n\rightarrow0$ and $c_n\rightarrow0$ with $b_n^{1/(2\beta)}c_n^{-1}\rightarrow0$ for $\beta$ defined in Condition \ref{iden} such that 
    $\max_{k\in\mathcal{S}}\sup_{0<t<|\theta_k^0|+c_n}P_{1,\pi}'(t)=O(\chi_n)$.
\end{as}


Here Condition \ref{as:smalleigenvalue} is actually a weaker one  than that in  (\ref{eq:c2}) in the sense that  it only requires the population covariance matrices of subsets of estimating functions need to well behave at the truth $\btheta_0$.
The first two bounds in Condition \ref{as:pa1} are used to characterize the behavior of the eigenvalues of $\widehat{\bV}_{\mathcal{F}}(\btheta)$ when $\btheta$ in a small neighborhood of $\btheta_0$; see Lemma \ref{la:covc} in Section \ref{se:1}.
We do not impose explicit rate on $\xi_n$, $\omega_n$, and $\varrho_n$, so that the conditions are generally not restrictive. Similar to our earlier discussion for $\varphi_n$ in (\ref{eq:p2}) in Condition \ref{conti}, we can actually choose $\xi_n=\omega_n=\varrho_n=1$ under some additional mild conditions provided that $\max\{\log r,\log p\}=o(n)$. 
Condition \ref{as:paa1} is similar to (\ref{eq:pena1}) in Proposition \ref{pn:00} with a differently defined $b_n$.
Similar to that in Proposition \ref{pn:00}, Condition \ref{as:paa1} can be replaced by (\ref{eq:p1}) if the minimal signal strength condition is satisfied for appropriately chosen tuning parameter $\pi$. Then $\chi_n=0$ when $n$ is large for those asymptotically unbiased penalty functions like SCAD and MCP.


We now present the following theorem for the consistency of $\widehat{\btheta}_n$.

\begin{theorem}\label{tm:consistency}
Let $P_{1,\pi}(\cdot), P_{2,\nu}(\cdot)\in\mathcal{P}$ for $\mathcal{P}$ defined in {\rm(\ref{eq:peclass})}, and $P_{2,\nu}(\cdot)$ be a convex function with bounded second derivative around $0$. Assume Conditions {\rm\ref{iden}--\ref{as:paa1}} hold. Let $b_n=\max\{a_n,\nu^2\}$ with $a_n=\sum_{k=1}^pP_{1,\pi}(|\theta_k^0|)$, and $\kappa_n=\max\{\ell_n^{1/2}n^{-1/2},s^{1/2}\chi_n^{1/2}b_n^{1/(4\beta)}\}$. If $\log r=o(n^{1/3})$, $\varrho_n=o(n^2)$, $s^2\ell_n\omega_nb_n^{1/\beta}=o(1)$, $\ell_n^2n^{-1}\varrho_n\log r=o(1)$, $\max\{b_n,\ell_n\kappa_n^2\}=o(n^{-2/\gamma})$, $\ell_n^{1/2}\varrho_n^{1/2}\kappa_n=o(\nu)$ and $\ell_n^{1/2}\xi_n^{1/2}\max\{\ell_n\nu,s^{1/2}\chi_n^{1/2}b_n^{1/(4\beta)}\}=o(\pi)$, then there exists a local minimizer $\widehat{\btheta}_n\in\bTheta$ for {\rm(\ref{eq:est})} such that $|\widehat{\btheta}_{n,\mathcal{S}}-\btheta_{0,\mathcal{S}}|_\infty=O_p\{b_n^{1/(2\beta)}\}$ and $\mathbb{P}(\widehat{\btheta}_{n,\mathcal{S}^c}=\bzero)\rightarrow1$ as $n\rightarrow\infty$.\end{theorem}

Theorem \ref{tm:consistency} establishes  the consistency of $\widehat{\btheta}_n$ in the sense that $|\widehat{\btheta}_n-\btheta_0|_\infty\xrightarrow{p}0$. The convergence rate $O_p\{b_n^{1/(2\beta)}\}$ is a conservative one before we establish the asymptotic normality of the penalized EL estimator $\widehat{\btheta}_{n,{\cal S}}$ later. Under additional regularity conditions, such a rate can be improved as $O_p(\nu)$.
Results in Theorem \ref{tm:consistency} holds for broad situations accommodating various cases of the estimating functions.
In   reasonable cases that we discussed earlier,  $\chi_n=0$ and $\xi_n=\omega_n=\varrho_n=1$.  
Theorem \ref{tm:consistency} holds provided that $\log r=o(n^{1/3})$, $\ell_n=o(\min\{n^{1/2}(\log r)^{-1/2},n^{1/2-1/\gamma}\})$, $a_n=o(\min\{s^{-2\beta}\ell_n^{-\beta},n^{-2/\gamma}\})$, and the tuning parameters $\nu$ and $\pi$ satisfy $\ell_nn^{-1/2}=o(\nu)$, $\nu=o(\min\{s^{-\beta}\ell_n^{-\beta/2},n^{-1/\gamma}\})$ and $\ell_n^{3/2}\nu=o(\pi)$.
Noticing that $a_n\lesssim s\pi$,  
by choosing $\pi=o(\min\{s^{-2\beta-1}\ell_n^{-\beta},s^{-1}n^{-2/\gamma}\})$ can ensure the consistency result.
%
Additionally, we note that  $s\leq\ell_n$.  Thus by letting
$\log r\asymp n^{\tau}$ and $\ell_n\asymp n^{\delta}$ for some $\tau\in[0,\frac{1}{3})$ and $\delta\in[0,\min\{\frac{\gamma-4}{7\gamma},\frac{1}{6\beta+7}\})$, $\widehat{\btheta}_n$ satisfies Theorem \ref{tm:consistency} if $\nu\asymp n^{-\phi_1}$ and $\pi\asymp n^{-\phi_2}$ with  $\phi_1\in (\max\{\frac{3\beta\delta}{2},\frac{1}{\gamma}\},\frac{1}{2}-\delta)$ and $\phi_2\in(\max\{(3\beta+1)\delta,\frac{2}{\gamma}+\delta\},\phi_1-\frac{3\delta}{2})$, which are reasonable choices for the tuning parameters.

To further establishing the limiting distribution of $\widehat{\btheta}_{n,\mathcal{S}}$, we need the following two additional conditions.

\begin{as}\label{as:partial2}
For each $j=1,\ldots,p$, $g_j(\bX;\btheta)$ is twice continuously differentiable with respect to $\btheta$ in $\bTheta$ for any $\bX$, and
\[
\sup_{\btheta\in\bTheta}\max_{1\leq j\leq r}\max_{k_1,k_2\in\mathcal{S}}\bigg\{\frac{1}{n}\sum_{i=1}^n\bigg|\frac{\partial^2g_j(\bX_i;\btheta)}{\partial\theta_{k_1}\partial\theta_{k_2}}\bigg|^2\bigg\}=O_p(\varpi_n)
\]
for some $\varpi_n\geq0$ that may diverge with $n$.
\end{as}

\begin{as}\label{as:eig2}
Let $\bQ_{\mathcal{F}}=[\mathbb{E}\{\nabla_{\btheta_{\mathcal{S}}}\bg_{\mathcal{F}}(\bX_i;\btheta_0)\}]^\T[\mathbb{E}\{\nabla_{\btheta_{\mathcal{S}}}\bg_{\mathcal{F}}(\bX_i;\btheta_0)\}]$ for any $\mathcal{F}\subset\{1,\ldots,r\}$. There exist uniform constants $0<K_6<K_7$ such that $K_6<\lambda_{\min}(\bQ_{\mathcal{F}})\leq\lambda_{\max}(\bQ_{\mathcal{F}})<K_7$ for any $\mathcal{F}$ with $s\leq |\mathcal{F}|\leq \ell_n$.
\end{as}

Following similar discussion for Condition \ref{as:pa1},  $\varpi_n=1$ in Condition \ref{as:partial2} for reasonable models in practice. Let $\mathcal{R}_n=\textrm{supp}\{\widehat{\blambda}(\widehat{\btheta}_n)\}$ and define
    \begin{equation}\label{eq:bias}
    \begin{split}
    \widehat{\bJ}_{\mathcal{R}_n}=&~\{\nabla_{\btheta_{\mathcal{S}}}\bar{\bg}_{\mathcal{R}_n}(\widehat{\btheta}_n)\}^\T\widehat{\bV}_{\mathcal{R}_n}^{-1}(\widehat{\btheta}_n)\{\nabla_{\btheta_{\mathcal{S}}}\bar{\bg}_{\mathcal{R}_n}(\widehat{\btheta}_n)\},\\
    \widehat{\bpsi}_{\mathcal{R}_n}=&~\widehat{\bJ}_{\mathcal{R}_n}^{-1}\{\nabla_{\btheta_{\mathcal{S}}}\bar{\bg}_{\mathcal{R}_n}(\widehat{\btheta}_n)\}^\T \widehat{\bV}_{\mathcal{R}_n}^{-1}(\widehat{\btheta}_n)\bigg\{\frac{1}{n}\sum_{i=1}^n\frac{\bg_{\mathcal{R}_n}(\bX_i;\widehat{\btheta}_n)}{1+\widehat{\blambda}(\widehat{\btheta}_n)^\T\bg(\bX_i;\widehat{\btheta}_n)}\bigg\}.
    \end{split}
    \end{equation}
We have the following limiting distribution for $\widehat{\btheta}_{n,\mathcal{S}}$.

\begin{theorem}\label{tm:2}
Let $P_{1,\pi}(\cdot), P_{2,\nu}(\cdot)\in\mathcal{P}$ for $\mathcal{P}$ defined in {\rm(\ref{eq:peclass})}, and $P_{2,\nu}(\cdot)$ be a convex function with bounded second derivative around $0$. Assume Conditions {\rm\ref{iden}--\ref{as:eig2}} hold. Let $b_n=\max\{a_n,\nu^2\}$ with $a_n=\sum_{k=1}^pP_{1,\pi}(|\theta_k^0|)$, and $\kappa_n=\max\{\ell_n^{1/2}n^{-1/2},s^{1/2}\chi_n^{1/2}b_n^{1/(4\beta)}\}$. If $\log r=o(n^{1/3})$, $\varrho_n=o(n^2)$, $b_n=o(n^{-2/\gamma})$, $ns\chi_n^2=o(1)$, $\ell_n^2\varrho_n^{1/2}(\log r)\max\{s^2(\omega_n+s\varpi_n)b_n^{1/\beta},n^{-1}(s\omega_n+\ell_n\varrho_n)\log r\}=o(1)$, $n\ell_n\kappa_n^4\max\{s\omega_n,n^{2/\gamma}\}=o(1)$, $n\ell_ns^2\varpi_n\max\{\ell_n^2\nu^4,s^2\chi_n^2b_n^{1/\beta}\}=o(1)$, $\ell_n^{1/2}\varrho_n^{1/2}\kappa_n=o(\nu)$ and $\ell_n^{1/2}\xi_n^{1/2}\max\{\ell_n\nu,s^{1/2}\chi_n^{1/2}b_n^{1/(4\beta)}\}=o(\pi)$, then local minimizer $\widehat{\btheta}_n\in\bTheta$ for {\rm(\ref{eq:est})} specified in Theorem {\rm\ref{tm:consistency}} satisfies
\begin{equation}\label{eq:asym}
{n}^{1/2}\balpha^\T\widehat{\bJ}_{\mathcal{R}_n}^{1/2}(\widehat{\btheta}_{n,\mathcal{S}}-\btheta_{0,\mathcal{S}}-\widehat{\bpsi}_{\mathcal{R}_n})\xrightarrow{d}N(0,1)
\end{equation}
as $n\rightarrow\infty$, where $\widehat{\bJ}_{\mathcal{R}_n}$ and $\widehat{\bpsi}_{\mathcal{R}_n}$ are defined in {\rm(\ref{eq:bias})}.
\end{theorem}

Theorem \ref{tm:2} shows that
subject to a bias correction, the penalized EL estimator for nonzero components is asymptotically normal in the sense of (\ref{eq:asym}).
The bias term $\widehat{\bpsi}_{\mathcal{R}_n}$ in (\ref{eq:asym}) is due to the penalty function $P_{2,\nu}(\cdot)$ used in (\ref{eq:est}); see also our discussion after the Proposition \ref{pn:2}. Write $\widehat{\blambda}(\widehat{\btheta}_n)=(\widehat{\lambda}_1,\ldots,\widehat{\lambda}_r)^\T$. Furthermore,   as shown in (\ref{eq:ke1}) in Section \ref{s7},  the correctable bias term is $\widehat{\bpsi}_{\mathcal{R}_n}=\widehat{\bJ}_{\mathcal{R}_n}^{-1}\{\nabla_{\btheta_{\mathcal{S}}}\bar{\bg}_{\mathcal{R}_n}(\widehat{\btheta}_n)\}^\T \widehat{\bV}_{\mathcal{R}_n}^{-1}(\widehat{\btheta}_n)\widehat{\bfeta}_{\mathcal{R}_n}$ where $\widehat{\bfeta}=(\widehat{\eta}_1,\ldots,\widehat{\eta}_r)^\T$ with $\widehat{\eta}_j=\nu \rho_{2}'(|\widehat{\lambda}_j|;\nu)\textrm{sgn}(\widehat{\lambda}_j)$ for $\widehat{\lambda}_j\neq0$ and $\widehat{\eta}_j\in[-\nu\rho_2'(0^+),\nu\rho_2'(0^+)]$ for $\widehat{\lambda}_j=0$.

Similar to that in Theorem \ref{tm:consistency},  with reasonable cases $\chi_n=0$ and $\xi_n=\omega_n=\varrho_n=\varpi_n=1$,  descriptions on the dimensionality in Theorem \ref{tm:2} can be simplified.
If $\ell_n\asymp s$,   
Theorem \ref{tm:2} holds provided that $\log r=o(n^{1/3})$, $s=o(\min\{n^{1/3}(\log r)^{-2/3},n^{1/(10\beta+7)}(\log r)^{-2\beta/(10\beta+7)},n^{(\gamma-4)/(7\gamma)}\})$, and $\nu$ and $\pi$ satisfying $sn^{-1/2}=o(\nu)$,
 $\nu=o(\min\{n^{-1/\gamma},s^{-5\beta/2}(\log r)^{-\beta/2},n^{-1/4}s^{-5/4}\})$, $s^{3/2}\nu=o(\pi)$ and
 $\pi=o(\min\{n^{-2/\gamma}s^{-1},s^{-5\beta-1}(\log r)^{-\beta}\})$.

Generally speaking, conditions in Theorem \ref{tm:2} is stronger than those in Theorem \ref{tm:consistency}, which can be viewed as the expense for the stronger asymptotic normality results.  In summary, we have established that the sparse penalized EL estimator (\ref{eq:est}) has desirable properties including consistency in estimating nonzero components and identifying zero components of $\btheta_0$, and asymptotic normality for the estimator of the nonzero components of $\btheta_0$.

%
%
%
%
%
%
%

\setcounter{equation}{0}
\section{Algorithms for implementations}
\label{s4}

For ease and stability in implementations, we calculate the penalized EL estimator $\widehat \btheta_n$ by minimizing the following slightly modified objective function:
    \bea
    {\widehat \btheta_n}= \arg \min_{\btheta\in\bTheta} \max_{\blambda\in\widehat{\Lambda}_n(\btheta)} \bigg[\sum_{i=1}^n \log_\star \{1+\blambda^\T\bg(\bX_i;\btheta)\}- n \sum_{j=1}^r P_{2,\nu}(|\lambda_j|) + n \sum_{k=1}^p P_{1,\pi}(|\theta_k|)\bigg], \label{eq:pel}
    \eea
where $\log_\star(z)$ is a twice differentiable pseudo-logarithm function with bounded support  adopted from \cite{Owen2001}:
    \begin{align}
    \log_\star(z)=\begin{cases}
    \log (z) & \mbox{ if } z\ge \epsilon; \\
    \log(\epsilon)-1.5+2z/\epsilon -z^2/(2\epsilon^2) & \mbox{ if } z\leq \epsilon;
    \end{cases}
    \end{align}
		where $\epsilon$ is chosen as $1/n$ in our implementations.
Here  $P_{1,\pi}(\cdot)$ and $P_{2,\nu}(\cdot)$ are two penalty functions with tuning parameters $\pi$ and $\nu$, respectively.  In the optimization, we apply the quadratic approximation \citep{FanLi2001} to the penalty functions $P_{1,\pi}(\cdot)$ and $P_{2,\nu}(\cdot)$. More specifically, for a penalty function $P_{\tau}(\cdot)$, the quadratic approximation states
    \bea \label{lqa}
    P_\tau(|t|)\approx P_\tau(|t_0|)+\dfrac{1}{2} {P_\tau'(|t_0|) \over |t_0|} (t^2-t_0^2)
    \eea
for $t$ being in a small neighborhood of $t_0$. The first and second derivatives are  approximated by
    $$P'_\tau(|t|) \approx {P'_\tau(|t_0|) \over |t_0|} \cdot t \ \ {\rm and \ \ } P''_\tau(|t|) \approx {P'_\tau(|t_0|) \over |t_0|}.$$

The computation of EL is a challenging aspect, especially with high-dimensional $p$ and $r$. To compute the penalized EL estimator $\widehat\btheta_n$, we propose to apply a modified two-layer coordinate decent algorithm extending the one in \cite{tangwu2014}. The inner layer of the algorithm solves for $\blambda$ with given $\btheta$ by maximizing $f(\blambda;\btheta)$ as given in Section \ref{s3}. This layer only involves maximizing a concave function, and hence is stable. The outer layer of the algorithm searches for the optimizer $\widehat \btheta_n$. Both layers can be solved using coordinate descent by cycling through and updating each of the coordinates; see \cite{tangwu2014}.

In the inner layer, $\blambda$ is solved at a given $\btheta$, which can be done by optimizing (\ref{eq:pel}) with respect to $\blambda$ using coordinate descent. Suppose that $\blambda$ starts at an initial value $\widehat{\blambda}^{(0)}$. With the other coordinates fixed, the $(m+1)$th Newton's update for $\lambda_j$  $(j=1,\dots,r)$, the $j$th component of $\blambda$, is given by
    \bea \label{cd_lambda}
    {\widehat \lambda}_j^{(m+1)} = {\widehat \lambda}_j^{(m)} - \frac{\sum_{i=1}^n \log_\star'(t_i^{(m)}) g_j(\bX_i;\btheta) - n P'_{2,\nu}(|\widehat{\lambda}_j^{(m)}|)}
    {\sum_{i=1}^n \log_\star''(t_i^{(m)}) \{g_j(\bX_i;\btheta)\}^2 - n P''_{2,\nu}(|\widehat{\lambda}_j^{(m)}|)},
    \eea
where $t_i^{(m)}=1+\bg(\bX_i;\btheta)^\T{\widehat \blambda}^{(m)}$ with $\widehat{\blambda}^{(m)}=(\widehat{\lambda}_1^{(m)},\ldots,\widehat{\lambda}_r^{(m)})^\T$. The procedure cycles through all the $r$ components of $\blambda$ and is repeated until convergence. During this process, the objective function needs to be checked to ensure it gets optimized in each step. If not, the step size continues to be halved until the objective function gets driven in the right direction. The iterative updating procedure (\ref{cd_lambda}) can be viewed as sequential univariate optimizations. The convergence rate and stability are studies in the optimization literature; see for example \cite{friedman07} and \cite{wu2008}.

The outer layer of the algorithm is to optimize (\ref{eq:pel}) with respect to the parameter $\btheta$, the main interest of the  penalized EL, using the coordinate descent algorithm. At a given $\blambda$, the algorithm updates $\theta_k$  $(k=1,\dots,p)$, by minimizing $S_n(\btheta)$ defined in Section \ref{s3} with respect to $\theta_k$ with other $\theta_l$ $(l \neq k)$ fixed. Suppose that $\btheta$ starts at an initial value $\widehat{\btheta}^{(0)}$. The $(m+1)$th update for $\theta_k$ is given by
    \bea \label{cd_theta}
    & & {\widehat \theta}_k^{(m+1)} = {\widehat \theta}_k^{(m)} -  \frac{\sum_{i=1}^n \log_\star'(s_i^{(m)}) w_{ik}^{(m)} + n P'_{1,\tau}(|\widehat{\theta}_k^{(m)}|) }
    { \sum_{i=1}^n [ \log_\star''(s_i^{(m)}) \{ w_{ik}^{(m)} \}^2 + \log_\star'(s_i^{(m)}) z_{ik}^{(m)} ] + n  P''_{1,\tau}(|\widehat{\theta}_k^{(m)}|) },
    \eea
where $s_i^{(m)}=1+\blambda^\T \bg(\bX_i;{\widehat \btheta}^{(m)})$, $w_{ik}^{(m)}=\blambda^\T \partial \bg(\bX_i;{\widehat \btheta}^{(m)})/\partial \theta_k$ and $z_{ik}^{(m)}=\blambda^\T \partial^2 \bg(\bX_i;{\widehat \btheta}^{(m)})/\partial \theta_k^2$ with $\widehat{\btheta}^{(m)}=(\widehat\theta_1^{(m)},\ldots,\widehat{\theta}_p^{(m)})^\T$. Since quadratic approximations are applied in the algorithms, we follow \cite{FanLi2001} and set a component $\widehat \lambda_j^{(m)}$ or $\widehat \theta_k^{(m)}$ as zero when it is less than a threshold level say $10^{-3}$ in an iteration.

We summarize the computation procedure for $\btheta$ and $\blambda$ in the following pseudo-code. Suppose $\xi$ is a pre-defined small number, say, $\xi=10^{-4}$.

\framebox{\parbox[t]{5in}{
    \begin{enumerate}
    \item Set the iteration counter $m=0$, and initialize $\widehat \btheta^{(0)}$ and $\widehat \blambda^{(0)}$;
    \item Define the $\bg(\bX_i;\btheta)$ function;
    \item (Outer layer) For $k=1,\dots,p$,
        \begin{enumerate}
        \item Calculate ${\widehat \theta}_k^{(m+1)}$ as in (\ref{cd_theta});
        \item (Inner layer) For $j=1,\dots,r$, update $\widehat \lambda_j^{(m)}$ as $\widehat \lambda_j^{(m+1)}$ defined in (\ref{cd_lambda});
        \end{enumerate}
    \item If $\max_{1 \leq k \leq p}|\widehat{\theta}_k^{(m+1)}-\widehat{\theta}_k^{(m)}| < \xi$, then stop;
    \item Otherwise repeat steps 3 through 4.
    \end{enumerate}
}} \\

\setcounter{equation}{0}
\section{Numerical examples}
\label{s5}

\subsection{Estimating high-dimensional mean parameter}

The first simulation study is to calculate the mean of a multivariate normal distribution in $\mathbb{R}^p$. Let $\bX=(X_1,\ldots,X_p)^\T\sim N(\btheta_0,\bSigma)$. Suppose only three elements, $X_1,X_2$, and $X_5$, have nonzero means and the rest $p-3$ elements have zero means, i.e., $\btheta_0=(5,4,0,0,1,0,\dots,0)^\T$. The covariance matrix $\bSigma=(\sigma_{kl})_{p\times p}$ is set as $\sigma_{kk}=1$ for each $k=1,\ldots,p$ and $\sigma_{kl}=0.9$ for any $k\neq l$. The estimating function is simply $\bg(\bX;\btheta)=\bX-\btheta$. In this case, the number of parameters $p$ is equal to the number of estimating equations $r$. We consider the underdetermined case where $p=r>n$. We generate 100 random samples. The SCAD penalty \citep{FanLi2001} is used for both the penalty functions $P_{1,\pi}(\cdot)$ and $P_{2,\nu}(\cdot)$ in (\ref{eq:est}) for all the numerical experiments in this paper. Since local quadratic approximation is applied in the algorithms, the convexity requirements of the results in Sections \ref{s2} and \ref{s3} are met.

Table \ref{table_pmean} summarizes the results for $(n,p)=(50,100)$, $(100,200)$, and $(100,500)$. The proposed penalized EL with two penalties (namely, PEL2) is compared to the single penalty approach (PEL) discussed in \cite{TangLeng2010}. Three information criteria for choosing  the tuning parameters $\pi$ and $\nu$ in the penalty functions -- BIC \citep{Schwarz1978}, BICC \citep{WangLiLeng2009JRSSB}, and EBIC \citep{ChenChen2008Bioka} -- are used.
In general, all the three BIC-type criteria work similarly, with the latter two yield slightly fewer nonzero parameters. The results from MLE for all $p$ variables and the three true variables (i.e., MLE-Oracle) are also reported. The column of $\btheta_{\rm nonzero}$ reports the average number of selected nonzero components. The column of $\btheta_{\rm true}$ reports the average number of true nonzero components that are selected. The difference is the average number of false predictors that get selected. The next column reports the model error (ME), which is defined by ${\rm ME}=({\widehat \btheta}-\btheta)^\T({\widehat \btheta}-\btheta)$
	for a given estimator $\widehat{\btheta}$.
A smaller ME means a better estimation and prediction. The last column reports the number of selected estimating equations. Obviously, in the single penalty approach, all equating equations are used since no selection is performed. In each cell, standard error appears in the parentheses.

It is clear from the table that the double-penalty approach outperforms the single-penalty approach, as expected. A much smaller subset of variables get selected with almost all the three true predictors identified by the double-penalty method. That says, the double-penalty approach yields lower false positives and higher true positives. While in the single-penalty approach, fewer true predictors are chosen in the larger set of selected variables or nothing can be picked out if $p\gg n$. What is the most interesting is that a small number (on average 5-8) of estimating equations are selected in the double-penalty approach. As a result, the double-penalty method yields a much smaller ME than the single-penalty method.

\begin{table}
\begin{center}
\tabcolsep 0.1in \arrayrulewidth 1pt \doublerulesep 2pt
\begin{tabular}{ l l c c c c c c } \hline
$(n,p,r)$ & Method & $\btheta_{\rm nonzeros}$ & $\btheta_{\rm true}$ & ME & No. EE's \\ \hline
$(50,100,200)$ & MLE-Oracle & 3 (0) & NA & 0.062 (0.009) & NA \\
& MLE & 100 (0) & 3 (0) & 2.096 (0.287) & NA \\
& PEL-BIC & 24.06 (4.13) & 0.72 (0.12) & 33.276 (1.507) & 100 (0) \\
& PEL-BICC & 23.15 (4.08) & 0.69 (0.12) & 33.635 (1.483) & 100 (0) \\
& PEL-EBIC & 23.15 (4.08) & 0.69 (0.12) & 33.635 (1.483) & 100 (0) \\
& PEL2-BIC & 3.41 (0.17) & 2.81 (0.04) & 0.332 (0.041) & 5.11 (0.34) \\
& PEL2-BICC & 3.29 (0.15) & 2.80 (0.04) & 0.302 (0.041) & 6.13 (0.33) \\
& PEL2-EBIC & 3.15 (0.13) & 2.76 (0.05) & 0.341 (0.052) & 8.20 (0.21) \\
\hline
$(100,200,400)$ & MLE-Oracle & 3 (0) & NA & 0.024 (0.003) & NA \\
& MLE & 200 (0) & 3 (0) & 1.743 (0.179) & NA \\
& PEL-BIC & 22.02 (6.02) & 0.33 (0.09) & 38.078 (1.073) & 199.98 (0.02) \\
& PEL-BICC & 22.02 (6.02) & 0.33 (0.09) & 38.078 (1.073) & 199.98 (0.02) \\
& PEL-EBIC & 22.02 (6.02) & 0.33 (0.09) & 38.078 (1.073) & 199.98 (0.02) \\
& PEL2-BIC & 6.41 (1.84) & 2.84 (0.04) & 0.333 (0.091) & 6.67 (0.23) \\
& PEL2-BICC & 6.18 (1.84) & 2.82 (0.04) & 0.352 (0.092) & 6.64 (0.23) \\
& PEL2-EBIC & 5.82 (1.86) & 2.80 (0.04) & 0.372 (0.094) & 6.69 (0.24) \\
\hline
$(100,500,1000)$ & MLE-Oracle & 3 (0) & NA & 0.031 (0.005) & NA \\
& MLE & NA & NA & NA & NA \\
& PEL-BIC & 85.71 (22.69) & 0.51 (0.14) & 37.585 (1.193) & 500 (0) \\
& PEL-BICC & 0 (0) & 0 (0) & 42 (0) & 500 (0) \\
& PEL-EBIC & 0 (0) & 0 (0) & 42 (0) & 500 (0) \\
& PEL2-BIC & 2.88 (0.11) & 2.70 (0.06) & 0.356 (0.057) & 6.40 (0.36) \\
& PEL2-BICC & 2.82 (0.09) & 2.70 (0.06) & 0.376 (0.058) & 6.53 (0.35) \\
& PEL2-EBIC & 2.83 (0.09) & 2.71 (0.06) & 0.369 (0.058) & 6.97 (0.32) \\
\hline
\end{tabular}
\caption{Simulation results for mean of a normal distribution based on 100 random samples. Here $\btheta_{\rm nonzero}$ is the average number of selected nonzero components,  $\btheta_{\rm true}$ is the average number of true nonzero components that are selected, ME reports the model error, and No.EE's reports the number of estimating equations selected.   \label{table_pmean}}
\end{center}
\end{table}

\subsection{Linear regression}\label{sim2}

In this simulation study, we consider a linear regression model
    $ Y_i = \bZ_i^\T\btheta_0 + \varepsilon_i, $
where $\btheta_0=(3,1.5,0,0,2,0,\dots,0)^\T$, $\bZ_i \in \mathbb{R}^p$ are generated from $N(\bzero,\bSigma)$ with  $\sigma_{kk}=1$ for any $k=1,\ldots,p$ and $\sigma_{kl}=0.5$ for any $k\neq l$, where $\bSigma=(\sigma_{kl})_{p\times p}$, and $\varepsilon_i$ is a standard normal distributed random variable. Write $\bX_i=(Y_i,\bZ_i^\T)^\T$. The estimating function is $\bg(\bX;\btheta)=\bZ(Y-\bZ^\T\btheta)$ with $p=r$.

The model error (ME) in the regression setting is defined by ${\rm ME}=({\widehat \btheta}-\btheta)^\T\bSigma({\widehat \btheta}-\btheta)$
	for a given estimator $\widehat{\btheta}$.
Table \ref{table_preg} reports the results for $(n,p,r)=(50,100,100), (100,200,200)$, and $(100,500,500)$ with the columns defined in the same way as those in Table \ref{table_pmean}.
 Similar to the previous example, the single-penalty approach (PEL)  of  \cite{TangLeng2010} is compared with the double-penalty approach (PEL2) together with the three BIC criteria for selecting the tuning parameter(s).
We also compare our method with the LASSO method with $L_1$ penalty.
Since the number of parameters $p$ doubles the number of subjects $n$, the MLE method does not work in this example. We only report the results from MLE-Oracle (i.e., the MLE method using the true predictors), which gives the smallest model error. In all the three settings, the single-penalty method fails to select any predictor when using all $r$ estimating equations. The double-penalty method identifies all true predictors from a handful of selected ones in most cases by using only a few estimating equations.
With the default tuning parameter selection method in the LASSO, we clearly see that the number of false inclusion of the predictors is high.  Hence,  compared with LASSO method,  we observe  that our method has better performance in recovering a sparse model.
\begin{table}
\begin{center}
\tabcolsep 0.1in \arrayrulewidth 1pt \doublerulesep 2pt
\begin{tabular}{ l l c c c c c c } \hline
$(n,p,r)$ & Method & $\btheta_{\rm nonzeros}$ & $\btheta_{\rm true}$ & ME & No. EE's \\ \hline
$(50,100,100)$ & MLE-Oracle & 3 (0) & NA & 0.069 (0.005) & NA \\
& LASSO & 15.21 (0.88) & 3 (0) & 0.439 (0.034) & NA \\
& PEL-BIC & 0 (0) & 0 (0) & 28.75 (0) & 100 (0) \\
& PEL-BICC & 0 (0) & 0 (0) & 28.75 (0) & 100 (0) \\
& PEL-EBIC & 0 (0) & 0 (0) & 28.75 (0) & 100 (0) \\
&PEL2-BIC & 6.39 (0.52) & 2.98 (0.02) & 0.497 (0.069) & 10.46 (0.46) \\
&PEL2-BICC & 6.33 (0.52) & 2.98 (0.02) & 0.498 (0.069) & 10.49 (0.46) \\
&PEL2-EBIC & 6.06 (0.52) & 2.97 (0.02) & 0.531 (0.07) & 10.43 (0.47) \\
\hline
$(100,200,200)$ & MLE-Oracle & 3 (0) & NA & 0.047 (0.005) & NA \\
& LASSO & 17.79 (0.87) & 3 (0) & 0.374 (0.019) & NA \\
& PEL-BIC & 0 (0) & 0 (0) & 28.75 (0) & 200 (0) \\
& PEL-BICC & 0 (0) & 0 (0) & 28.75 (0) & 200 (0) \\
& PEL-EBIC & 0 (0) & 0 (0) & 28.75 (0) & 200 (0) \\
& PEL2-BIC & 9.22 (1.27) & 3 (0) & 0.647 (0.118) & 5.38 (0.17) \\
& PEL2-BICC & 9.28 (1.28) & 3 (0) & 0.651 (0.119) & 5.39 (0.17) \\
& PEL2-EBIC & 8.38 (1.03) & 3 (0) & 0.632 (0.119) & 5.34 (0.17) \\
\hline
$(100,500,500)$ & MLE-Oracle & 3 (0) & NA & 0.039 (0.003) & NA \\
& LASSO & 23.79 (1.23) & 3 (0) & 0.507 (0.028) & NA \\
& PEL-BIC & 0 (0) & 0 (0) & 28.75 (0) & 500 (0) \\
& PEL-BICC & 0 (0) & 0 (0) & 28.75 (0) & 500 (0) \\
& PEL-EBIC & 0 (0) & 0 (0) & 28.75 (0) & 500 (0) \\
& PEL2-BIC & 6.28 (1.31) & 3 (0) & 0.601 (0.083) & 5.48 (0.16) \\
& PEL2-BICC & 5.96 (1.31) & 3 (0) & 0.593 (0.085) & 5.38 (0.17) \\
& PEL2-EBIC & 6.04 (1.32) & 3 (0) & 0.602 (0.086) & 5.41 (0.16) \\
\hline
\end{tabular}
\caption{Simulation results for linear regression based on 100 replicates. Here $\btheta_{\rm nonzero}$ is the average number of selected nonzero components,  $\btheta_{\rm true}$ is the average number of true nonzero components that are selected, ME reports the model error, and No.EE's reports the number of estimating equations selected.  \label{table_preg}}
\end{center}
\end{table}

\subsection{Regression model with repeated measures}\label{sim3}

This is an example with more estimating equations than the number of parameters, i.e., $r>p$.
Now we consider a repeated measures model such that
    $y_{ij} = \bz_{ij}^\T\btheta_0+\epsilon_{ij}~(i=1,\dots,n; j=1,2),$
where $\btheta_0=(3,1.5,0,0,2,0,\dots,0)^\T \in \mathbb{R}^p$, $\bz_{ij}$ are generated from  $N(0,\bSigma)$ with $\sigma_{kl}=0.5^{|k-l|}$, where $\bSigma=(\sigma_{kl})_{p\times p}$. The random errors $(\epsilon_{i1},\epsilon_{i2})^\T$ are generated from a two-dimensional normal distribution with mean zero and unit marginal compound symmetry covariance matrix with $\rho=0.7$.

Let $\bY_i=(y_{i1},y_{i2})^\T$ and $\bZ_i=(\bz_{i1}^\T,\bz_{i2}^\T)^\T$ respectively collect the response and predictor variables, and write $\bX_i=(\bY_i^\T,\bZ_i^\T)^\T$. To incorporate the dependence among the repeated measures from the same subject when estimating $\btheta_0$, we use the quadratic estimating equations proposed by \cite{Qu2000}:
    $$\bg(\bX_i;\btheta)=\left( \begin{array}{c}
    \bZ_i^\T \bv_i^{-1/2}\bM_1\bv_i^{-1/2}(\bY_i-\bZ_i^\T\btheta) \\
    \vdots \\
    \bZ_i^T \bv_i^{-1/2}\bM_m\bv_i^{-1/2}(\bY_i-\bZ_i^\T\btheta)
    \end{array} \right),
    $$
where  $\bv_i$ is a diagonal matrix of the conditional variances of subject $i$, and $\bM_j$ $(j=1,\dots,m)$ are working correlation matrices. Note that when $m=1$, i.e., using only one working correlation matrix $\bM_1$, the model becomes the one in \cite{LiangeZeger1986} and we have $r=p$. Here we choose two sets of basis matrices with $\bM_1$ being the identity matrix of size $n_i$ and $\bM_2$ being the compound symmetry with the diagonal elements of 1 and off-diagonal elements of $\rho$. In our setting, $n_i=2$ and therefore $r=2p$ estimating equations to estimate $p$ parameters. For each simulation, we repeat the experiment 100 times.

We obtain the same quantities as those in the example of Section \ref{sim2}, and report them in  Table \ref{longitudinal_table1}.
In comparison of the single-penalty method, we can conclude from Table \ref{longitudinal_table1}, with the columns defined in the same way as those in Table \ref{table_preg},  that the proposed double-penalty method has much better performance.  
This confirms the efficacy and efficiency of adding the additional penalty on the Lagrange multiplier $\blambda$, which performs the selection of estimating equations by reducing the number of estimating equations to less than 10.

\begin{table}
\begin{center}
\tabcolsep 0.1in \arrayrulewidth 1pt \doublerulesep 2pt
\begin{tabular}{ l l c c c c c c } \hline
$(n,p,r)$ & Method & $\btheta_{\rm nonzeros}$ & $\btheta_{\rm true}$ & ME & No. EE's \\ \hline
$(50,100,200)$ & MLE-Oracle & 3 (0) & NA & 0.023 (0.002) & NA \\
& MLE & 100 (0) & 3 (0) & 3.446 (0.106) & NA \\
& PEL-BIC & 0 (0) & 0 (0) & 15.25 (0) & 200 (0) \\
& PEL-BICC & 0 (0) & 0 (0) & 15.25 (0) & 200 (0) \\
& PEL-EBIC & 0 (0) & 0 (0) & 15.25 (0) & 200 (0) \\
& PEL2-BIC & 27.92 (2.51) & 2.95 (0.04) & 5.252 (0.871) & 5.29 (0.23) \\
& PEL2-BICC & 27.00 (2.69) & 2.95 (0.04) & 4.532 (0.552) & 5.21 (0.24) \\
& PEL2-EBIC & 24.80 (2.87) & 2.94 (0.04) & 4.657 (0.625) & 5.26 (0.25) \\
\hline
$(100,200,400)$ & MLE-Oracle & 3 (0) & NA & 0.014 (0.001) & NA \\
& MLE & 200 (0) & 3 (0) & 3.438 (0.068) & NA \\
& PEL-BIC & 0 (0) & 0 (0) & 15.25 (0) & 400 (0) \\
& PEL-BICC & 0 (0) & 0 (0) & 15.25 (0) & 400 (0) \\
& PEL-EBIC & 0 (0) & 0 (0) & 15.25 (0) & 400 (0) \\
& PEL2-BIC & 45.46 (4.37) & 3 (0) & 5.241 (0.793) & 5.51 (0.19) \\
& PEL2-BICC & 43.00 (4.25) & 2.99 (0.01) & 4.736 (0.659) & 5.50 (0.18) \\
& PEL2-EBIC & 42.40 (4.33) & 2.99 (0.01) & 4.546 (0.649) & 5.52 (0.19) \\
\hline
$(100,500,1000)$ & MLE-Oracle & 3 (0) & NA & 0.011 (0.001) & NA \\
& MLE & NA & NA & NA & NA \\
& PEL-BIC & 0 (0) & 0 (0) & 15.25 (0) & 1000 (0) \\
& PEL-BICC & 0 (0) & 0 (0) & 15.25 (0) & 1000 (0) \\
& PEL-EBIC & 0 (0) & 0 (0) & 15.25 (0) & 1000 (0) \\
& PEL2-BIC & 30.02 (6.11) & 2.93 (0.03) & 2.300 (0.359) & 6.70 (0.16) \\
& PEL2-BICC & 26.73 (6.02) & 2.93 (0.03) & 2.430 (0.377) & 6.62 (0.16) \\
& PEL2-EBIC & 25.09 (5.91) & 2.93 (0.03) & 2.415 (0.377) & 6.59 (0.16) \\
\hline
\end{tabular}
\caption{Simulation results for regression model for longitudinal data with repeated measures based on 100 replicates. Here $\btheta_{\rm nonzero}$ is the average number of selected nonzero components,  $\btheta_{\rm true}$ is the average number of true nonzero components that are selected, ME reports the model error, and No.EE's reports the number of estimating equations selected. \label{longitudinal_table1}}
\end{center}
\end{table}

\subsection{Trial of activity for adolescent girls 2 (TAAG2)}

We apply the penalized EL with two penalties  to examine the individual-, social-, and neighborhood-level factors associated with adolescent girls' physical activity over time in the Trial of Activity for Adolescent Girls 2 (TAAG2) \citep{Young2014,Grant2015}. The 589 girls in the Maryland site from TAAG2 were collected data at 8th grade (2009) and 11th (2011) grade. The response variable, moderate to vigorous physical activity (MVPA) minutes, were assessed from accelerometers. Forty-two variables to be considered include: (1) demographic and psychosocial information (individual- and social-level variables) that were obtained from questionnaires; (2) height, weight, and triceps skinfold to assess body composition; and (3) geographical information systems and self-report for neighborhood-level variables. There are 554 girls have complete information for all 42 variables and are used in this analysis.

A two-time point longitudinal linear mixed effects model is used to identify factors that are most relevant to MVPA. A similar model as in Section \ref{sim3} is used with two working correlation structure matrices. Our double-penalty EL method identifies four variables are related to MVPA: {\em Self-management strategies}, {\em Self-efficacy}, {\em Perceived barriers}, and {\em Social support}. In particular, higher {\em Self-management strategies}, {\em Self-efficacy}, {\em Social support} and lower {\em Perceived barriers} are associated with higher MVPA. Our finding confirms the previous results in \cite{Young2014,Grant2015}.


\setcounter{equation}{0}
\section{Discussion}
\label{s6}

We study a new penalized EL approach with two penalties, with one encouraging sparsity of the estimator and the other encouraging sparsity of the Lagrange multiplier in the optimizations associated with the EL.  Such an approach utilizes sparsity in the target parameters and effectively achieves a moment selection procedure for estimating the sparse parameter.   Both theory and numerical examples confirm the merits of the new penalized EL.

One interesting extension of the approach is to explore inferences with estimating equations after the variable selection procedure.  Such a direction is a suitable stage for EL method with estimating equations who takes advantage of adaptivity  to various moment conditions with less stringent distributional assumptions.  The other interesting and challenging  problem is to explore the optimality of the sparse estimator using estimating equations with high data dimensionality.  Semiparametric efficiency of EL with estimating equations is shown in \cite{QinLawless1994}.  However,  when the paradigm shifts to high-dimensional statistical problems, the efficiency of the sparse estimator respecting its nonzero components remains open for further investigations.  We plan to address the problems in future works.

\section*{Acknowledgments}
We are grateful to the Co-Editor, the Associate Editor and three
referees for very constructive comments and suggestions that  have greatly improved our
paper. Chang was supported in part by 
a grant from the Australian Research Council.
Tang acknowledges supports from NSF Grants IIS-1546087 and SES-1533956.
Wu's  research was partially supported by NIH grants R01HL094572 and R01HL119058.

\setcounter{equation}{0}
\section{Proofs}\label{s7}

In the sequel, we use the abbreviations ``w.p.a.1" and ``w.r.t" to denote, respectively, ``with probability approaching one" and ``with respect to", and $C$ denotes a generic positive finite constant that may be different in different uses. For simplicity and when no confusion arises, we use notation $\bh_i(\btheta)$ as equivalent to $\bh(\bX_i;\btheta)$ for a generic $q$-dimensional multivariate function $\bh(\cdot;\cdot)$ and denote by $h_{i,k}(\btheta)$ the $k$th component of $\bh_i(\btheta)$.  Let $\bar \bh(\btheta)=n^{-1}\sum_{i=1}^n \bh_i (\btheta)$, and $\bar h_k(\btheta)=n^{-1}\sum_{i=1}^nh_{i,k}(\btheta)$ be the $k$th component of $\bar \bh(\btheta)$. For a given set $\mathcal{L}\subset\{1,\ldots,q\}$, we denote by $\bh_{\mathcal{L}}(\cdot;\cdot)$ the subvector of $\bh(\cdot;\cdot)$ collecting the components indexed by $\mathcal{L}$. Analogously, we let $\bh_{i,\mathcal{L}}(\btheta)=\bh_{\mathcal{L}}(\bX_i;\btheta)$ and $\bar{\bh}_{\mathcal{L}}(\btheta)=n^{-1}\sum_{i=1}^n\bh_{i,\mathcal{L}}(\btheta)$. For an $s_1\times s_2$ matrix $\bB=(b_{ij})$, let $|\bB|_\infty=\max_{1\leq i\leq s_1,1\leq j\leq s_2}|b_{ij}|$, $\|\bB\|_1=\max_{1\leq j\leq s_2}\sum_{i=1}^{s_1}|b_{ij}|$, $\|\bB\|_\infty=\max_{1\leq i\leq s_1}\sum_{j=1}^{s_2}|b_{ij}|$ and $\|\bB\|_2=\lambda_{\max}^{1/2}(\bB\bB^\T)$ where $\lambda_{\max}(\bB\bB^\T)$ denotes the largest eigenvalue of $\bB\bB^\T$. Specifically, if $s_2=1$, we use $|\bB|_1=\sum_{i=1}^{s_1}|b_{i1}|$ and $|\bB|_2=(\sum_{i=1}^{s_1}b_{i1}^2)^{1/2}$ to denote the $L_1$-norm and $L_2$-norm of the $s_1$-dimensional vector $\bB$, respectively.

\subsection{Proof of Proposition \ref{pn:0}}\label{se:p1}
Define $A_n(\btheta,\blambda)=n^{-1}\sum_{i=1}^n\log\{1+\blambda^\T\bg_i(\btheta)\}$ for any $\btheta\in\bTheta$ and $\blambda\in\widehat{\Lambda}_n(\btheta)$. We first prove that $\max_{\blambda\in\widehat{\Lambda}_n(\btheta_0)}A_n(\btheta_0,\blambda)=O_p(rn^{-1})$. Let $\widetilde{\blambda}=\arg\max_{\blambda\in\widehat{\Lambda}_n(\btheta_0)}A_n(\btheta_0,\blambda)$. Pick $\delta_n=o(r^{-1/2}n^{-1/\gamma})$ and $r^{1/2}n^{-1/2}=o(\delta_n)$, which is guaranteed by $r^2n^{2/\gamma-1}=o(1)$. Let $\bar{\blambda}=\arg\max_{\blambda\in\Lambda_n}A_n(\btheta_0,\blambda)$ where $\Lambda_n=\{\blambda\in\mathbb{R}^r:|\blambda|_2\leq \delta_n\}$. It follows from Markov inequality that $\max_{1\leq i\leq n}|\bg_i(\btheta_0)|_2=O_p(r^{1/2}n^{1/\gamma})$. Then $\max_{1\leq i\leq n,\blambda\in\Lambda_n}|{\blambda}^\T\bg_i(\btheta_0)|=o_p(1)$. By Taylor expansion, it holds w.p.a.1 that
\begin{equation}\label{eq:0}
\begin{split}
0=A_n(\btheta_0,\bzero)\leq A_n(\btheta_0,\bar{\blambda})=&~\bar{\blambda}^\T\bar{\bg}(\btheta_0)-\frac{1}{2n}\sum_{i=1}^n\frac{\bar{\blambda}^\T\bg_i(\btheta_0)\bg_i(\btheta_0)^\T\bar{\blambda}}{\{1+c\bar{\blambda}^\T\bg_i(\btheta_0)\}^2}\\
\leq&~|\bar{\blambda}|_2|\bar{\bg}(\btheta_0)|_2-C|\bar{\blambda}|_2^2\{1+o_p(1)\},
\end{split}
\end{equation}
for some $|c|<1$. Notice that $|\bar{\bg}(\btheta_0)|_2=O_p(r^{1/2}n^{-1/2})$, (\ref{eq:0}) yields that $|\bar{\blambda}|_2=O_p(r^{1/2}n^{-1/2})=o_p(\delta_n)$. Therefore, $\bar{\blambda}\in\textrm{int}(\Lambda_n)$ w.p.a.1. Since $\Lambda_n\subset\widehat{\Lambda}_n(\btheta_0)$ w.p.a.1, $\widetilde{\blambda}=\bar{\blambda}$ w.p.a.1 by the concavity of $A_n(\btheta_0,\blambda)$ and $\widehat{\Lambda}_n(\btheta_0)$. Hence, by (\ref{eq:0}), we have $\max_{\blambda\in\widehat{\Lambda}_n(\btheta_0)}A_n(\btheta_0,\blambda)=O_p(rn^{-1})$.

We then show $|\bar{\bg}(\widehat{\btheta})|_2=O_p(r^{1/2}n^{-1/2})$. For $\delta_n$ specified above, let $\blambda^*=\delta_n\bar{\bg}(\widehat{\btheta})/|\bar{\bg}(\widehat{\btheta})|_2$, then $\blambda^*\in\Lambda_n$. By Taylor expansion, it holds w.p.a.1 that
\begin{equation}\label{eq:cterm1}
\begin{split}
A_n(\widehat{\btheta},\blambda^*)=&~\blambda^{*,\T}\bar{\bg}(\widehat{\btheta})-\frac{1}{2n}\sum_{i=1}^n\frac{{\blambda}^{\ast,\T}\bg_i(\widehat{\btheta})\bg_i(\widehat{\btheta})^\T{\blambda}^\ast}{\{1+c{\blambda}^{\ast,\T}\bg_i(\widehat{\btheta})\}^2}\\
\geq&~\delta_n|\bar{\bg}(\widehat{\btheta})|_2-C\delta_n^2\{1+o_p(1)\},
\end{split}
\end{equation}
for some $|c|<1$. Notice that $A_n(\widehat{\btheta},\blambda^*)\leq \max_{\blambda\in\widehat{\Lambda}_n(\widehat{\btheta})}A_n(\widehat{\btheta},\blambda)\leq \max_{\blambda\in\widehat{\Lambda}_n(\btheta_0)}A_n(\btheta_0,\blambda)=O_p(rn^{-1})$, thus $|\bar{\bg}(\widehat{\btheta})|_2=O_p(\delta_n)$. Consider any $\epsilon_n\rightarrow0$ and let $\blambda^{**}=\epsilon_n \bar{\bg}(\widehat{\btheta})$, then $|\blambda^{**}|_2=o_p(\delta_n)$. Using the same arguments above, we can obtain $ \epsilon_n|\bar{\bg}(\widehat{\btheta})|_2^2-C\epsilon_n^2|\bar{\bg}(\widehat{\btheta})|_2^2\{1+o_p(1)\}=O_p(rn^{-1}). $ Then $\epsilon_n|\bar{\bg}(\widehat{\btheta})|_2^2=O_p(rn^{-1})$. Notice that we can select arbitrary slow $\epsilon_n\rightarrow0$, following a standard result from probability theory, we have $|\bar{\bg}(\widehat{\btheta})|_2^2=O_p(rn^{-1})$. Hence, we complete the proof. $\hfill\Box$

\subsection{Proof of Proposition \ref{pn:00}}\label{se:pn00}

Define $F_n(\btheta)=\max_{\blambda\in\widehat{\Lambda}_n(\btheta)}A_n(\btheta,\blambda)+\sum_{k=1}^pP_{1,\pi}(|\theta_k|)$ where $\btheta=(\theta_1,\ldots,\theta_p)^\T$ and $A_n(\btheta,\blambda)=n^{-1}\sum_{i=1}^n\log\{1+\blambda^\T\bg_i(\btheta)\}$. Recall $a_n=\sum_{k=1}^pP_{1,\pi}(|\theta_k^0|)$ and $b_n=\max\{rn^{-1},a_n\}$. As shown in the proof of Proposition \ref{pn:0}, $\max_{\blambda\in\widehat{\Lambda}_n(\btheta_0)}A_n(\btheta_0,\blambda)=O_p(rn^{-1})$ which implies $F_n(\btheta_0)=O_p(rn^{-1})+a_n$. Define $\bTheta_*=\{\btheta=(\btheta_{\mathcal{S}}^\T,\btheta_{\mathcal{S}^c}^\T)^\T:|\btheta_{\mathcal{S}}-\btheta_{0,\mathcal{S}}|_\infty\leq\varepsilon, |\btheta_{\mathcal{S}^c}|_1\leq n^{-1/2}\varphi_n^{-1}\}$ for some fixed $\varepsilon>0$. Let $\widetilde{\btheta}_n=\arg\min_{\btheta\in\bTheta_*}F_n(\btheta)$. As ${F}_n(\widetilde{\btheta}_n)\leq {F}_n(\btheta_0)$, we have $F_n(\widetilde{\btheta}_n)\leq O_p(rn^{-1})+a_n=O_p(b_n)$.
We will first show that $\widetilde{\btheta}_n\in\textrm{int}(\bTheta_*)$ w.p.a.1. To do this, our proof includes two steps: (i) to show that for any $\epsilon_n\rightarrow\infty$ satisfying $b_n\epsilon_n^{2\beta}n^{2/\gamma}=o(1)$,  there exists a uniform constant $K>0$ independent of $\btheta$ such that $\mathbb{P}\{{F}_n(\btheta)>Kb_n\epsilon_n^{2\beta}\}\rightarrow1$ as $n\rightarrow\infty$ for any
$\btheta=(\btheta_{\mathcal{S}}^\T,\btheta_{\mathcal{S}^c}^\T)^\T\in\bTheta_*$ satisfying $|\btheta_{\mathcal{S}}-\btheta_{0,\mathcal{S}}|_\infty> \epsilon_nb_n^{1/(2\beta)}$. Thus $|\widetilde{\btheta}_{n,\mathcal{S}}-\btheta_{0,\mathcal{S}}|_\infty=O_p\{\epsilon_nb_n^{1/(2\beta)}\}$. Notice that we can select arbitrary slow diverging $\epsilon_n$, following a standard result from probability theory, we have $|\widetilde{\btheta}_{n,\mathcal{S}}-\btheta_{0,\mathcal{S}}|_\infty=O_p\{b_n^{1/(2\beta)}\}$, (ii) to show that $|\widetilde{\btheta}_{n,\mathcal{S}^c}|_1<n^{-1/2}\varphi_n^{-1}$.

For (i), we will use the technique developed for the proof of Theorem 1 in \cite{ChangTangWu2013}. For any ${\btheta}=(\btheta_{\mathcal{S}}^\T,\btheta_{\mathcal{S}^c}^\T)^\T\in\bTheta_*$ satisfying $|{\btheta}_{\mathcal{S}}-\btheta_{0,\mathcal{S}}|_\infty> \epsilon_nb_n^{1/(2\beta)}$, define ${\btheta}^*=(\btheta_{\mathcal{S}}^\T,\bzero^\T)^\T$ and let $j_0=\arg\max_{1\leq j\leq r}|\mathbb{E}\{g_{i,j}({\btheta}^*)\}|$. Define $\mu_{j_0}=\mathbb{E}\{g_{i,j_0}({\btheta})\}$, $\mu_{j_0}^*=\mathbb{E}\{g_{i,j_0}({\btheta}^*)\}$, and $\widetilde{\blambda}=\delta b_n^{1/2}\epsilon_n^{\beta}\bfe_{j_0}$ where $\delta>0$ is a constant to be determined later, and $\bfe_{j_0}$ is an $r$-dimensional vector with the $j_0$-th component being $1$ and other components being $0$. Without lose of generality, we assume $\mu_{j_0}^*>0$. (\ref{eq:c1}) and Markov inequality yield that $\max_{1\leq i\leq n}|g_{i,j_0}(\btheta)|=O_p(n^{1/\gamma})$, which implies $\max_{1\leq i\leq n}|\widetilde{\blambda}^\T\bg_i(\btheta)|=O_p(b_n^{1/2}\epsilon_n^{\beta}n^{1/\gamma})=o_p(1)$. Then $\widetilde{\blambda}\in\widehat{\Lambda}_n({\btheta})$ w.p.a.1. Write $\btheta=(\theta_1,\ldots,\theta_p)^\T$ and $\widetilde{\blambda}=(\widetilde{\lambda}_1,\ldots,\widetilde{\lambda}_r)^\T$. By the definition of ${F}_n(\btheta)$, it
holds w.p.a.1 that
\begin{equation*}\label{eq:equ2}
\begin{split}
{F}_n(\btheta)\geq&~\frac{1}{n}\sum_{i=1}^n\log\{1+\widetilde{\blambda}^\T\bg_i({\btheta})\}+\sum_{k=1}^pP_{1,\pi}(|{\theta}_k|)\\
\geq&~\frac{1}{n}\sum_{i=1}^n\widetilde{\lambda}_{j_0}g_{i,j_0}({\btheta})-\frac{1}{2n}\sum_{i=1}^n\frac{\{\widetilde{\lambda}_{j_0}g_{i,j_0}({\btheta})\}^2}{\{1+c\widetilde{\lambda}_{j_0}g_{i,j_0}({\btheta})\}^2}\\
\geq&~\frac{1}{n}\sum_{i=1}^n\widetilde{\lambda}_{j_0}g_{i,j_0}({\btheta})-\frac{1}{n}\sum_{i=1}^n\{\widetilde{\lambda}_{j_0}g_{i,j_0}({\btheta})\}^2\\
\end{split}
\end{equation*}
for some $|c|<1$ and $\widetilde{\lambda}_{j_0}=\delta b_n^{1/2}\epsilon_n^\beta$. Therefore, it holds that
\[
\begin{split}
&~\mathbb{P}\big\{{F}_n({\btheta})\leq Kb_n\epsilon_n^{2\beta}\big\}\\
\leq&~\mathbb{P} \bigg[\frac{1}{n}\sum_{i=1}^n\{g_{i,j_0}({\btheta})-\mu_{j_0}\}\leq b_n^{1/2}\epsilon_n^{\beta}\bigg\{\frac{K}{\delta}+\frac{\delta}{n}\sum_{i=1}^ng_{i,j_0}^2(\btheta)\bigg\}-\mu_{j_0}\bigg]+o(1).
\end{split}
\]
From (\ref{eq:c1}) and Markov inequality, there exists a uniform positive constant $L$ independent of $\btheta$ such that $\mathbb{P}\{n^{-1}\sum_{i=1}^ng_{i,j_0}^2(\btheta)>L\}\rightarrow0$. Thus, with $\delta=(K/L)^{1/2}$, we have
\[
\mathbb{P}\big\{{F}_n({\btheta})\leq Kb_n\epsilon_n^{2\beta}\big\}\leq\mathbb{P} \bigg[\frac{1}{n}\sum_{i=1}^n\{g_{i,j_0}({\btheta})-\mu_{j_0}\}\leq2b_n^{1/2}\epsilon_n^{\beta}(KL)^{1/2}-\mu_{j_0}\bigg]+o(1).
\]
From (\ref{eq:idenlocal}) and (\ref{eq:ide1}), we know that $\mu_{j_0}^*\geq \Delta(\epsilon_nb_n^{1/(2\beta)})\geq K_1\epsilon_n^\beta b_n^{1/2}/2$ with $K_1$ specified in (\ref{eq:idenlocal}) for sufficiently large $n$, and
 \[
 |\mu_{j_0}-\mu_{j_0}^*|\leq \sum_{k\notin\mathcal{S}}\mathbb{E}\bigg\{\sup_{\btheta\in\bTheta_*}\bigg|\frac{\partial g_{i,j_0}(\btheta)}{\partial \theta_k}\bigg|\bigg\}|\theta_k|\leq K_2|\btheta_{\mathcal{S}^c}|_1=o(b_n^{1/2})
 \]
 for $K_2$ specified in (\ref{eq:ide1}). Therefore, $\mu_{j_0}\geq K_1\epsilon_n^\beta b_n^{1/2}/3$ for sufficiently large $n$.
 For sufficiently small $K$ independent of $\btheta$, we have
$
2b_n^{1/2}\epsilon_n^{\beta}(KL)^{1/2}-\mu_{j_0}\leq-c\mu_{j_0}
$
for some $0<c<1$, which implies that
$
n^{1/2}\{2b_n^{1/2}\epsilon_n^{\beta}(KL)^{1/2}-\mu_{j_0}\}\leq -cn^{1/2}\mu_{j_0}\lesssim -\epsilon_n^\beta b_n^{1/2}n^{1/2}\rightarrow-\infty.
$
As $n^{-1/2}\sum_{i=1}^n\{g_{i,j_0}({\btheta})-\mu_{j_0}\}\xrightarrow{d}N(0,\sigma^2)$ for some $\sigma>0$, it holds that $\mathbb{P}\{{F}_n({\btheta})\leq Kb_n\epsilon_n^{2\beta}\}\rightarrow0$. Hence, we complete the proof for (i). 

For (ii), if $|\widetilde{\btheta}_{n,\mathcal{S}^c}|_1=n^{-1/2}\varphi_n^{-1}$, we define $\widetilde{\btheta}_n^*=(\widetilde{\btheta}_{n,\mathcal{S}}^\T,\tau\widetilde{\btheta}_{n,\mathcal{S}^c}^\T)^\T$ for some $\tau\in(0,1)$ and will show $F_n(\widetilde{\btheta}_n^*)<F_n(\widetilde{\btheta}_n)$ w.p.a.1. Notice that $\widetilde{\btheta}_n=\arg\min_{\btheta\in\bTheta_*}F_n(\btheta)$. This will be a contradiction. Therefore, $|\widetilde{\btheta}_{n,(2)}|_1<n^{-1/2}\varphi_n^{-1}$. Write $\widetilde{\btheta}_n=(\widetilde{\theta}_{n,1},\ldots,\widetilde{\theta}_{n,p})^\T$ and $\widetilde{\btheta}_n^*=(\widetilde{\theta}_{n,1}^*,\ldots,\widetilde{\theta}_{n,p}^*)^\T$. By the definition of ${F}_n(\btheta)$ and the inequality ${F}_n(\widetilde{\btheta}_n)\leq {F}_n(\btheta_0)$, it holds that
\[
\max_{\blambda\in\widehat{\Lambda}_n(\widetilde{\btheta}_n)}A_n(\widetilde{\btheta}_n,\blambda)\leq \max_{\blambda\in\widehat{\Lambda}_n(\btheta_0)}A_n(\btheta_0,\blambda)+\sum_{k=1}^pP_{1,\pi}(|\theta_k^0|)-\sum_{k=1}^pP_{1,\pi}(|\widetilde{\theta}_{n,k}|).
\]
On the other hand, it holds that
\begin{equation}\label{eq:pdiff1}
\begin{split}
\sum_{k=1}^pP_{1,\pi}(|\theta_k^0|)-\sum_{k=1}^pP_{1,\pi}(|\widetilde{\theta}_{n,k}|)\leq&~\sum_{k=1}^{s}P_{1,\pi}(|\theta_k^0|)-\sum_{k=1}^{s}P_{1,\pi}(|\widetilde{\theta}_{n,k}|)\\
\leq&~\sum_{k=1}^{s}P_{1,\pi}'\{c_k|\widetilde{\theta}_{n,k}|+(1-c_k)|\theta_k^0|\}|\widetilde{\theta}_{n,k}-\theta_k^0|\\
=&~O_p\{s\chi_nb_n^{1/(2\beta)}\}
\end{split}
\end{equation}
for some $c_k\in(0,1)$.
As we have shown in Section \ref{se:p1}, $\max_{\blambda\in\widehat{\Lambda}_n(\btheta_0)}A_n(\btheta_0,\blambda)=O_p(rn^{-1})$. Therefore, $\max_{\blambda\in\widehat{\Lambda}_n(\widetilde{\btheta}_n)}A_n(\widetilde{\btheta}_n,\blambda)=O_p(rn^{-1})+O_p\{s\chi_nb_n^{1/(2\beta)}\}$. Pick $\delta_n$ satisfying $\delta_n=o(r^{-1/2}n^{-1/\gamma})$ and $\max\{rn^{-1},s\chi_nb_n^{1/(2\beta)}\}=o(\delta_n^2)$, which can be guaranteed by $r^2n^{2/\gamma-1}=o(1)$ and $rs\chi_n b_n^{1/(2\beta)}n^{2/\gamma}=o(1)$. Same as (\ref{eq:cterm1}), we have
\begin{equation*}
\begin{split}
o_p(\delta_n^2)=\max_{\blambda\in\widehat{\Lambda}_n(\widetilde{\btheta}_n)}A_n(\widetilde{\btheta}_n,\blambda)\geq \delta_n|\bar{\bg}(\widetilde{\btheta}_n)|_2-C\delta_n^2\{1+o_p(1)\},
\end{split}
\end{equation*}
which implies $|\bar{\bg}(\widetilde{\btheta}_n)|_2=O_p(\delta_n)$. Following the same arguments in Section \ref{se:p1} below (\ref{eq:cterm1}), we have $|\bar{\bg}(\widetilde{\btheta}_n)|_2=O_p(r^{1/2}n^{-1/2})+O_p\{s^{1/2}\chi_n^{1/2}b_n^{1/(4\beta)}\}$. Notice that
$|\bar{\bg}(\widetilde{\btheta}_n^*)|_2\leq |\bar{\bg}(\widetilde{\btheta}_n)|_2+|\{\nabla_{\btheta}\bar{\bg}(\bar{\btheta})\}(\widetilde{\btheta}_n^*-\widetilde{\btheta}_n)|_2$ for some $\bar{\btheta}$ lying on the jointing line between $\widetilde{\btheta}_n$ and $\widetilde{\btheta}_n^*$. Since $\widetilde{\btheta}_{n,\mathcal{S}}=\widetilde{\btheta}_{n,\mathcal{S}}^*$, by (\ref{eq:p2}), it holds that $|\{\nabla_{\btheta}\bar{\bg}(\bar{\btheta})\}(\widetilde{\btheta}_n^*-\widetilde{\btheta}_n)|_2=O_p(r^{1/2}n^{-1/2})$. Hence, $|\bar{\bg}(\widetilde{\btheta}_n^*)|_2=O_p(r^{1/2}n^{-1/2})+O_p\{s^{1/2}\chi_n^{1/2}b_n^{1/(4\beta)}\}$. Write $\blambda^*=\arg\max_{\blambda\in\widehat{\Lambda}_n(\widetilde{\btheta}_n^*)}A_n(\widetilde{\btheta}_n^*,\blambda)$. Following the same arguments for (\ref{eq:0}), it holds that $|\blambda^*|_2=O_p(r^{1/2}n^{-1/2})+O_p\{s^{1/2}\chi_n^{1/2}b_n^{1/(4\beta)}\}$. Since $\widetilde{\btheta}_n^*=(\widetilde{\btheta}_{n,\mathcal{S}}^\T,\tau\widetilde{\btheta}_{n,\mathcal{S}^c}^\T)^\T$ and $F_n(\widetilde{\btheta}_n)\geq A_n(\widetilde{\btheta}_n,\blambda^*)+\sum_{k=1}^pP_{1,\pi}(|\widetilde{\theta}_{n,k}|)$, then
\begin{equation}\label{eq:diff21}
\begin{split}
F_n(\widetilde{\btheta}_n^*)=&~\frac{1}{n}\sum_{i=1}^n\log\{1+\blambda^{*,\T}\bg_i(\widetilde{\btheta}_n^*)\}+\sum_{k=1}^pP_{1,\pi}(|\widetilde{\theta}_{n,k}^*|)\\
=&~\frac{1}{n}\sum_{i=1}^n\log\{1+\blambda^{*,\T}\bg_i(\widetilde{\btheta}_n)\}+\bigg\{\frac{1}{n}\sum_{i=1}^n\frac{\blambda^{*,\T}\nabla_{\btheta}\bg_i(\check{\btheta})}{1+\blambda^{*,\T}\bg_i(\check{\btheta})}\bigg\}(\widetilde{\btheta}_n^*-\widetilde{\btheta}_n)+\sum_{k=1}^pP_{1,\pi}(|\widetilde{\theta}_{n,k}^*|)\\
\leq&~F_n(\widetilde{\btheta}_n)+\bigg\{\frac{1}{n}\sum_{i=1}^n\frac{\blambda^{*,\T}\nabla_{\btheta}\bg_i(\check{\btheta})}{1+\blambda^{*,\T}\bg_i(\check{\btheta})}\bigg\}(\widetilde{\btheta}_n^*-\widetilde{\btheta}_n)+\sum_{k=s+1}^pP_{1,\pi}(\tau|\widetilde{\theta}_{n,k}|)-\sum_{k=s+1}^pP_{1,\pi}(|\widetilde{\theta}_{n,k}|),
\end{split}
\end{equation}
for some $\check{\btheta}$ lying on the jointing line between $\widetilde{\btheta}_n$ and $\widetilde{\btheta}_n^*$. Notice that $\max_{1\leq i\leq n}|\blambda^{*,\T}\bg_i(\check{\btheta})|=o_p(1)$, then
\[
\begin{split}
\bigg|\bigg\{\frac{1}{n}\sum_{i=1}^n\frac{\blambda^{*,\T}\nabla_{\btheta}\bg_i(\check{\btheta})}{1+\blambda^{*,\T}\bg_i(\check{\btheta})}\bigg\}(\widetilde{\btheta}_n^*-\widetilde{\btheta}_n)\bigg|\leq&~ |\blambda^*|_2\bigg|\bigg\{\frac{1}{n}\sum_{i=1}^n\frac{\nabla_{\btheta}\bg_i(\check{\btheta})}{1+\blambda^{*,\T}\bg_i(\check{\btheta})}\bigg\}(\widetilde{\btheta}_n^*-\widetilde{\btheta}_n)\bigg|_2\\
\leq&~|\blambda^*|_2|\widetilde{\btheta}_{n,\mathcal{S}^c}|_1O_p(r^{1/2}\varphi_n).
\end{split}
\]
On the other hand,
\[
\begin{split}
\sum_{k=s+1}^pP_{1,\pi}(\tau|\widetilde{\theta}_{n,k}|)-\sum_{k=s+1}^pP_{1,\pi}(|\widetilde{\theta}_{n,k}|)=&-(1-\tau)\sum_{k=s+1}^pP_{1,\pi}'\{(c_k\tau+1-c_k)|\widetilde{\theta}_{n,k}|\}|\widetilde{\theta}_{n,k}|\\
\leq&-(1-\tau)C\pi\sum_{k=s+1}^p|\widetilde{\theta}_{n,k}|=-(1-\tau)C\pi|\widetilde{\btheta}_{n,\mathcal{S}^c}|_1
\end{split}
\]
for some $c_k\in(0,1)$. If $r^{1/2}\varphi_n\max\{r^{1/2}n^{-1/2},s^{1/2}\chi_n^{1/2}b_n^{1/(4\beta)}\}=o(\pi)$, (\ref{eq:diff21}) implies $F_n(\widetilde{\btheta}_n^*)<F_n(\widetilde{\btheta}_n)$ w.p.a.1. Hence, we complete the proof of (ii).

Nextly, we will show $\mathbb{P}(\widetilde{\btheta}_{n,\mathcal{S}^c}=\bzero)\rightarrow1$. Define
\begin{equation*}\label{eq:G}
\widehat{G}_n(\btheta,\blambda)=\frac{1}{n}\sum_{i=1}^n\log\{1+\blambda^\T\bg_i({\btheta})\}+\sum_{k=1}^pP_{1,\pi}(|\theta_k|)
\end{equation*}
for $\btheta=(\theta_1,\ldots,\theta_p)^\T$. Then $\widetilde{\btheta}_n$ and its Lagrange multiplier $\widehat{\blambda}$ satisfy the score equation $\nabla_{\blambda}\widehat{G}_n(\widetilde{\btheta}_n,\widehat{\blambda})=\bzero$. By the implicit theorem [Theorem 9.28 of \cite{Rudin1976}], for all $\btheta$ in a $|\cdot|_2$-neighborhood of $\widetilde{\btheta}_n$, there is a $\widehat{\blambda}(\btheta)$ such that $\nabla_{\blambda}\widehat{G}_n\{\btheta,\widehat{\blambda}(\btheta)\}=\bzero$ and $\widehat{\blambda}(\btheta)$ is continuously differentiable in $\btheta$. By the concavity of $\widehat{G}_n(\btheta,\blambda)$ w.r.t $\blambda$, $\widehat{G}_n\{\btheta,\widehat{\blambda}(\btheta)\}=\max_{\blambda\in\widehat{\Lambda}_n(\btheta)}\widehat{G}_n(\btheta,\blambda)$. Write $\widehat{\blambda}=(\widehat{\lambda}_1,\ldots,\widehat{\lambda}_r)^\T$. From the envelope theorem,
\begin{equation*}
\begin{split}
\bzero=\nabla_{\btheta}\widehat{G}_n\{\btheta,\widehat{\blambda}(\btheta)\}\Big|_{\btheta=\widetilde{\btheta}_n}=\frac{1}{n}\sum_{i=1}^n\frac{\{\nabla_{\btheta}\bg_i(\widetilde{\btheta}_n)\}^\T \widehat{\blambda}}{1+\widehat{\blambda}^\T\bg_i(\widetilde{\btheta}_n)}+\bigg\{\sum_{k=1}^p\nabla_{\btheta}P_{1,\pi}(|\theta_k|)\bigg\}\bigg|_{\btheta=\widetilde{\btheta}_n}.
\end{split}
\end{equation*}
Write $\widehat{\bh}=(\widehat{h}_1,\ldots,\widehat{h}_p)^\T=\nabla_{\btheta}\widehat{G}_n\{\btheta,\widehat{\blambda}(\btheta)\}|_{\btheta=\widetilde{\btheta}_n}$. 
Let $\rho_1(t;\pi)=\pi^{-1}P_{1,\pi}(t)$. Since $P_{1,\pi}(\cdot)\in\mathcal{P}$, $\rho_1'(0^+;\pi)$ is independent of $\pi$. We write it as $\rho_1'(0^+)$ for simplicity. Therefore, for each $k=1,\ldots,p$,
\[
\widehat{h}_k=\frac{1}{n}\sum_{i=1}^n\sum_{j=1}^r\frac{\widehat{\lambda}_j}{1+\widehat{\blambda}^\T\bg_i(\widetilde{\btheta}_n)}\frac{\partial g_{i,j}(\widetilde{\btheta}_n)}{\partial \theta_k}+\widehat{\kappa}_k,
\]
where $\widehat{\kappa}_k=\pi \rho_{1}'(|\widetilde{\theta}_k|;\pi)\textrm{sgn}(\widetilde{\theta}_k)$ for $\widetilde{\theta}_k\neq0$ and $\widehat{\kappa}_k\in[-\pi \rho_{1}'(0^+),\pi \rho_{1}'(0^+)]$ otherwise. From Triangle inequality, it holds that
\[
\begin{split}
\sup_{k\notin \mathcal{S}}\bigg|\frac{1}{n}\sum_{i=1}^n\sum_{j=1}^r\frac{\widehat{\lambda}_j}{1+\widehat{\blambda}^\T\bg_i(\widetilde{\btheta}_n)}\frac{\partial g_{i,j}(\widetilde{\btheta}_n)}{\partial \theta_k}\bigg|\leq&~\bigg[\sum_{j=1}^r|\widehat{\lambda}_j|\sup_{k\notin \mathcal{S}}\bigg\{\frac{1}{n}\sum_{i=1}^n\bigg|\frac{\partial g_{i,j}(\widetilde{\btheta}_n)}{\partial \theta_k}\bigg|\bigg\}\bigg]\{1+o_p(1)\}\\
\leq&~ O_p(\varphi_n)\cdot\sum_{j=1}^r|\widehat{\lambda}_j|\\
=&~O_p\big(r^{1/2}\varphi_n\max\{r^{1/2}n^{-1/2},s^{1/2}\chi_n^{1/2}b_n^{1/(4\beta)}\}\big).
\end{split}
\]
As $r^{1/2}\varphi_n\max\{r^{1/2}n^{-1/2},s^{1/2}\chi_n^{1/2}b_n^{1/(4\beta)}\}=o(\pi)$, if $\widetilde{\theta}_k\neq0$ for some $k\notin\mathcal{S}$, then $\pi \rho_{1}'(|\widetilde{\theta}_k|;\pi)\textrm{sgn}(\widetilde{\theta}_k)$ will dominates the sign of $\widehat{h}_k$. According to the arguments for the proof of Lemma 1 in \cite{FanLi2001}, we know $\widetilde{\btheta}_{n,\mathcal{S}^c}=\bzero$ w.p.a.1. Hence, we complete the proof of Proposition \ref{pn:00}. $\hfill\Box$

\subsection{Proof of Proposition \ref{pn:2}}

Recall $\mathcal{M}_{\btheta_n}=\{1\leq j\leq r:|\bar{g}_j(\btheta_n)|\geq\nu\rho_2'(0^+)\}$ and $\mathcal{M}_{\btheta_n}^*=\{1\leq j\leq r:|\bar{g}_j(\btheta_n)|\geq C\nu\rho_2'(0^+)\}$ for some $C\in(0,1)$. Clearly, $\mathcal{M}_{\btheta_n}\subset\mathcal{M}_{\btheta_n}^*$. Recall $m_n=|\mathcal{M}_{\btheta_n}^*|$. Given $\mathcal{M}_{\btheta_n}$, we select $\delta_n$ satisfying $\delta_n=o(m_n^{-1/2}n^{-1/\gamma})$ and $u_n=o(\delta_n)$. Let $\bar{\blambda}_n=\arg\max_{\blambda\in\Lambda_n}f(\blambda;\btheta_n)$ where $\Lambda_{n}=\{\blambda=(\blambda_{\mathcal{M}_{\btheta_n}}^\T,\blambda_{\mathcal{M}_{\btheta_n}^c}^\T)^\T\in\mathbb{R}^r:|\blambda_{\mathcal{M}_{\btheta_n}}|_2\leq \delta_n~\textrm{and}~\blambda_{\mathcal{M}_{\btheta_n}^c}=\bzero\}$. For given $\mathcal{M}_{\btheta_n}$, Condition \ref{as:moment} and Markov inequality imply that $\max_{1\leq i\leq n}|\bg_{i,\mathcal{M}_{\btheta_n}}(\btheta_n)|_2=O_p(m_n^{1/2}n^{1/\gamma})$, which leads to $ \max_{1\leq i\leq n}|\bar{\blambda}_n^\T\bg_i(\btheta_n)|=o_p(1)$. Write $\bar{\blambda}_n=(\bar{\lambda}_{n,1},\ldots,\bar{\lambda}_{n,r})^\T$. By the definition of $\bar{\blambda}_n$ and Taylor expansion, noting $P_{2,\nu}(t)=\nu\rho_2(t;\nu)$ and $\rho_2'(t;\nu)\geq \rho_2'(0^+)$ for any $t>0$, we have
\begin{equation*}\label{eq:equ1}
\begin{split}
0=&~f(\bzero;\btheta_n)\leq f(\bar{\blambda}_n;\btheta_n)\\
=&~\frac{1}{n}\sum_{i=1}^n\bar{\blambda}_n^\T\bg_i(\btheta_n)-\frac{1}{2n}\sum_{i=1}^n\frac{\bar{\blambda}_n^\T\bg_i(\btheta_0)\bg_i(\btheta_n)^\T\bar{\blambda}_n}{\{1+c\bar{\blambda}_n^\T\bg_i(\btheta_n)\}^2}-\sum_{j=1}^rP_{2,\nu}(|\bar{\lambda}_{n,j}|)\\
\leq&~\bar{\blambda}_{n,\mathcal{M}_{\btheta_n}}^\T\{\bar{\bg}_{\mathcal{M}_{\btheta_n}}(\btheta_n)-\nu\rho_2'(0^+)\sgn(\bar{\blambda}_{n,\mathcal{M}_{\btheta_n}})\}-\frac{1}{2}\lambda_{\min}\{\widehat{\bV}_{\mathcal{M}_{\btheta_n}}(\btheta_n)\}|\bar{\blambda}_{n,\mathcal{M}_{\btheta_n}}|_2^2\{1+o_p(1)\}\\
\leq&~\bar{\blambda}_{n,\mathcal{M}_{\btheta_n}}^\T[\bar{\bg}_{\mathcal{M}_{\btheta_n}}(\btheta_n)-\nu\rho_2'(0^+)\sgn\{\bar{\bg}_{\mathcal{M}_{\btheta_n}}(\btheta_n)\}]-\frac{1}{2}\lambda_{\min}\{\widehat{\bV}_{\mathcal{M}_{\btheta_n}}(\btheta_n)\}|\bar{\blambda}_{n,\mathcal{M}_{\btheta_n}}|_2^2\{1+o_p(1)\}\\
\end{split}
\end{equation*}
 Notice that $|\bar{\bg}_{\mathcal{M}_{\btheta_n}}(\btheta_n)-\nu\rho_2'(0^+)\sgn\{\bar{\bg}_{\mathcal{M}_{\btheta_n}}(\btheta_n)\}|_2=O_p(u_n)$ and $\mathbb{P}[\lambda_{\min}\{\widehat{\bV}_{\mathcal{M}_{\btheta_n}}(\btheta_n)\}\geq C]\rightarrow1$, then $|\bar{\blambda}_{n,\mathcal{M}_{\btheta_n}}|_2=O_p(u_n)=o_p(\delta_n)$. Write $\bar{\blambda}_{n,\mathcal{M}_{\btheta_n}}=(\bar{\lambda}_1,\ldots,\bar{\lambda}_{|\mathcal{M}_{\btheta_n}|})^\T$. We have w.p.a.1 that
 \begin{equation}\label{eq:firstorder}
 \bzero=\frac{1}{n}\sum_{i=1}^n\frac{\bg_{i,\mathcal{M}_{\btheta_n}}(\btheta_n)}{1+\bar{\blambda}_{n,\mathcal{M}_{\btheta_n}}^\T\bg_{i,\mathcal{M}_{\btheta_n}}(\btheta_n)}-\widehat{\bfeta}
 \end{equation}
 where $\widehat{\bfeta}=(\widehat{\eta}_1,\ldots,\widehat{\eta}_{|\mathcal{M}_{\btheta_n}|})^\T$ with $\widehat{\eta}_j=\nu\rho_2'(|\bar{\lambda}_j|;\nu)\textrm{sgn}(\bar{\lambda}_j)$ for $\bar{\lambda}_j\neq0$ and $\widehat{\eta}_j\in[-\nu\rho_2'(0^+),\nu\rho_2'(0^+)]$ for $\bar{\lambda}_j=0$. (\ref{eq:firstorder}) implies that
$
\widehat{\bfeta}=\bar{\bg}_{\mathcal{M}_{\btheta_n}}(\btheta_n)+\bR
$ with $|\bR|_\infty=O_p(\varsigma_n^{1/2}u_n)$. Since $\varsigma_n^{1/2}u_n=o(\nu)$,
  then w.p.a.1 $\textrm{sgn}(\bar{\lambda}_j)=\textrm{sgn}\{\bar{g}_{j}(\btheta_n)\}$ for any $\bar{\lambda}_j\neq0$. 

 We will show that $\bar{\blambda}_n$ is a local maximizer for $f(\blambda;\btheta_n)$ w.p.a.1. We first show that $\bar{\blambda}_n=\arg\max_{\blambda\in\Lambda_n^*(\btheta_n)}f(\blambda;\btheta_n)$ w.p.a.1, where $\Lambda_n^*(\btheta_n)=\{\blambda=(\blambda_{\mathcal{M}_{\btheta_n}^*}^{\T},\blambda_{\mathcal{M}_{\btheta_n}^{*,c}}^\T)^\T\in\mathbb{R}^r:|\blambda_{\mathcal{M}_{\btheta_n}^*}|_2\leq \epsilon, \blambda_{\mathcal{M}_{\btheta_n}^{*,c}}=\bzero\}$ for some $\epsilon>0$. Notice that $f(\blambda;\btheta_n)$ is concave w.r.t $\blambda$. To do this, it suffices to show that $\bw=\bar{\blambda}_{n,\mathcal{M}_{\btheta_n}^*}^\T=:(w_1,\ldots,w_{m_n})^\T\in\mathbb{R}^{m_n}$ satisfies the equation
 \[
 \bzero=\frac{1}{n}\sum_{i=1}^n\frac{\bg_{i,\mathcal{M}_{\btheta_n}^*}(\btheta_n)}{1+\bw^\T\bg_{i,\mathcal{M}_{\btheta_n}^*}(\btheta_n)}-\widehat{\bfeta}^*
 \]
 w.p.a.1, where $\widehat{\bfeta}^*=(\widehat{\eta}_1^*,\ldots,\widehat{\eta}_{m_n}^*)^\T$ with $\widehat{\eta}_j^*=\nu\rho_2'(|w_j|;\nu)\textrm{sgn}(w_j)$ for $w_j\neq0$ and $\widehat{\eta}_j^*\in[-\nu\rho_2'(0^+),\nu\rho_2'(0^+)]$ for $w_j=0$.
Based on (\ref{eq:firstorder}), we know
$
 0=n^{-1}\sum_{i=1}^ng_{i,j}(\btheta_n)/\{1+\bw^\T\bg_{i,\mathcal{M}_{\btheta_n}^*}(\btheta_n)\}-\widehat{\eta}_j^*
$
 holds for any $j\in\mathcal{M}_{\btheta_n}$. For each $j\in\mathcal{M}_{\btheta_n}^*\backslash\mathcal{M}_{\btheta_n}$, it holds that
$
 n^{-1}\sum_{i=1}^ng_{i,j}(\btheta_n)/\{1+\bw^\T\bg_{i,\mathcal{M}_{\btheta_n}^*}(\btheta_n)\}=\bar{g}_j(\btheta_n)+O_p(\varsigma_n^{1/2}u_n)
$
 where $O_p(\varsigma_n^{1/2}u_n)$ is uniform for any $j\in\mathcal{M}_{\btheta_n}^*\backslash\mathcal{M}_{\btheta_n}$. Since $C\nu\rho_2'(0^+)\leq|\bar{g}_j(\btheta_n)|<\nu\rho_2'(0^+)$ for $j\in\mathcal{M}_{\btheta_n}^*\backslash\mathcal{M}_{\btheta_n}$, if $\varsigma_n^{1/2}u_n=o(\nu)$, then $
 |n^{-1}\sum_{i=1}^ng_{i,j}(\btheta_n)/\{1+\bw^\T\bg_{i,\mathcal{M}_{\btheta_n}^*}(\btheta_n)\}|<\nu\rho_2'(0^+)
$ w.p.a.1 for any $j\in\mathcal{M}_{\btheta_n}^*\backslash\mathcal{M}_{\btheta_n}$.
  This implies that there exists $\widehat{\eta}_j^*$ such that $
 0=n^{-1}\sum_{i=1}^ng_{i,j}(\btheta_n)/\{1+\bw^\T\bg_{i,\mathcal{M}_{\btheta_n}^*}(\btheta_n)\}-\widehat{\eta}_j^*
$ holds for any $j\in\mathcal{M}_{\btheta_n}^*\backslash\mathcal{M}_{\btheta_n}$.

 Secondly, we prove $\bar{\blambda}_n$
 is a local maximizer for $f(\blambda;\btheta_n)$ over $\blambda\in\widetilde{\Lambda}_n(\btheta_n)$ w.p.a.1, where $\widetilde{\Lambda}_n(\btheta_n)=\{\blambda=(\blambda_{\mathcal{M}_{\btheta_n}^*}^\T,\blambda_{\mathcal{M}_{\btheta_n}^{*,c}}^\T)^\T\in\mathbb{R}^r:|\blambda_{\mathcal{M}_{\btheta_n}^*}-\bar{\blambda}_{n,\mathcal{M}_{\btheta_n}^*}|_2\leq o(u_n),|\blambda_{\mathcal{M}_{\btheta_n}^{*,c}}|_1=o(r^{-1/\gamma}n^{-1/\gamma})\}$. Notice that
$\max_{1\leq i\leq n,\blambda\in\widetilde{\Lambda}_n(\btheta_n)}|{\blambda}^\T\bg_i(\btheta_n)|=o_p(1)$. For any $\blambda\in\widetilde{\Lambda}_n(\btheta_n)$, we write $\blambda=(\blambda_{\mathcal{M}_{\btheta_n}^*}^\T,\blambda_{\mathcal{M}_{\btheta_n}^{*,c}}^\T)^\T$ and denote by $\widetilde{\blambda}=(\blambda_{\mathcal{M}_{\btheta_n}^*}^{\T},\bzero^\T)^\T$ the projection of $\blambda$ onto the subspace $\Lambda_n^*(\btheta_n)$. We only need to show \begin{equation}\label{eq:toshow}\mathbb{P}\bigg[\sup_{\blambda\in\widetilde{\Lambda}_n(\btheta_n)}\{f(\blambda;\btheta_n)-f(\widetilde{\blambda};\btheta_n)\}\leq0\bigg]\rightarrow1.\end{equation}
By Taylor expansion, it holds that
\begin{equation*}
\sup_{\blambda\in\widetilde{\Lambda}_n(\btheta_n)}\{f(\blambda;\btheta_n)-f(\widetilde{\blambda};\btheta_n)\}=\sup_{\blambda\in\widetilde{\Lambda}_n(\btheta_n)}\bigg\{\frac{1}{n}\sum_{i=1}^n\frac{\bg_i(\btheta_n)^\T(\blambda-\widetilde{\blambda})}{1+\blambda_*^\T\bg_i(\btheta_n)}-\sum_{j\in\mathcal{M}_{\btheta_n}^{*,c}}P_{2,\nu}(|\lambda_j|)\bigg\},
\end{equation*}
for some $\blambda_*$ lying on the jointing line between $\blambda$ and $\widetilde{\blambda}$. 
We have that \[
\bigg|\frac{1}{n}\sum_{i=1}^n\frac{\bg_i(\btheta_n)^\T(\blambda-\widetilde{\blambda})}{1+\blambda_*^\T\bg_i(\btheta_n)}\bigg|\leq C\nu\rho_2'(0^+)\sum_{j\in\mathcal{M}_{\btheta_n}^{*,c}}|\lambda_j|+O_p(m_n^{1/2}u_n\varsigma_n)\cdot\sum_{j\in\mathcal{M}_{\btheta_n}^{*,c}}|\lambda_j|.
\]
where the term $O_p(m_n^{1/2}u_n\varsigma_n)$ is uniformly for any $\blambda\in\widetilde{\Lambda}_n(\btheta_n)$. On the other hand,  we have
\[
\sum_{j\in\mathcal{M}_{\btheta_n}^{*,c}}P_{2,\nu}(|\lambda_j|)\geq \nu\rho_2'(0^+)\sum_{j\in\mathcal{M}_{\btheta_n}^{*,c}}|\lambda_j|.
\]
Hence,
\[
\begin{split}
&\frac{1}{n}\sum_{i=1}^n\frac{\bg_i(\btheta_n)^\T(\blambda-\widetilde{\blambda})}{1+\blambda_*^\T\bg_i(\btheta_n)}-\sum_{j\in\mathcal{M}_{\btheta_n}^{*,c}}P_{2,\nu}(|\lambda_j|)\leq \bigg\{-(1-C)\nu\rho_2'(0^+)+O_p(m_n^{1/2}u_n\varsigma_n)\bigg\}\sum_{j\in\mathcal{M}_{\btheta_n}^{*,c}}|\lambda_j|.
\end{split}
\]
Notice that $m_n^{1/2}u_n\varsigma_n/\nu\rightarrow 0$, then $-(1-C)\nu\rho_2'(0^+)+O_p(m_n^{1/2}u_n\varsigma_n)\leq 0$ w.p.a.1 which implies (\ref{eq:toshow}) holds. Hence, $\bar{\blambda}_n$ w.p.a.1 is a local maximizer of $f(\blambda;\btheta_n)$. We complete the proof of Proposition \ref{pn:2}. $\hfill\Box$

\subsection{Proof of Theorem \ref{tm:consistency}}\label{se:1}

Let $\mathcal{G}_0=\textrm{supp}\{\widehat{\blambda}(\btheta_0)\}$. It holds that
\[
\begin{split}
\max_{\blambda\in\widehat{\Lambda}_n(\btheta_0)}f(\blambda;\btheta_0)=&~\max_{\bfeta\in\widehat{\Lambda}_n^\dag(\btheta_0)}\bigg[\frac{1}{n}\sum_{i=1}^n\log\{1+\bfeta^\T\bg_{i,\mathcal{G}_0}(\btheta_0)\}-\sum_{j=1}^{|\mathcal{G}_0|}P_{2,\nu}(|\eta_j|)\bigg]\\
\leq&~\max_{\bfeta\in\widehat{\Lambda}_n^\dag(\btheta_0)}\frac{1}{n}\sum_{i=1}^n\log\{1+\bfeta^\T\bg_{i,\mathcal{G}_0}(\btheta_0)\},
\end{split}
\]
where $\widehat{\Lambda}_n^\dag(\btheta_0)=\{\bfeta\in\mathbb{R}^{m_0}:\bfeta^\T\bg_{i,\mathcal{G}_0}(\btheta_0)\in\mathcal{V},i=1,\ldots,n\}$ for some open interval $\mathcal{V}$ containing zero. Given $\mathcal{G}_{0}$, since $|\mathcal{G}_0|\leq \ell_n$, following the proof of Proposition \ref{pn:0}, we have $\max_{\bfeta\in\widehat{\Lambda}_n^\dag(\btheta_0)}n^{-1}\sum_{i=1}^n\log\{1+\bfeta^\T\bg_{i,\mathcal{G}_0}(\btheta_0)\}=O_p(\ell_nn^{-1})$ which implies
$\max_{\blambda\in\widehat{\Lambda}_n(\btheta_0)}f(\blambda;\btheta_0)=O_p(\ell_nn^{-1})$.

Recall $a_n=\sum_{k=1}^pP_{1,\pi}(|\theta_k^0|)$, $b_n=\max\{\ell_nn^{-1},a_n,\nu^2\}$ and $S_n(\btheta)=\max_{\blambda\in\widehat{\Lambda}_n(\btheta)}f(\blambda;\btheta)+\sum_{k=1}^pP_{1,\pi}(|\theta_k|)$ for any $\btheta=(\theta_1,\ldots,\theta_p)^\T$. Define $\bTheta_*=\{\btheta=(\btheta_{\mathcal{S}}^\T,\btheta_{\mathcal{S}^c}^\T)^\T:|\btheta_{\mathcal{S}}-\btheta_{0,\mathcal{S}}|_\infty\leq\varepsilon, |\btheta_{\mathcal{S}^c}|_1\leq \aleph_n\}$ for some fixed $\varepsilon>0$ and $\aleph_n=\min\{s\omega_n^{1/2}b_n^{1/(2\beta)}\xi_n^{-1/2},o(b_n^{1/2}),o(\nu\varrho_n^{-1/2}\ell_n^{-3/2}\xi_n^{-1/2})\}$. Let $\widehat{\btheta}_n=\arg\min_{\btheta\in\bTheta_*} S_n(\btheta)$. As we have shown above, $\mathbb{P}\{{S}_n(\btheta_0)\leq a_n+O_p(\ell_nn^{-1})\}\rightarrow1$ as $n\rightarrow\infty$. As ${S}_n(\widehat{\btheta}_n)\leq {S}_n(\btheta_0)$, we have $\mathbb{P}\{{S}_n(\widehat{\btheta}_n)\leq a_n+O_p(\ell_nn^{-1})\}\rightarrow1$ as $n\rightarrow\infty$. We will show that $\widehat{\btheta}_n\in\textrm{int}(\bTheta_*)$ w.p.a.1. Same as the proof of Proposition \ref{pn:00} stated in Section \ref{se:pn00}, our proof includes two steps: (i) to show that for any $\epsilon_n\rightarrow\infty$ satisfying $b_n\epsilon_n^{2\beta}n^{2/\gamma}=o(1)$,  there exists a uniform constant $K>0$ independent of $\btheta$ such that $\mathbb{P}\{{S}_n(\btheta)>Kb_n\epsilon_n^{2\beta}\}\rightarrow1$ as $n\rightarrow\infty$ for any
$\btheta=(\btheta_{\mathcal{S}}^\T,\btheta_{\mathcal{S}^c}^\T)^\T\in\bTheta_*$ satisfying $|\btheta_{\mathcal{S}}-\btheta_{0,\mathcal{S}}|_\infty> \epsilon_nb_n^{1/(2\beta)}$, which leads to $|\widehat{\btheta}_{n,\mathcal{S}}-\btheta_{0,\mathcal{S}}|_\infty=O_p\{b_n^{1/(2\beta)}\}$. (ii) to show that $|\widehat{\btheta}_{n,\mathcal{S}^c}|_1<\aleph_n$. The proof of (i) is the same as that stated in Section \ref{se:pn00}, thus we omit its proof and only show (ii) here. We need the following lemma whose proof is given in the supplementary material.
\begin{lemma}\label{la:covc}
Let $\mathscr{F}=\{\mathcal{F}\subset\{1,\ldots,r\}:|\mathcal{F}|\leq\ell_n\}$ and $\bTheta_n=\{\btheta=(\btheta_{\mathcal{S}}^\T,\btheta_{\mathcal{S}^c}^\T)^\T:|\btheta_{\mathcal{S}}-\btheta_{0,\mathcal{S}}|_\infty=O_p\{b_n^{1/(2\beta)}\},|\btheta_{\mathcal{S}^c}|_1\leq \aleph_n\}$. Assume that Conditions {\rm \ref{as:smalleigenvalue}} and {\rm\ref{as:pa1}}, then
\[
\sup_{\btheta\in\bTheta_n}\sup_{\mathcal{F}\in\mathscr{F}}\|\widehat{\bV}_{\mathcal{F}}(\btheta)-\bV_{\mathcal{F}}(\btheta_0)\|_2=O_p\{s(\ell_n\omega_nb_n^{1/\beta})^{1/2}\}+O_p\{\ell_n(n^{-1}\varrho_n\log r)^{1/2}\}
\]
provided that $\log r=o(n^{1/3})$, $s^2\ell_n\omega_nb_n^{1/\beta}=o(1)$ and $\ell_n^2n^{-1}\varrho_n\log r=o(1)$.
\end{lemma}

We begin to prove (ii) now. If $|\widehat{\btheta}_{n,\mathcal{S}^c}|_1=\aleph_n$, we define $\widehat{\btheta}_n^*=(\widehat{\btheta}_{n,\mathcal{S}}^\T,\tau\widehat{\btheta}_{n,\mathcal{S}^c}^\T)^\T$ for some $\tau\in(0,1)$ and will show $S_n(\widehat{\btheta}_n^*)<S_n(\widehat{\btheta}_n)$ w.p.a.1. Notice that $\widehat{\btheta}_n=\arg\min_{\btheta\in\bTheta_*}S_n(\btheta)$. This will be a contradiction. Therefore, $|\widehat{\btheta}_{n,\mathcal{S}^c}|_1<\aleph_n$. Write $\widehat{\btheta}_n=(\widehat{\theta}_{n,1},\ldots,\widehat{\theta}_{n,p})^\T$. Notice that
\[
\max_{\blambda\in\widehat{\Lambda}_n(\widehat{\btheta}_n)}f(\blambda;\widehat{\btheta}_n)\leq \max_{\blambda\in\widehat{\Lambda}_n(\btheta_0)}f(\blambda;\btheta_0)+\sum_{k=1}^pP_{1,\pi}(|\theta_k^0|)-\sum_{k=1}^pP_{1,\pi}(|\widehat{\theta}_{n,k}|),
\]
by (\ref{eq:pdiff1}), we have $\max_{\blambda\in\widehat{\Lambda}_n(\widehat{\btheta}_n)}f(\blambda;\widehat{\btheta}_n)=O_p(\ell_nn^{-1})+O_p\{s\chi_n b_n^{1/(2\beta)}\}$. Pick $\delta_n$ satisfying $\delta_n=o(\ell_n^{-1/2}n^{-1/\gamma})$ and $\max\{\ell_nn^{-1},s\chi_nb_n^{1/(2\beta)}\}=o(\delta_n^2)$, which can be guaranteed by $\ell_n s\chi_nb_n^{1/(2\beta)}n^{2/\gamma}=o(1)$ and $\ell_n^2n^{2/\gamma-1}=o(1)$. Select $\blambda^*$ such that $\blambda^*_{\mathcal{M}_{\widehat{\btheta}_n}}=\delta_n[\bar{\bg}_{\mathcal{M}_{\widehat{\btheta}_n}}(\widehat{\btheta}_n)-\nu\rho_2'(0^+)\sgn\{\bar{\bg}_{\mathcal{M}_{\widehat{\btheta}_n}}(\widehat{\btheta}_n)\}]/|\bar{\bg}_{\mathcal{M}_{\widehat{\btheta}_n}}(\widehat{\btheta}_n)-\nu\rho_2'(0^+)\sgn\{\bar{\bg}_{\mathcal{M}_{\widehat{\btheta}_n}}(\widehat{\btheta}_n)\}|_2$ and $\blambda^*_{\mathcal{M}_{\widehat{\btheta}_n}^c}=\bzero$. Write $\blambda^*=(\lambda_1^*,\ldots,\lambda_r^*)^\T$. Then
\[
\begin{split}
o_p(\delta_n^2)=&~\max_{\blambda\in\widehat{\Lambda}_n(\widehat{\btheta}_n)}f(\blambda;\widehat{\btheta}_n)\\
\geq&~\frac{1}{n}\sum_{i=1}^n\log\{1+\blambda^{*,\T}\bg_i(\widehat{\btheta}_n)\}-\sum_{j=1}^rP_{2,\nu}(|\lambda_j^*|)\\
=&~\blambda^{*,\T}_{\mathcal{M}_{\widehat{\btheta}_n}}\bar{\bg}_{\mathcal{M}_{\widehat{\btheta}_n}}(\widehat{\btheta}_n)-\frac{1}{2n}\sum_{i=1}^n\frac{\blambda^{*,\T}_{\mathcal{M}_{\widehat{\btheta}_n}}\bg_{i,\mathcal{M}_{\widehat{\btheta}_n}}(\widehat{\btheta}_n)\bg_{i,\mathcal{M}_{\widehat{\btheta}_n}}(\widehat{\btheta}_n)^\T\blambda^{*}_{\mathcal{M}_{\widehat{\btheta}_n}}}{\{1+c\blambda^{*,\T}_{\mathcal{M}_{\widehat{\btheta}_n}}\bg_{i,\mathcal{M}_{\widehat{\btheta}_n}}(\widehat{\btheta}_n)\}^2}-\sum_{j\in\mathcal{M}_{\widehat{\btheta}_n}}P_{2,\nu}(|\lambda_j^*|)\\
\geq&~\blambda^{*,\T}_{\mathcal{M}_{\widehat{\btheta}_n}}\bar{\bg}_{\mathcal{M}_{\widehat{\btheta}_n}}(\widehat{\btheta}_n)-C\delta_n^2\{1+o_p(1)\}-\nu\sum_{j\in\mathcal{M}_{\widehat{\btheta}_n}}\rho_2'(c_j|\lambda_j^*|;\nu)|\lambda_j^*|\\
=&~\blambda^{*,\T}_{\mathcal{M}_{\widehat{\btheta}_n}}\bar{\bg}_{\mathcal{M}_{\widehat{\btheta}_n}}(\widehat{\btheta}_n)-\nu\rho_2'(0^+)\sum_{j\in\mathcal{M}_{\widehat{\btheta}_n}}|\lambda_j^*|-C\delta_n^2\{1+o_p(1)\}-\nu\sum_{j\in\mathcal{M}_{\widehat{\btheta}_n}}c_j\rho_2''(c_j^*|\lambda_j^*|;\nu)|\lambda_j^*|^2\\
\geq&~\blambda^{*,\T}_{\mathcal{M}_{\widehat{\btheta}_n}}\{\bar{\bg}_{\mathcal{M}_{\widehat{\btheta}_n}}(\widehat{\btheta}_n)-\nu\rho_2'(0^+)\sgn(\blambda_{\mathcal{M}_{\widehat{\btheta}_n}}^*)\}-C\delta_n^2\{1+o_p(1)\}
\end{split}
\]
for some $c, c_j, c_j^*\in(0,1)$. Recall $\mathcal{M}_{\widehat{\btheta}_n}=\{1\leq j\leq r:|\bar{g}_j(\widehat{\btheta}_n)|\geq \nu\rho_2'(0^+)\}$, then $\sgn(\blambda_{\mathcal{M}_{\widehat{\btheta}_n}}^*)=\sgn\{\bar{\bg}_{\mathcal{M}_{\widehat{\btheta}_n}}(\widehat{\btheta}_n)\}$. Thus $|\bar{\bg}_{\mathcal{M}_{\widehat{\btheta}_n}}(\widehat{\btheta}_n)-\nu\rho_2'(0^+)\sgn\{\bar{\bg}_{\mathcal{M}_{\widehat{\btheta}_n}}(\widehat{\btheta}_n)\}|_2=O_p(\delta_n)$. Using the technique developed in Section \ref{se:p1}, we have $|\bar{\bg}_{\mathcal{M}_{\widehat{\btheta}_n}}(\widehat{\btheta}_n)-\nu\rho_2'(0^+)\sgn\{\bar{\bg}_{\mathcal{M}_{\widehat{\btheta}_n}}(\widehat{\btheta}_n)\}|_2=O_p(\ell_n^{1/2}n^{-1/2})+O_p\{s^{1/2}\chi_n^{1/2}b_n^{1/(4\beta)}\}$.

By Lemma \ref{la:covc} and Condition \ref{as:smalleigenvalue}, we know $\lambda_{\min}\{\widehat{\bV}_{\mathcal{M}_{\widehat{\btheta}_n}}(\widehat{\btheta}_n)\}\geq C$ w.p.a.1. Therefore Proposition \ref{pn:2} leads to $|\widehat{\blambda}(\widehat{\btheta}_n)|_2=O_p(\ell_n^{1/2}n^{-1/2})+O_p\{s^{1/2}\chi_n^{1/2}b_n^{1/(4\beta)}\}$. Based on this property of the Lagrange multiplier $\widehat{\blambda}(\widehat{\btheta}_n)$, we can follow the same arguments stated in Section \ref{se:pn00} to construct (ii). Specifically, write $\widehat{\blambda}(\widehat{\btheta}_n)$ and $\widehat{\blambda}(\widehat{\btheta}_n^*)$ as $\widehat{\blambda}=(\widehat{\lambda}_1,\ldots,\widehat{\lambda}_r)^\T$ and $\widehat{\blambda}^*=(\widehat{\lambda}_1^*,\ldots,\widehat{\lambda}_r^*)^\T$, respectively. In the sequel, we use $\check{\btheta}$ to denote a generic vector lying on the jointing line between $\widehat{\btheta}_n$ and $\widehat{\btheta}_n^*$ that may be different in different uses. Write $\widehat{\btheta}_n^*=(\widehat{\theta}_{n,1}^*,\ldots,\widehat{\theta}_{n,p}^*)^\T$. By Taylor expansion, it holds that
\begin{equation}\label{eq:diff2}
\begin{split}
S_n(\widehat{\btheta}_n^*)=&~\frac{1}{n}\sum_{i=1}^n\log\{1+\widehat{\blambda}^{*,\T}\bg_i(\widehat{\btheta}_n^*)\}-\sum_{j=1}^rP_{2,\nu}(|\widehat{\lambda}_j^*|)+\sum_{k=1}^pP_{1,\pi}(|\widehat{\theta}_{n,k}^*|)\\
\leq&~S_n(\widehat{\btheta}_n)+\underbrace{\sum_{j=1}^rP_{2,\nu}(|\widehat{\lambda}_j|)-\sum_{j=1}^rP_{2,\nu}(|\widehat{\lambda}_j^*|)}_{\textrm{I}}+\underbrace{\frac{1}{n}\sum_{i=1}^n\frac{\widehat{\blambda}^{*,\T}\nabla_{\btheta_{\mathcal{S}^c}}\bg_i(\check{\btheta})}{1+\widehat{\blambda}^{*,\T}\bg_i(\check{\btheta})}(\widehat{\btheta}_{n,\mathcal{S}^c}^*-\widehat{\btheta}_{n,\mathcal{S}^c})}_{\textrm{II}}\\
&+\underbrace{\sum_{k=s+1}^pP_{1,\pi}(\tau|\widehat{\theta}_{n,k}|)-\sum_{k=s+1}^pP_{1,\pi}(|\widehat{\theta}_{n,k}|)}_{\textrm{III}}.
\end{split}
\end{equation}
 We will show $\textrm{I}+\textrm{II}+\textrm{III}<0$ w.p.a.1 as follows.

For $\textrm{I}$, we will first specify the convergence rate of $|\widehat{\blambda}^*-\widehat{\blambda}|_1$. Define
\begin{equation}\label{eq:H}
\widehat{H}_n(\btheta,\blambda)=\frac{1}{n}\sum_{i=1}^n\log\{1+\blambda^\T\bg_i({\btheta})\}+\sum_{k=1}^pP_{1,\pi}(|\theta_k|)-\sum_{j=1}^rP_{2,\nu}(|\lambda_j|)
\end{equation}
for any $\btheta=(\theta_1,\ldots,\theta_p)^\T$ and $\blambda=(\lambda_1,\ldots,\lambda_r)^\T$. Then $\widehat{\btheta}_n$ and its Lagrange multiplier $\widehat{\blambda}$ satisfy the score equation $\nabla_{\blambda}\widehat{H}_n(\widehat{\btheta}_n,\widehat{\blambda})=\bzero$, i.e.
 \begin{equation}\label{eq:ke11}
\bzero=\frac{1}{n}\sum_{i=1}^n\frac{\bg_i(\widehat{\btheta}_n)}{1+\widehat{\blambda}^\T\bg_i(\widehat{\btheta}_n)}-\widehat{\bfeta},
\end{equation}
where $\widehat{\bfeta}=(\widehat{\eta}_1,\ldots,\widehat{\eta}_r)^\T$ with $\widehat{\eta}_j=\nu \rho_{2}'(|\widehat{\lambda}_j|;\nu)\textrm{sgn}(\widehat{\lambda}_j)$ for $\widehat{\lambda}_j\neq0$ and $\widehat{\eta}_j\in[-\nu \rho_{2}'(0^+),\nu \rho_{2}'(0^+)]$ for $\widehat{\lambda}_j=0$. Let $\mathcal{R}_n=\textrm{supp}\{\widehat{\blambda}(\widehat{\btheta}_n)\}$.
 Restricted on $\mathcal{R}_n$, for any $\btheta\in\mathbb{R}^p$ and $\bzeta=(\zeta_1,\ldots,\zeta_{|\mathcal{R}_n|})^\T\in\mathbb{R}^{|\mathcal{R}_n|}$ with each $\zeta_j\neq0$, define
 \[
 \bm(\bzeta,\btheta)=\frac{1}{n}\sum_{i=1}^n\frac{\bg_{i,\mathcal{R}_n}(\btheta)}{1+\bzeta^\T\bg_{i,\mathcal{R}_n}(\btheta)}-\bw,
 \]
where $\bw=(w_1,\ldots,w_{|\mathcal{R}_n|})^\T$ with $w_j=\nu\rho_{2}'(|\zeta_j|;\nu)\textrm{sgn}(\zeta_j)$. Then, $\widehat{\blambda}_{\mathcal{R}_n}$ and $\widehat{\btheta}_n$ satisfy $\bm(\widehat{\blambda}_{\mathcal{R}_n},\widehat{\btheta}_n)=\bzero$. By the implicit theorem [Theorem 9.28 of \cite{Rudin1976}], for all $\btheta$ in a $|\cdot|_2$-neighborhood of $\widehat{\btheta}_n$, there is a $\bzeta(\btheta)$ such that $\bm\{\bzeta(\btheta),\btheta\}=\bzero$ and $\bzeta(\btheta)$ is continuously differentiable in $\btheta$. Since $\widehat{\btheta}_{n,\mathcal{S}}^*=\widehat{\btheta}_{n,\mathcal{S}}$, we have
\[
\begin{split}
|\bzeta(\widehat{\btheta}_n^*)-\widehat{\blambda}_{\mathcal{R}_n}|_1=\big|\{\nabla_{\btheta}\bzeta(\btheta)|_{\btheta=\check{\btheta}}\}(\widehat{\btheta}_n^*-\widehat{\btheta}_n)\big|_1\leq\big\|\nabla_{\btheta_{\mathcal{S}^c}}\bzeta(\btheta)|_{\btheta=\check{\btheta}}\big\|_1|\widehat{\btheta}_{n,\mathcal{S}^c}^*-\widehat{\btheta}_{n,\mathcal{S}^c}|_1.
\end{split}
\]
 Notice that
\[
\begin{split}
\nabla_{\btheta_{\mathcal{S}^c}} \bzeta(\btheta)\big|_{\btheta=\check{\btheta}}=&-(\nabla_{\bzeta}\bm)^{-1}(\nabla_{\btheta_{\mathcal{S}^c}}\bm)\big|_{\btheta=\check{\btheta}}\\
=&\bigg(\frac{1}{n}\sum_{i=1}^n\frac{\bg_{i,\mathcal{R}_n}(\check{\btheta})\bg_{i,\mathcal{R}_n}(\check{\btheta})^\T}{\{1+\bzeta(\check{\btheta})^\T\bg_{i,\mathcal{R}_n}(\check{\btheta})\}^2}+\nu\textrm{diag}[\rho''_2\{|\zeta_1(\check{\btheta})|;\nu\},\ldots,\rho''_2\{|\zeta_{|\mathcal{R}_n|}(\check{\btheta})|;\nu\}]\bigg)^{-1}\\
&~~~~~~~~\times\bigg\{\frac{1}{n}\sum_{i=1}^n\frac{\nabla_{\btheta_{\mathcal{S}^c}}\bg_{i,\mathcal{R}_n}(\check{\btheta})}{1+\bzeta(\check{\btheta})^\T\bg_{i,\mathcal{R}_n}(\check{\btheta})}-\frac{1}{n}\sum_{i=1}^n\frac{\bg_{i,\mathcal{R}_n}(\check{\btheta})\bzeta(\check{\btheta})^\T\nabla_{\btheta_{\mathcal{S}^c}}\bg_{i,\mathcal{R}_n}(\check{\btheta})}{\{1+\bzeta(\check{\btheta})^\T\bg_{i,\mathcal{R}_n}(\check{\btheta})\}^2}\bigg\}\\
=&:\bA(\check{\btheta})\times\bB(\check{\btheta}).
\end{split}
\]
Since $\max_{1\leq i\leq n}|\bzeta(\check{\btheta})^\T\bg_{i,\mathcal{R}_n}(\check{\btheta})|=o_p(1)$, from Lemma \ref{la:covc}, we know $\|\bA(\check{\btheta})\|_1\leq |\mathcal{R}_n|^{1/2}\|\bA(\check{\btheta})\|_2=O_p(\ell_n^{1/2})$. Meanwhile, we have $|\bB(\check{\btheta})|_\infty=O_p(\xi_n^{1/2})$ which implies $\|\bB(\check{\btheta})\|_1=O_p(\xi_n^{1/2}\ell_n)$. Therefore, it holds that $\|\nabla_{\btheta_{\mathcal{S}^c}} \bzeta(\btheta)|_{\btheta=\check{\btheta}}\|_1\leq \|\bA(\check{\btheta})\|_1\|\bB(\check{\btheta})\|_1=O_p(\ell_n^{3/2}\xi_n^{1/2})$, which implies $|\bzeta(\widehat{\btheta}_n^*)-\widehat{\blambda}_{\mathcal{R}_n}|_1=O_p(\ell_n^{3/2}\xi_n^{1/2})|\widehat{\btheta}_{n,\mathcal{S}^c}|_1$. Let $\widetilde{\blambda}$ satisfy $\widetilde{\blambda}_{\mathcal{R}_n}=\bzeta(\widehat{\btheta}_n^*)$ and $\widetilde{\blambda}_{\mathcal{R}_n^c}=\bzero$. For any $j\in\mathcal{R}_n^c$, we have
\[
\begin{split}
&~\frac{1}{n}\sum_{i=1}^n\frac{g_{i,j}(\widehat{\btheta}_n^*)}{1+\widetilde{\blambda}^\T\bg_i(\widehat{\btheta}_n^*)}\\
=&~\frac{1}{n}\sum_{i=1}^n\frac{g_{i,j}(\widehat{\btheta}_n)}{1+\widetilde{\blambda}^\T\bg_i(\widehat{\btheta}_n)}+\bigg[\frac{1}{n}\sum_{i=1}^n\frac{\nabla_{\btheta_{\mathcal{S}^c}}g_{i,j}(\check{\btheta})}{1+\widetilde{\blambda}^\T\bg_i(\check{\btheta}_n)}-\frac{1}{n}\sum_{i=1}^n\frac{g_{i,j}(\check{\btheta})\widetilde{\blambda}^\T\nabla_{\btheta_{\mathcal{S}^c}}\bg_i(\check{\btheta})}{\{1+\widetilde{\blambda}^\T\bg_i(\check{\btheta}_n)\}^2}\bigg](\widehat{\btheta}_{n,\mathcal{S}^c}^*-\widehat{\btheta}_{n,\mathcal{S}^c})\\
=&~\frac{1}{n}\sum_{i=1}^n\frac{g_{i,j}(\widehat{\btheta}_n)}{1+\widehat{\blambda}^\T\bg_i(\widehat{\btheta}_n)}-\bigg[\frac{1}{n}\sum_{i=1}^n\frac{g_{i,j}(\widehat{\btheta}_n)\bg_i(\widehat{\btheta}_n)^\T}{\{1+\check{\blambda}^\T\bg_i(\widehat{\btheta}_n)\}^2}\bigg](\widetilde{\blambda}-\widehat{\blambda})+O_p(\xi_n^{1/2})|\widehat{\btheta}_{n,\mathcal{S}^c}^*-\widehat{\btheta}_{n,\mathcal{S}^c}|_1\\
=&~\frac{1}{n}\sum_{i=1}^n\frac{g_{i,j}(\widehat{\btheta}_n)}{1+\widehat{\blambda}^\T\bg_i(\widehat{\btheta}_n)}+O_p(\varrho_n^{1/2})|\widetilde{\blambda}-\widehat{\blambda}|_1+O_p(\xi_n^{1/2})|\widehat{\btheta}_{n,\mathcal{S}^c}^*-\widehat{\btheta}_{n,\mathcal{S}^c}|_1\\
=&~\frac{1}{n}\sum_{i=1}^n\frac{g_{i,j}(\widehat{\btheta}_n)}{1+\widehat{\blambda}^\T\bg_i(\widehat{\btheta}_n)}+o_p(\nu),
\end{split}
\]
where the term $o_p(\nu)$ holds uniformly for any $j\in\mathcal{R}_n^c$. Write $\widetilde{\blambda}=(\widetilde{\lambda}_1,\ldots,\widetilde{\lambda}_r)^\T$.
Recall that $\bzeta(\widehat{\btheta}_n^*)$ and $\widehat{\btheta}_n^*$ satisfy $\bm\{\bzeta(\widehat{\btheta}_n^*),\widehat{\btheta}_n^*\}=\bzero$, and (\ref{eq:ke11}) holds, then it holds w.p.a.1 that
\[
\bzero=\frac{1}{n}\sum_{i=1}^n\frac{\bg_i(\widehat{\btheta}_n^*)}{1+\widetilde{\blambda}^\T\bg_i(\widehat{\btheta}_n^*)}-\widehat{\bfeta}^*
\]
for $\widehat{\bfeta}^*=(\widehat{\eta}_1^*,\ldots,\widehat{\eta}_r^*)^\T$ with $\widehat{\eta}_j^*=\nu \rho_{2}'(|\widetilde{\lambda}_j|;\nu)\textrm{sgn}(\widetilde{\lambda}_j)$ for $\widetilde{\lambda}_j\neq0$ and $\widehat{\eta}_j^*\in[-\nu \rho_{2}'(0^+),\nu \rho_{2}'(0^+)]$ for $\widetilde{\lambda}_j=0$. By the concavity of $f(\blambda;\btheta)=n^{-1}\sum_{i=1}^n\log\{1+\blambda^\T\bg_i(\btheta)\}-\sum_{j=1}^rP_{2,\nu}(|\lambda_j|)$, we know $\widehat{\blambda}^*=\widetilde{\blambda}$ w.p.a.1. Hence, $|\widehat{\blambda}^*-\widehat{\blambda}|_1=O_p(\ell_n^{3/2}\xi_n^{1/2})|\widehat{\btheta}_{n,\mathcal{S}^c}|_1$. This implies $\textrm{I}=O_p(\ell_n^{3/2}\xi_n^{1/2}\nu)|\widehat{\btheta}_{n,\mathcal{S}^c}|_1$.

Let $\mathcal{J}_*=\textrm{supp}(\widehat{\blambda}^*)$. Notice that $\max_{1\leq i\leq n}|\widehat{\blambda}^{*,\T}\bg_i(\check{\btheta})|=o_p(1)$, then
\[
\begin{split}
|\textrm{II}|\leq |\widehat{\blambda}^*|_2\bigg|\bigg\{\frac{1}{n}\sum_{i=1}^n\frac{\nabla_{\btheta_{\mathcal{S}^c}}\bg_{i,\mathcal{J}_*}(\check{\btheta})}{1+\widehat{\blambda}^{*,\T}\bg_i(\check{\btheta})}\bigg\}(\widehat{\btheta}_{n,\mathcal{S}^c}^*-\widehat{\btheta}_{n,\mathcal{S}^c})\bigg|_2\leq|\widehat{\blambda}^*|_2|\widehat{\btheta}_{n,\mathcal{S}^c}|_1O_p(\ell_n^{1/2}\xi_n^{1/2}),
\end{split}
\]
which implies $\textrm{II}=\max\{\ell_n^{1/2}n^{-1/2},s^{1/2}\chi_n^{1/2}b_n^{1/(4\beta)}\}|\widehat{\btheta}_{n,\mathcal{S}^c}|_1O_p(\ell_n^{1/2}\xi_n^{1/2})$. On the other hand, by Taylor expansion, we have
\[
\textrm{III}=-(1-\tau)\sum_{k=s+1}^pP_{1,\pi}'\{(c_k\tau+1-c_k)|\widehat{\theta}_{n,k}|\}|\widehat{\theta}_{n,k}|\leq-(1-\tau)C\pi|\widehat{\btheta}_{n,\mathcal{S}^c}|_1
\]
for some $c_k\in(0,1)$. Since $\max\{\ell_n^{3/2}\xi_n^{1/2}\nu,\ell_n\xi_n^{1/2}n^{-1/2},\ell_n^{1/2}\xi_n^{1/2}s^{1/2}\chi_n^{1/2}b_n^{1/(4\beta)}\}=o(\pi)$, (\ref{eq:diff2}) implies $S_n(\widehat{\btheta}_n^*)<S_n(\widehat{\btheta}_n)$ w.p.a.1. Hence, we complete the proof of (ii). Together with (i), we know such defined $\widehat{\btheta}_n$ is a local minimizer of $S_n(\btheta)$. Following the same arguments stated in Section \ref{se:pn00}, we can prove $\mathbb{P}(\widehat{\btheta}_{n,\mathcal{S}^c}=\bzero)\rightarrow1$. We complete the proof of Theorem \ref{tm:consistency}. $\hfill\Box$

\subsection {Proof of Theorem \ref{tm:2}}

Recall $\mathcal{R}_n=\textrm{supp}\{\widehat{\blambda}(\widehat{\btheta}_n)\}$. We still write $\widehat{\blambda}=\widehat{\blambda}(\widehat{\btheta}_n)=(\widehat{\lambda}_1,\ldots,\widehat{\lambda}_r)^\T$. For $\widehat{H}_n(\btheta,\blambda)$ defined in (\ref{eq:H}), we have $\nabla_{\blambda}\widehat{H}_n(\widehat{\btheta}_n,\widehat{\blambda})=\bzero$, i.e. 
\begin{equation}\label{eq:ke1}
\bzero=\frac{1}{n}\sum_{i=1}^n\frac{\bg_i(\widehat{\btheta}_n)}{1+\widehat{\blambda}^\T\bg_i(\widehat{\btheta}_n)}-\widehat{\bfeta},
\end{equation}
where $\widehat{\bfeta}=(\widehat{\eta}_1,\ldots,\widehat{\eta}_r)^\T$ with $\widehat{\eta}_j=\nu \rho_{2}'(|\widehat{\lambda}_j|;\nu)\textrm{sgn}(\widehat{\lambda}_j)$ for $\widehat{\lambda}_j\neq0$ and $\widehat{\eta}_j\in[-\nu \rho_{2}'(0^+),\nu \rho_{2}'(0^+)]$ for $\widehat{\lambda}_j=0$. By Taylor expansion, we have
\begin{equation*}\label{eq:expan1}
\bzero=\frac{1}{n}\sum_{i=1}^n\bg_{i,\mathcal{R}_n}(\widehat{\btheta}_n)-\frac{1}{n}\sum_{i=1}^n\frac{\bg_{i,\mathcal{R}_n}(\widehat{\btheta}_n)\bg_{i,\mathcal{R}_n}(\widehat{\btheta}_n)^\T\widehat{\blambda}_{\mathcal{R}_n}}{\{1+c\widehat{\blambda}_{\mathcal{R}_n}^\T\bg_{i,\mathcal{R}_n}(\widehat{\btheta}_n)\}^2}-\widehat{\bfeta}_{\mathcal{R}_n},
\end{equation*}
for some $|c|<1$, which implies
\[
\begin{split}
\widehat{\blambda}_{\mathcal{R}_n}=&~\bigg[\frac{1}{n}\sum_{i=1}^n\frac{\bg_{i,\mathcal{R}_n}(\widehat{\btheta}_n)\bg_{i,\mathcal{R}_n}(\widehat{\btheta}_n)^\T}{\{1+c\widehat{\blambda}_{\mathcal{R}_n}^\T\bg_{i,\mathcal{R}_n}(\widehat{\btheta}_n)\}^2}\bigg]^{-1}\{\bar{\bg}_{\mathcal{R}_n}(\widehat{\btheta}_n)-\widehat{\bfeta}_{\mathcal{R}_n}\}.
\end{split}
\]
On the other hand, together with
\[
\bzero=\nabla_{\btheta}\widehat{H}_n(\btheta,\widehat{\blambda}(\btheta))\big|_{\btheta=\widehat{\btheta}_n}=\bigg\{\frac{1}{n}\sum_{i=1}^n\frac{\nabla_{\btheta}\bg_i(\widehat{\btheta}_n)}{1+\widehat{\blambda}^\T\bg_i(\widehat{\btheta}_n)}\bigg\}^\T\widehat{\blambda}+\bigg\{\sum_{k=1}^p\nabla_{\btheta}P_{1,\pi}(|\theta_k|)\bigg\}\bigg|_{\btheta=\widehat{\btheta}_n},
\]
it holds that
\begin{equation}\label{eq:expan3}
\begin{split}
\bzero=\bigg\{\frac{1}{n}\sum_{i=1}^n\frac{\nabla_{\btheta_{\mathcal{S}}}\bg_{i,\mathcal{R}_n}(\widehat{\btheta}_n)}{1+\widehat{\blambda}_{\mathcal{R}_n}^\T\bg_{i,\mathcal{R}_n}(\widehat{\btheta}_n)}\bigg\}^\T\bigg[\frac{1}{n}\sum_{i=1}^n\frac{\bg_{i,\mathcal{R}_n}(\widehat{\btheta}_n)\bg_{i,\mathcal{R}_n}(\widehat{\btheta}_n)^\T}{\{1+c\widehat{\blambda}_{\mathcal{R}_n}^\T\bg_{i,\mathcal{R}_n}(\widehat{\btheta}_n)\}^2}\bigg]^{-1}\{\bar{\bg}_{\mathcal{R}_n}(\widehat{\btheta}_n)-\widehat{\bfeta}_{\mathcal{R}_n}\}+\widehat{\bkappa}_{\mathcal{S}},
\end{split}
\end{equation}
where $\widehat{\bkappa}_{\mathcal{S}}=\{\sum_{k=1}^p\nabla_{\btheta_{\mathcal{S}}}P_{1,\pi}(|\theta_k|)\}|_{\btheta_{\mathcal{S}}=\widehat{\btheta}_{n,\mathcal{S}}}$. From Condition \ref{as:paa1}, it holds that
$
|\widehat{\bkappa}_{\mathcal{S}}|_\infty=O_p(\chi_n).
$
We will use (\ref{eq:expan3}) to derive the limiting distribution of $\widehat{\btheta}_{n,\mathcal{S}}$. Before this, we need the following lemmas.

\begin{lemma}\label{la:1}
Assume the conditions of Theorem {\rm\ref{tm:consistency}} hold. Then
\[
\begin{split}
\bigg\|\frac{1}{n}\sum_{i=1}^n\frac{\bg_{i,\mathcal{R}_n}(\widehat{\btheta}_n)\bg_{i,\mathcal{R}_n}(\widehat{\btheta}_n)^\T}{\{1+c\widehat{\blambda}_{\mathcal{R}_n}^\T\bg_{i,\mathcal{R}_n}(\widehat{\btheta}_n)\}^2}-\widehat{\bV}_{\mathcal{R}_n}(\btheta_0)\bigg\|_2=&~O_p(\ell_nn^{-1/2+1/\gamma})+O_p\{\ell_n^{1/2}s^{1/2}\chi_n^{1/2}b_n^{1/(4\beta)}n^{1/\gamma}\},
\end{split}
\]
and
\[
\begin{split}
\bigg|\bigg\{\frac{1}{n}\sum_{i=1}^n\frac{\nabla_{\btheta_{\mathcal{S}}}\bg_{i,\mathcal{R}_n}(\widehat{\btheta}_n)}{1+\widehat{\blambda}^\T_{\mathcal{R}_n}\bg_{i,\mathcal{R}_n}(\widehat{\btheta}_n)}-\nabla_{\btheta_{\mathcal{S}}}\bar{\bg}_{\mathcal{R}_n}(\widehat{\btheta}_n)\bigg\}\bz\bigg|_2=|\bz|_2\big[O_p(\ell_ns^{1/2}\omega_n^{1/2}n^{-1/2})+O_p\{\ell_n^{1/2}s\omega_n^{1/2}\chi_n^{1/2}b_n^{1/(4\beta)}\}\big]
\end{split}
\]
holds uniformly for any $\bz\in\mathbb{R}^{s}$.

\end{lemma}

\begin{lemma}\label{la:2}
Assume the conditions of Theorem {\rm\ref{tm:consistency}} and Condition {\rm\ref{as:partial2}} hold. Then
\[
\sup_{\mathcal{F}\in\mathscr{F}}\big|[\nabla_{\btheta_{\mathcal{S}}}\bar{\bg}_{\mathcal{F}}(\widehat{\btheta}_n)-\mathbb{E}\{\nabla_{\btheta_{\mathcal{S}}}\bg_{i,\mathcal{F}}(\btheta_0)\}]\bz\big|_2=|\bz|_2\big[O_p\{s^{3/2}\ell_n^{1/2}\varpi_n^{1/2}b_n^{1/(2\beta)}\}+O_p\{(n^{-1}s\ell_n\omega_n\log r)^{1/2}\}\big]
\]
holds uniformly for any $\bz\in\mathbb{R}^s$, where $\mathscr{F}$ is defined in Lemma {\rm\ref{la:covc}}.
\end{lemma}

\begin{lemma}\label{la:3}
Let $\widehat{\bJ}_{\mathcal{F}}=\{\nabla_{\btheta_{\mathcal{S}}}\bar{\bg}_{\mathcal{F}}(\widehat{\btheta}_n)\}^\T\widehat{\bV}_{\mathcal{F}}^{-1}(\widehat{\btheta}_n)\{\nabla_{\btheta_{\mathcal{S}}}\bar{\bg}_{\mathcal{F}}(\widehat{\btheta}_n)\}$ for any $\mathcal{F}\in\mathscr{F}$, where $\mathscr{F}$ is defined in Lemma {\rm\ref{la:covc}}. Assume the conditions for Lemma {\rm\ref{la:2}} and Condition {\rm\ref{as:eig2}} hold. If $s^2\ell_n^2b_n^{1/\beta}\varrho_n^{1/2}\max\{\omega_n,s\varpi_n\}\log r=o(1)$, $n^{-1}\ell_n^2s\omega_n \varrho_n^{1/2}(\log r)^2=o(1)$ and $n^{-1}\ell_n^{3}\varrho_n^{3/2}(\log r)^2=o(1)$, we have
\[
\sup_{\mathcal{F}\in\mathscr{F}}\Big|\mathbb{P}\big[n^{1/2}\balpha^\T\widehat{\bJ}_{\mathcal{F}}^{-1/2}\{\nabla_{\btheta_{\mathcal{S}}}\bar{\bg}_{\mathcal{F}}(\widehat{\btheta}_n)\}^\T\widehat{\bV}_{\mathcal{F}}^{-1}(\widehat{\btheta}_n)\bar{\bg}_{\mathcal{F}}(\btheta_0)\leq u\big]-\Phi(u)\Big|\rightarrow0,~~~\textrm{as}~n\rightarrow\infty,
\]
for any $u\in\mathbb{R}$ and $\balpha\in\mathbb{R}^s$, where $\Phi(\cdot)$ is the cumulative distribution function of the standard normal distribution.
\end{lemma}

Now we begin the proof of Theorem \ref{tm:2}. Recall $\widehat{\bJ}_{\mathcal{R}_n}=\{\nabla_{\btheta_{\mathcal{S}}}\bar{\bg}_{\mathcal{R}_n}(\widehat{\btheta}_n)\}^\T \widehat{\bV}_{\mathcal{R}_n}^{-1}(\widehat{\btheta}_n)\{\nabla_{\btheta_{\mathcal{S}}}\bar{\bg}_{\mathcal{R}_n}(\widehat{\btheta}_n)\}$. For any $\balpha\in\mathbb{R}^s$ with unit $L_2$-norm, let $\bdelta=\widehat{\bJ}_{\mathcal{R}_n}^{-1/2}\balpha,$ then
\[
\begin{split}
|\{\nabla_{\btheta_{\mathcal{S}}}\bar{\bg}_{\mathcal{R}_n}(\widehat{\btheta}_n)\}\bdelta|_2^2=&~\balpha^\T(\bU^\T\bU)^{-1/2}\bU^\T \widehat{\bV}_{\mathcal{R}_n}(\widehat{\btheta}_n)\bU(\bU^\T\bU)^{-1/2}\balpha\\
\leq&~\lambda_{\max}\{\widehat{\bV}_{\mathcal{R}_n}(\widehat{\btheta}_n)\}\cdot|\bU(\bU^\T\bU)^{-1/2}\balpha|_2^2\\
=&~\lambda_{\max}\{\widehat{\bV}_{\mathcal{R}_n}(\widehat{\btheta}_n)\},
\end{split}
\]
where $\bU=\widehat{\bV}_{\mathcal{R}_n}^{-1/2}(\widehat{\btheta}_n)\{\nabla_{\btheta_{\mathcal{S}}}\bar{\bg}_{\mathcal{R}_n}(\widehat{\btheta}_n)\}$. Thus, by Lemma \ref{la:covc}, $|\{\nabla_{\btheta_{\mathcal{S}}}\bar{\bg}_{\mathcal{R}_n}(\widehat{\btheta}_n)\}\bdelta|_2=O_p(1)$. Meanwhile, notice that
$|\bdelta|_2=O_p(1)$. 
Lemma \ref{la:1} yields that \[\bigg|\bigg\{\frac{1}{n}\sum_{i=1}^n\frac{\nabla_{\btheta_{\mathcal{S}}}\bg_{i,\mathcal{R}_n}(\widehat{\btheta}_n)}{1+\widehat{\blambda}^\T_{\mathcal{R}_n}\bg_{i,\mathcal{R}_n}(\widehat{\btheta}_n)}\bigg\}\bdelta\bigg|_2=O_p(1).\]
As shown in Section \ref{se:1}, $|\bar{\bg}_{\mathcal{M}_{\widehat{\btheta}_n}}(\widehat{\btheta}_n)-\nu\rho_2'(0^+)\sgn\{\bar{\bg}_{\mathcal{M}_{\widehat{\btheta}_n}}(\widehat{\btheta}_n)\}|_2=O_p(\ell_n^{1/2}n^{-1/2})+O_p\{s^{1/2}\chi_n^{1/2}b_n^{1/(4\beta)}\}$. From Proposition \ref{pn:2}, we have $|\bar{\bg}_{\mathcal{R}_n}(\widehat{\btheta}_n)-\widehat{\bfeta}_{\mathcal{R}_n}|_2=O_p(\ell_n^{1/2}n^{-1/2})+O_p\{s^{1/2}\chi_n^{1/2}b_n^{1/(4\beta)}\}$.
Following Lemmas \ref{la:1} and \ref{la:2}, (\ref{eq:expan3}) leads to
\[
\begin{split}
&~\bdelta^\T\{\nabla_{\btheta_{\mathcal{S}}}\bar{\bg}_{\mathcal{R}_n}(\widehat{\btheta}_n)\}^\T \widehat{\bV}_{\mathcal{R}_n}^{-1}(\widehat{\btheta}_n)\{\bar{\bg}_{\mathcal{R}_n}(\widehat{\btheta}_n)-\widehat{\bfeta}_{\mathcal{R}_n}\}\\
=&~O_p\big(\ell_n^{1/2}\max\{\ell_nn^{-1},s\chi_nb_n^{1/(2\beta)}\}\max\{s^{1/2}\omega_n^{1/2},n^{1/\gamma}\}\big)+O_p(s^{1/2}\chi_n).
\end{split}
\]
Expanding $\bar{\bg}_{\mathcal{R}_n}(\widehat{\btheta}_n)$ around $\btheta=\btheta_0$, it holds w.p.a.1 that
\begin{equation}\label{eq:expan4}
\begin{split}
&~\bdelta^\T\{\nabla_{\btheta_{\mathcal{S}}}\bar{\bg}_{\mathcal{R}_n}(\widehat{\btheta}_n)\}^\T \widehat{\bV}_{\mathcal{R}_n}^{-1}(\widehat{\btheta}_n)[\{\nabla_{\btheta_{\mathcal{S}}}\bar{\bg}_{\mathcal{R}_n}(\widetilde{\btheta})\}(\widehat{\btheta}_{n,\mathcal{S}}-\btheta_{0,\mathcal{S}})-\widehat{\bfeta}_{\mathcal{R}_n}]\\
=&-\bdelta^\T\{\nabla_{\btheta_{\mathcal{S}}}\bar{\bg}_{\mathcal{R}_n}(\widehat{\btheta}_n)\}^\T \widehat{\bV}_{\mathcal{R}_n}^{-1}(\widehat{\btheta}_n)\bar{\bg}_{\mathcal{R}_n}(\btheta_0)+O_p(s^{1/2}\chi_n)\\
&+O_p\big(\ell_n^{1/2}\max\{\ell_nn^{-1},s\chi_nb_n^{1/(2\beta)}\}\max\{s^{1/2}\omega_n^{1/2},n^{1/\gamma}\}\big),
\end{split}
\end{equation}
where $\widetilde{\btheta}$ is on the line joining $\btheta_0$ and $\widehat{\btheta}_n$. Notice that $|\bar{\bg}_{\mathcal{R}_n}(\widehat{\btheta}_n)-\bar{\bg}_{\mathcal{R}_n}(\btheta_0)|_2\leq |\bar{\bg}_{\mathcal{R}_n}(\widehat{\btheta}_n)|_2+|\bar{\bg}_{\mathcal{R}_n}(\btheta_0)|_2=O_p(\ell_n^{1/2}\nu)+O_p\{s^{1/2}\chi_n^{1/2}b_n^{1/(4\beta)}\}$. By Taylor expansion, $|\bar{\bg}_{\mathcal{R}_n}(\widehat{\btheta}_n)-\bar{\bg}_{\mathcal{R}_n}(\btheta_0)|_2\geq \lambda_{\min}([\nabla_{\btheta_{\mathcal{S}}}\bar{\bg}_{\mathcal{R}_n}(\dot{\btheta})]^\T[\nabla_{\btheta_{\mathcal{S}}}\bar{\bg}_{\mathcal{R}_n}(\dot{\btheta})])|\widehat{\btheta}_{n,\mathcal{S}}-\btheta_{0,\mathcal{S}}|_2$ for some $\dot{\btheta}$ lying on the line jointing $\btheta_0$ and $\widehat{\btheta}_n$. Same as Lemma \ref{la:2}, $\lambda_{\min}([\nabla_{\btheta_{\mathcal{S}}}\bar{\bg}_{\mathcal{R}_n}(\dot{\btheta})]^\T[\nabla_{\btheta_{\mathcal{S}}}\bar{\bg}_{\mathcal{R}_n}(\dot{\btheta})])$ is bounded away from zero w.p.a.1, which implies $|\widehat{\btheta}_{n,\mathcal{S}}-\btheta_{0,\mathcal{S}}|_2=O_p(\ell_n^{1/2}\nu)+O_p\{s^{1/2}\chi_n^{1/2}b_n^{1/(4\beta)}\}$. Together with Condition \ref{as:partial2}, it holds that
$
|\{\nabla_{\btheta_{\mathcal{S}}}\bar{\bg}_{\mathcal{R}_n}(\widetilde{\btheta})-\nabla_{\btheta_{\mathcal{S}}}\bar{\bg}_{\mathcal{R}_n}(\widehat{\btheta}_n)\}(\widehat{\btheta}_{n,\mathcal{S}}-\btheta_{0,\mathcal{S}})|_2=O_p(\ell_n^{3/2}s\varpi_n^{1/2}\nu^2)+O_p\{\ell_n^{1/2}s^2\varpi_n^{1/2}\chi_nb_n^{1/(2\beta)}\}.
$
Therefore, (\ref{eq:expan4}) leads to
\[
\begin{split}
&~\bdelta^\T\widehat{\bJ}_{\mathcal{R}_n}\big[\widehat{\btheta}_{n,\mathcal{S}}-\btheta_{0,\mathcal{S}}-\widehat{\bJ}_{\mathcal{R}_n}^{-1}\{\nabla_{\btheta_{\mathcal{S}}}\bar{\bg}_{\mathcal{R}_n}(\widehat{\btheta}_n)\}^\T \widehat{\bV}_{\mathcal{R}_n}^{-1}(\widehat{\btheta}_n)\widehat{\bfeta}_{\mathcal{R}_n}\big]\\
=&-\bdelta^\T\{\nabla_{\btheta_{\mathcal{S}}}\bar{\bg}_{\mathcal{R}_n}(\widehat{\btheta}_n)\}^\T \widehat{\bV}_{\mathcal{R}_n}^{-1}(\widehat{\btheta}_n)\bar{\bg}_{\mathcal{R}_n}(\btheta_0)+O_p(\ell_n^{3/2}s\varpi_n^{1/2}\nu^2)+O_p\{\ell_n^{1/2}s^2\varpi_n^{1/2}\chi_nb_n^{1/(2\beta)}\}\\
&+O_p\big(\ell_n^{1/2}\max\{\ell_nn^{-1},s\chi_nb_n^{1/(2\beta)}\}\max\{s^{1/2}\omega_n^{1/2},n^{1/\gamma}\}\big)+O_p(s^{1/2}\chi_n)\\
=&-\balpha^\T\widehat{\bJ}_{\mathcal{R}_n}^{-1/2}\{\nabla_{\btheta_{\mathcal{S}}}\bar{\bg}_{\mathcal{R}_n}(\widehat{\btheta}_n)\}^\T \widehat{\bV}_{\mathcal{R}_n}^{-1}(\widehat{\btheta}_n)\bar{\bg}_{\mathcal{R}_n}(\btheta_0)+O_p(\ell_n^{3/2}s\varpi_n^{1/2}\nu^2)+O_p\{\ell_n^{1/2}s^2\varpi_n^{1/2}\chi_nb_n^{1/(2\beta)}\}\\
&+O_p\big(\ell_n^{1/2}\max\{\ell_nn^{-1},s\chi_nb_n^{1/(2\beta)}\}\max\{s^{1/2}\omega_n^{1/2},n^{1/\gamma}\}\big)+O_p(s^{1/2}\chi_n).
\end{split}
\]
Lemma \ref{la:3} leads to $n^{1/2}\balpha^\T\widehat{\bJ}_{\mathcal{R}_n}^{-1/2}\{\nabla_{\btheta_{\mathcal{S}}}\bar{\bg}_{\mathcal{R}_n}(\widehat{\btheta}_n)\}^\T \widehat{\bV}_{\mathcal{R}_n}^{-1}(\widehat{\btheta}_n)\bar{\bg}_{\mathcal{R}_n}(\btheta_0)\rightarrow_d N(0,1)$ as $n\rightarrow\infty$. We complete the proof of Theorem \ref{tm:2}. $\hfill\Box$

\newpage

\setcounter{page}{1}
\pagestyle{fancy}
\fancyhf{}
\rhead{\bfseries\thepage}
\lhead{\bfseries SUPPLEMENTARY MATERIAL}

\setcounter{page}{1}
\begin{center}
{\bf\Large Supplementary Material for ``A New Scope of Penalized Empirical Likelihood with High-dimensional Estimating Equations'' by Chang, Tang and Wu.}
\end{center}

\section*{Proof of Lemma \ref{la:covc}} Notice that $\|\widehat{\bV}_{\mathcal{F}}(\btheta)-\bV_{\mathcal{F}}(\btheta_0)\|_2\leq \|\widehat{\bV}_{\mathcal{F}}(\btheta)-\widehat{\bV}_{\mathcal{F}}(\btheta_0)\|_2+\|\widehat{\bV}_{\mathcal{F}}({\btheta}_0)-\bV_{\mathcal{F}}(\btheta_0)\|_2$ for any $\mathcal{F}\in\mathscr{F}$ and $\btheta\in\bTheta_n$. Following the moderate deviation of self-normalized sums \citep{JingShaoWang2003} and Condition \ref{as:pa1}, it holds that $\max_{1\leq j_1,j_2\leq r}|n^{-1}\sum_{i=1}^ng_{i,j_1}(\btheta_0)g_{i,j_2}(\btheta_0)-\mathbb{E}\{g_{i,j_1}(\btheta_0)g_{i,j_2}(\btheta_0)\}|=O_p\{(n^{-1}\varrho_n\log r)^{1/2}\}$, which implies $\sup_{\mathcal{F}\in\mathscr{F}}\|\widehat{\bV}_{\mathcal{F}}(\btheta_0)-\bV_{\mathcal{F}}(\btheta_0)\|_2=O_p\{\ell_n(n^{-1}\varrho_n\log r)^{1/2}\}$ provided that $\log r=o(n^{1/3})$. For any $\bz\in\mathbb{R}^{|\mathcal{F}|}$ with unit $L_2$-norm, we have
\[
\begin{split}
\big|\bz^\T\{\widehat{\bV}_{\mathcal{F}}(\btheta)-\widehat{\bV}_{\mathcal{F}}(\btheta_0)\}\bz\big|\leq&~\frac{1}{n}\sum_{i=1}^n|\bg_{i,\mathcal{F}}(\btheta)-\bg_{i,\mathcal{F}}(\btheta_0)|_2^2\\
&+2\lambda_{\max}^{1/2}\{\widehat{\bV}_{\mathcal{F}}(\btheta_0)\}\bigg\{\frac{1}{n}\sum_{i=1}^n|\bg_{i,\mathcal{F}}(\btheta)-\bg_{i,\mathcal{F}}(\btheta_0)|_2^2\bigg\}^{1/2},
\end{split}
\]
which implies
\[
\begin{split}
\sup_{\mathcal{F}\in\mathscr{F}}\|\widehat{\bV}_{\mathcal{F}}(\btheta)-\widehat{\bV}_{\mathcal{F}}(\btheta_0)\|_2\leq&~\sup_{\mathcal{F}\in\mathscr{F}}\bigg\{\frac{1}{n}\sum_{i=1}^n|\bg_{i,\mathcal{F}}(\btheta)-\bg_{i,\mathcal{F}}(\btheta_0)|_2^2\bigg\}\\
&+2\sup_{\mathcal{F}\in\mathscr{F}}\lambda_{\max}^{1/2}\{\widehat{\bV}_{\mathcal{F}}(\btheta_0)\}\cdot\sup_{\mathcal{F}\in\mathscr{F}}\bigg\{\frac{1}{n}\sum_{i=1}^n|\bg_{i,\mathcal{F}}(\btheta)-\bg_{i,\mathcal{F}}(\btheta_0)|_2^2\bigg\}^{1/2}.
\end{split}
\]
Write $\btheta=(\btheta_{\mathcal{S}}^\T,\btheta_{\mathcal{S}^c}^\T)^\T$ with $\btheta_{\mathcal{S}}\in\mathbb{R}^s$. By Taylor expansion and Cauchy-Schwarz inequality, we have
\[
\begin{split}
\frac{1}{n}\sum_{i=1}^n|\bg_{i,\mathcal{F}}(\btheta)-\bg_{i,\mathcal{F}}(\btheta_0)|_2^2\leq&~\frac{2}{n}\sum_{i=1}^n\bigg|\frac{\partial\bg_{i,\mathcal{F}}(\widetilde{\btheta})}{\partial \btheta_{\mathcal{S}}}(\btheta_{\mathcal{S}}-\btheta_{0,\mathcal{S}})\bigg|_2^2+\frac{2}{n}\sum_{i=1}^n\bigg|\frac{\partial\bg_{i,\mathcal{F}}(\widetilde{\btheta})}{\partial \btheta_{\mathcal{S}^c}}{\btheta}_{\mathcal{S}^c}\bigg|_2^2\\
\leq&~2|{\btheta}_{\mathcal{S}}-\btheta_{0,\mathcal{S}}|_1^2\max_{1\leq k_1,k_2\leq s}\bigg|\frac{1}{n}\sum_{i=1}^n\bigg\{\frac{\partial \bg_{i,\mathcal{F}}(\widetilde{\btheta})}{\partial\theta_{k_1}}\bigg\}^\T\bigg\{\frac{\partial\bg_{i,\mathcal{F}}(\widetilde{\btheta})}{\partial\theta_{k_2}}\bigg\}\bigg|\\
&+2|\btheta_{\mathcal{S}^c}|_1^2\max_{s+1\leq k_1,k_2\leq p}\bigg|\frac{1}{n}\sum_{i=1}^n\bigg\{\frac{\partial \bg_{i,\mathcal{F}}(\widetilde{\btheta})}{\partial\theta_{k_1}}\bigg\}^\T\bigg\{\frac{\partial\bg_{i,\mathcal{F}}(\widetilde{\btheta})}{\partial\theta_{k_2}}\bigg\}\bigg|,
\end{split}
\]
for some $\widetilde{\btheta}$ lying on the jointing line between $\btheta_0$ and $\btheta$. By Condition \ref{as:pa1},
\[
\begin{split}
\max_{1\leq k_1,k_2\leq s}\bigg|\frac{1}{n}\sum_{i=1}^n\bigg\{\frac{\partial \bg_{i,\mathcal{F}}(\widetilde{\btheta})}{\partial\theta_{k_1}}\bigg\}^\T\bigg\{\frac{\partial\bg_{i,\mathcal{F}}(\widetilde{\btheta})}{\partial\theta_{k_2}}\bigg\}\bigg|\leq&~ \sum_{j\in\mathcal{F}}\max_{k\in\mathcal{S}}\bigg\{\frac{1}{n}\sum_{i=1}^n\bigg|\frac{\partial g_{i,j}(\widetilde{\btheta})}{\partial \theta_k}\bigg|^2\bigg\}\\
\leq&~|\mathcal{F}|\max_{1\leq j\leq r}\max_{k\in\mathcal{S}}\bigg\{\frac{1}{n}\sum_{i=1}^n\bigg|\frac{\partial g_{i,j}(\widetilde{\btheta})}{\partial \theta_k}\bigg|^2\bigg\}\\
=&~O_p(\ell_n\omega_n).
\end{split}
\]
Similarly, we have
\[
\max_{s+1\leq k_1,k_2\leq p}\bigg|\frac{1}{n}\sum_{i=1}^n\bigg\{\frac{\partial \bg_{i,\mathcal{F}}(\widetilde{\btheta})}{\partial\theta_{k_1}}\bigg\}^\T\bigg\{\frac{\partial\bg_{i,\mathcal{F}}(\widetilde{\btheta})}{\partial\theta_{k_2}}\bigg\}\bigg|=O_p(\ell_n\xi_n).
\]
 Therefore,
\[
\frac{1}{n}\sum_{i=1}^n|\bg_{i,\mathcal{F}}(\btheta)-\bg_{i,\mathcal{F}}(\btheta_0)|_2^2=O_p(s^2\ell_n\omega_nb_n^{1/\beta})
\]
holds uniformly for $\btheta\in\bTheta_n$. Meanwhile, by Condition \ref{as:smalleigenvalue}, it holds that $\sup_{\mathcal{F}\in\mathscr{F}}\lambda_{\max}\{\widehat{\bV}_{\mathcal{F}}(\btheta_0)\}\leq C$ w.p.a.1. Then $ \sup_{\btheta\in\bTheta_n}\sup_{\mathcal{F}\in\mathscr{F}}\|\widehat{\bV}_{\mathcal{F}}(\btheta)-\widehat{\bV}_{\mathcal{F}}(\btheta_0)\|_2=O_p\{s(\ell_n\omega_nb_n^{1/\beta})^{1/2}\}$. Thus we complete the proof of Lemma \ref{la:covc}. $\hfill\Box$

\section*{Proof of Lemma \ref{la:1}}

As shown in Section \ref{se:1}, $|\widehat{\blambda}|_2=O_p(\ell_n^{1/2}n^{-1/2})+O_p\{s^{1/2}\chi_n^{1/2}b_n^{1/(4\beta)}\}$ and  $\max_{1\leq i\leq n}|\widehat{\blambda}_{\mathcal{R}_n}^\T\bg_{i,\mathcal{R}_n}(\widehat{\btheta}_n)|=O_p(\ell_nn^{-1/2+1/\gamma})+O_p\{\ell_n^{1/2}s^{1/2}\chi_n^{1/2}b_n^{1/(4\beta)}n^{1/\gamma}\}=o_p(1)$. Notice that $|(1+x)^{-2}-1|\leq 5|x|$ for any $|x|<1/2$, by Lemma \ref{la:covc}, it holds that w.p.a.1
\[
\begin{split}
\bigg\|\frac{1}{n}\sum_{i=1}^n\frac{\bg_{i,\mathcal{R}_n}(\widehat{\btheta}_n)\bg_{i,\mathcal{R}_n}(\widehat{\btheta}_n)^\T}{\{1+c\widehat{\blambda}_{\mathcal{R}_n}^\T\bg_{i,\mathcal{R}_n}(\widehat{\btheta}_n)\}^2}-\widehat{\bV}_{\mathcal{R}_n}(\widehat{\btheta}_n)\bigg\|_2\leq&~5\lambda_{\max}\{\widehat{\bV}_{\mathcal{R}_n}(\widehat{\btheta}_n)\}\max_{1\leq i\leq n}|\widehat{\blambda}_{\mathcal{R}_n}^\T\bg_{i,\mathcal{R}_n}(\widehat{\btheta}_n)|\\
=&~O_p(\ell_nn^{-1/2+1/\gamma})+O_p\{\ell_n^{1/2}s^{1/2}\chi_n^{1/2}b_n^{1/(4\beta)}n^{1/\gamma}\}.
\end{split}
\]
For the second result, by Taylor expansion and Cauchy-Schwarz inequality, it holds that w.p.a.1
\begin{equation}\label{eq:d2}
\begin{split}
&~\bigg|\bigg\{\frac{1}{n}\sum_{i=1}^n\frac{\nabla_{\btheta_{\mathcal{S}}}\bg_{i,\mathcal{R}_n}(\widehat{\btheta}_n)}{1+\widehat{\blambda}_{\mathcal{R}_n}^\T\bg_{i,\mathcal{R}_n}(\widehat{\btheta}_n)}-\nabla_{\btheta_{\mathcal{S}}}\bar{\bg}_{\mathcal{R}_n}(\widehat{\btheta}_n)\bigg\}\bz\bigg|_2^2\\
\leq&~\bigg[\frac{1}{n}\sum_{i=1}^n\frac{\widehat{\blambda}_{\mathcal{R}_n}^\T\bg_{i,\mathcal{R}_n}(\widehat{\btheta}_n)\bg_{i,\mathcal{R}_n}(\widehat{\btheta}_n)^\T\widehat{\blambda}_{\mathcal{R}_n}}{\{1+c\widehat{\blambda}_{\mathcal{R}_n}^\T\bg_{i,\mathcal{R}_n}(\widehat{\btheta}_n)\}^4}\bigg]\bigg[\frac{1}{n}\sum_{i=1}^n\bz^\T\{\nabla_{\btheta_{\mathcal{S}}}\bg_{i,\mathcal{R}_n}(\widehat{\btheta}_n)\}^\T\{\nabla_{\btheta_{\mathcal{S}}}\bg_{i,\mathcal{R}_n}(\widehat{\btheta}_n)\}\bz\bigg]\\
\leq&~\widehat{\blambda}_{\mathcal{R}_n}^\T\widehat{\bV}_{\mathcal{R}_n}(\widehat{\btheta}_n)\widehat{\blambda}_{\mathcal{R}_n}\bigg[\frac{1}{n}\sum_{i=1}^n\bz^\T\{\nabla_{\btheta_{\mathcal{S}}}\bg_{i,\mathcal{R}_n}(\widehat{\btheta}_n)\}^\T\{\nabla_{\btheta_{\mathcal{S}}}\bg_{i,\mathcal{R}_n}(\widehat{\btheta}_n)\}\bz\bigg]\{1+o_p(1)\}
\end{split}
\end{equation}
for some $|c|<1$. By Lemma \ref{la:covc}, it holds that $\widehat{\blambda}^\T_{\mathcal{R}_n}\widehat{\bV}_{\mathcal{R}_n}(\widehat{\btheta}_n)\widehat{\blambda}_{\mathcal{R}_n}\leq \lambda_{\max}\{\widehat{\bV}_{\mathcal{R}_n}(\widehat{\btheta}_n)\}|\widehat{\blambda}_{\mathcal{R}_n}|_2^2=O_p(\ell_nn^{-1})+O_p\{s\chi_nb_n^{1/(2\beta)}\}$. Meanwhile, write $\bz=(z_1,\ldots,z_{s})^\T$, by Cauchy-Schwarz inequality and Condition \ref{as:pa1},
\[
\begin{split}
&\frac{1}{n}\sum_{i=1}^n\bz^\T\{\nabla_{\btheta_{\mathcal{S}}}\bg_{i,\mathcal{R}_n}(\widehat{\btheta}_n)\}^\T\{\nabla_{\btheta_{\mathcal{S}}}\bg_{i,\mathcal{R}_n}(\widehat{\btheta}_n)\}\bz\leq\frac{|\bz|^2_2}{n}\sum_{i=1}^n\sum_{j\in\mathcal{R}_n}\sum_{k=1}^{s}\bigg|\frac{\partial g_{i,j}(\widehat{\btheta}_n)}{\partial \theta_k}\bigg|^2=|\bz|_2^2\cdot O_p(\ell_ns\omega_n).
\end{split}
\]
Therefore, (\ref{eq:d2}) leads to
\begin{equation}\label{eq:term4}
\begin{split}
&\bigg|\bigg\{\frac{1}{n}\sum_{i=1}^n\frac{\nabla_{\btheta_{\mathcal{S}}}\bg_{i,\mathcal{R}_n}(\widehat{\btheta}_n)}{1+\widehat{\blambda}^\T_{\mathcal{R}_n}\bg_{i,\mathcal{R}_n}(\widehat{\btheta}_n)}-\nabla_{\btheta_{\mathcal{S}}}\bar{\bg}_{\mathcal{R}_n}(\widehat{\btheta}_n)\bigg\}\bz\bigg|_2\\
=&~|\bz|_2 [O_p(\ell_ns^{1/2}\omega_n^{1/2}n^{-1/2})+O_p\{\ell_n^{1/2}s\omega_n^{1/2}\chi_n^{1/2}b_n^{1/(4\beta)}\}].
\end{split}
\end{equation}
We complete the proof of Lemma \ref{la:1}. $\hfill\Box$

\section*{Proof of Lemma \ref{la:2}}

Notice that
\begin{equation}\label{eq:term2}
\begin{split}
&~\big|\big[\nabla_{\btheta_{\mathcal{S}}}\bar{\bg}_{\mathcal{F}}(\widehat{\btheta}_n)-\mathbb{E}\{\nabla_{\btheta_{\mathcal{S}}} \bg_{i,\mathcal{F}}(\btheta_0)\}\big]\bz\big|_2\\
\leq&~\big|\{\nabla_{\btheta_{\mathcal{S}}}\bar{\bg}_{\mathcal{F}}(\widehat{\btheta}_n)-\nabla_{\btheta_{\mathcal{S}}}\bar{\bg}_{\mathcal{F}}(\btheta_0)\}\bz\big|_2+\big|[\nabla_{\btheta_{\mathcal{S}}}\bar{\bg}_{\mathcal{F}}(\btheta_0)-\mathbb{E}\{\nabla_{\btheta_{\mathcal{S}}} \bg_{i,\mathcal{F}}(\btheta_0)\}]\bz\big|_2
\end{split}
\end{equation}
for any $\bz\in\mathbb{R}^s$.
By Taylor expansion, Jensen's inequality and Cauchy-Schwarz inequality, it holds that w.p.a.1
\[
\begin{split}
\big|\{\nabla_{\btheta_{\mathcal{S}}}\bar{\bg}_{\mathcal{F}}(\widehat{\btheta}_n)-\nabla_{\btheta_{\mathcal{S}}}\bar{\bg}_{\mathcal{F}}(\btheta_0)\}\bz\big|_2^2=&\sum_{j\in\mathcal{F}}\bigg\{\frac{1}{n}\sum_{i=1}^n\sum_{k=1}^{s}z_k\sum_{l=1}^{s}\frac{\partial^2 g_{i,j}(\widetilde{\btheta})}{\partial\theta_k\partial\theta_l}(\widehat{\theta}_l-\theta_l^0)\bigg\}^2\\
\leq&~\frac{|\bz|_2^2}{n}\sum_{j\in\mathcal{F}}\sum_{i=1}^n\sum_{k=1}^{s}\sum_{l=1}^{s}\bigg|\frac{\partial^2 g_{i,j}(\widetilde{\btheta})}{\partial\theta_k\partial\theta_l}\bigg|^2|\widehat{\btheta}_{n,\mathcal{S}}-\btheta_{0,\mathcal{S}}|_2^2,
\end{split}
\]
where $\widetilde{\btheta}$ lies on the jointing line between $\btheta_0$ and $\widehat{\btheta}_n$. It follows from Condition \ref{as:partial2} that
\begin{equation}\label{eq:term3}
\begin{split}
&\sup_{\mathcal{F}\in\mathscr{F}}\big|\{\nabla_{\btheta_{\mathcal{S}}}\bar{\bg}_{\mathcal{F}}(\widehat{\btheta}_n)-\nabla_{\btheta_{\mathcal{S}}}\bar{\bg}_{\mathcal{F}}(\btheta_0)\}\bz\big|_2=|\bz|_2\cdot O_p\{s^{3/2}\ell_n^{1/2}\varpi_n^{1/2}b_n^{1/(2\beta)}\}.
\end{split}
\end{equation}
On the other hand, by Cauchy-Schwarz inequality, it holds that
\[
\begin{split}
&\big|[\nabla_{\btheta_{\mathcal{S}}}\bar{\bg}_{\mathcal{F}}(\btheta_0)-\mathbb{E}\{\nabla_{\btheta_{\mathcal{S}}} \bg_{i,\mathcal{F}}(\btheta_0)\}]\bz\big|_2^2\\
&~~~~~~~~~~~~~=\sum_{j\in\mathcal{F}}\bigg(\frac{1}{n}\sum_{i=1}^n\sum_{k=1}^{s}z_k\bigg[\frac{\partial g_{i,j}(\btheta_0)}{\partial \theta_k}-\mathbb{E}\bigg\{\frac{\partial g_{i,j}(\btheta_0)}{\partial \theta_k}\bigg\}\bigg]\bigg)^2\\
&~~~~~~~~~~~~~\leq |\bz|_2^2\sum_{j\in\mathcal{F}}\sum_{k=1}^{s}\bigg(\frac{1}{n}\sum_{i=1}^n\bigg[\frac{\partial g_{i,j}(\btheta_0)}{\partial \theta_k}-\mathbb{E}\bigg\{\frac{\partial g_{i,j}(\btheta_0)}{\partial \theta_k}\bigg\}\bigg]\bigg)^2.
\end{split}
\]
Notice that
\[
\sup_{1\leq j\leq r}\sup_{1\leq k\leq s}\bigg|\frac{1}{n}\sum_{i=1}^n\bigg[\frac{\partial g_{i,j}(\btheta_0)}{\partial \theta_k}-\mathbb{E}\bigg\{\frac{\partial g_{i,j}(\btheta_0)}{\partial \theta_k}\bigg\}\bigg]\bigg|= O_p\{(n^{-1}\omega_n\log r)^{1/2}\},
\]
therefore
\[
\sup_{\mathcal{F}\in\mathscr{F}}\big|[\nabla_{\btheta_{\mathcal{S}}}\bar{\bg}_{\mathcal{F}}(\btheta_0)-\mathbb{E}\{\nabla_{\btheta_{\mathcal{S}}} \bg_{i,\mathcal{F}}(\btheta_0)\}]\bz\big|_2=|\bz|_2\cdot O_p\{(n^{-1}s\ell_n\omega_n\log r)^{1/2}\}.
\]
Together with (\ref{eq:term3}), (\ref{eq:term2}) yields that
\begin{equation*}\label{eq:te2}
\begin{split}
&~\sup_{\mathcal{F}\in\mathscr{F}}\big|[\nabla_{\btheta_{\mathcal{S}}}\bar{\bg}_{\mathcal{F}}(\widehat{\btheta}_n)-\mathbb{E}\{\nabla_{\btheta_{\mathcal{S}}} \bg_{i,\mathcal{F}}(\btheta_0)\}]\bz\big|_2\\
=&~|\bz|_2\big[O_p\{s^{3/2}\ell_n^{1/2}\varpi_n^{1/2}b_n^{1/(2\beta)}\}+O_p\{(n^{-1}s\ell_n\omega_n\log r)^{1/2}\}\big].
\end{split}
\end{equation*}
We complete the proof of Lemma \ref{la:2}. $\hfill\Box$

\section*{Proof of Lemma \ref{la:3}}

For any $\mathcal{F}\in\mathscr{F}$, let $\bJ_{\mathcal{F}}=[\mathbb{E}\{\nabla_{\btheta_{\mathcal{S}}}{\bg}_{i,\mathcal{F}}(\btheta_0)\}]^\T{\bV}_{\mathcal{F}}^{-1}(\btheta_0)[\mathbb{E}\{\nabla_{\btheta_{\mathcal{S}}}{\bg}_{i,\mathcal{F}}(\btheta_0)\}]$. Given $\mathcal{F}$, by Lindeberg-Feller Central Limit Theorem, we have \[
n^{1/2}\balpha^\T\bJ_{\mathcal{F}}^{-1/2}[\mathbb{E}\{\nabla_{\btheta_{\mathcal{S}}}{\bg}_{i,\mathcal{F}}(\btheta_0)\}]^\T{\bV}_{\mathcal{F}}^{-1}(\btheta_0)\bar{\bg}_{\mathcal{F}}(\btheta_0)\xrightarrow{d}N(0,1).
\]
Let $Z_{i,\mathcal{F}}=\balpha^\T\bJ_{\mathcal{F}}^{-1/2}[\mathbb{E}\{\nabla_{\btheta_{\mathcal{S}}}{\bg}_{i,\mathcal{F}}(\btheta_0)\}]^\T{\bV}_{\mathcal{F}}^{-1}(\btheta_0){\bg}_{i,\mathcal{F}}(\btheta_0)$. Applying Berry-Esseen inequality, we have
\[
\sup_{u\in\mathbb{R}}\Big|\mathbb{P}\big[n^{1/2}\balpha^\T\bJ_{\mathcal{F}}^{-1/2}[\mathbb{E}\{\nabla_{\btheta_{\mathcal{S}}}{\bg}_{i,\mathcal{F}}(\btheta_0)\}]^\T{\bV}_{\mathcal{F}}^{-1}(\btheta_0)\bar{\bg}_{\mathcal{F}}(\btheta_0)\leq u\big]-\Phi(u)\Big|\leq Cn^{-1/2}\mathbb{E}(|Z_{i,\mathcal{F}}|^3),
\]
where $C$ is a uniform positive constant independent of $\mathcal{F}$. By Cauchy-Schwarz inequality,
\[
\begin{split}
|Z_{i,\mathcal{F}}|^2\leq&~ |{\bV}_{\mathcal{F}}^{-1/2}(\btheta_0)[\mathbb{E}\{\nabla_{\btheta_{\mathcal{S}}}{\bg}_{i,\mathcal{F}}(\btheta_0)\}]\bJ_{\mathcal{F}}^{-1/2}\balpha|_2^2|{\bV}_{\mathcal{F}}^{-1/2}(\btheta_0){\bg}_{i,\mathcal{F}}(\btheta_0)|_2^2\\
\leq&~\lambda_{\min}^{-1}\{\bV_{\mathcal{F}}(\btheta_0)\}|\bg_{i,\mathcal{F}}(\btheta_0)|_2^2,
\end{split}
\]
which implies
\[
\mathbb{E}(|Z_{i,\mathcal{F}}|^3)\leq \lambda_{\min}^{-3/2}\{\bV_{\mathcal{F}}(\btheta_0)\}\mathbb{E}\{|\bg_{i,\mathcal{F}}(\btheta_0)|_2^3\}\leq C\lambda_{\min}^{-3/2}\{\bV_{\mathcal{F}}(\btheta_0)\}\ell_n^{3/2}
\]
for a uniform positive constant $C$ independent of $\mathcal{F}$. Therefore, if $\ell_n=o(n^{1/3})$, we have
\begin{equation}\label{eq:uc}
\sup_{\mathcal{F}\in\mathscr{F}}\sup_{u\in\mathbb{R}}\Big|\mathbb{P}\big[n^{1/2}\balpha^\T\bJ_{\mathcal{F}}^{-1/2}[\mathbb{E}\{\nabla_{\btheta_{\mathcal{S}}}{\bg}_{i,\mathcal{F}}(\btheta_0)\}]^\T{\bV}_{\mathcal{F}}^{-1}(\btheta_0)\bar{\bg}_{\mathcal{F}}(\btheta_0)\leq u\big]-\Phi(u)\Big|\rightarrow0.
\end{equation}
Write $\Psi_{\mathcal{F}}=\balpha^\T\bJ_{\mathcal{F}}^{-1/2}[\mathbb{E}\{\nabla_{\btheta_{\mathcal{S}}}{\bg}_{i,\mathcal{F}}(\btheta_0)\}]^\T{\bV}_{\mathcal{F}}^{-1}(\btheta_0)\bar{\bg}_{\mathcal{F}}(\btheta_0)$ and $\widehat{\Psi}_{\mathcal{F}}=\balpha^\T\widehat{\bJ}_{\mathcal{F}}^{-1/2}\{\nabla_{\btheta_{\mathcal{S}}}\bar{\bg}_{\mathcal{F}}(\widehat{\btheta}_n)\}^\T\widehat{\bV}_{\mathcal{F}}^{-1}(\widehat{\btheta}_n)\bar{\bg}_{\mathcal{F}}(\btheta_0)$. By Lemmas \ref{la:1} and \ref{la:2}, noting $\sup_{\mathcal{F}\in\mathscr{F}}|\bar{\bg}_{\mathcal{F}}(\btheta_0)|_2=n^{-1/2}\ell_n^{1/2}\varrho_n^{1/4}(\log r)^{1/2}$, we have
\[
\begin{split}
\sup_{\mathcal{F}\in\mathscr{F}}|n^{1/2}(\widehat{\Psi}_{\mathcal{F}}-\Psi_{\mathcal{F}})|=&~O_p\{s\ell_n\omega_n^{1/2}b_n^{1/(2\beta)}\varrho_n^{1/4}(\log r)^{1/2}\}+O_p(\ell_n^{3/2}n^{-1/2}\varrho_n^{3/4}\log r)\\
&+O_p\{s^{3/2}\ell_n\varpi_n^{1/2}b_n^{1/(2\beta)}\varrho_n^{1/4}(\log r)^{1/2}\}+O_p(n^{-1/2}s^{1/2}\ell_n\omega_n^{1/2}\varrho_n^{1/4}\log r)\\
=&~o_p(1).
\end{split}
\]
Hence, for any $u\in\mathbb{R}$, (\ref{eq:uc}) leads to the result. $\hfill\Box$

\end{document}